\documentclass[11pt]{article}
\usepackage[T1]{fontenc}
\usepackage{lmodern,amsmath,amssymb,amsthm,mathtools}
\usepackage[margin=1.08in]{geometry}
\usepackage{enumitem,microtype}
\usepackage[hidelinks]{hyperref}
\numberwithin{equation}{section}
\newtheorem{theorem}{Theorem}[section]
\newtheorem{proposition}[theorem]{Proposition}
\newtheorem{lemma}[theorem]{Lemma}
\theoremstyle{remark}

\newcommand{\R}{\mathbb R}
\newcommand{\Sph}{\mathbb S}
\newcommand{\eps}{\varepsilon}
\newcommand{\dd}{\,\mathrm d}
\newcommand{\tr}{\operatorname{tr}}
\newcommand{\diver}{\operatorname{div}}
\newcommand{\supp}{\operatorname{supp}}
\newcommand{\dist}{\operatorname{dist}}
\newcommand{\Span}{\operatorname{span}}
\newcommand{\diag}{\operatorname{diag}}
\newcommand{\norm}[1]{\lVert#1\rVert}
\setlist[enumerate,1]{label=\textup{(\roman*)},leftmargin=2.3em}
\allowdisplaybreaks[2]
\title{Isolated singularities of solutions to the Yamabe equation
with boundary in higher dimensions II}
\author{Yuxuan Liao\thanks{School of Mathematical Sciences,
Beijing Normal University.
\href{mailto:yxliao@mail.bnu.edu.cn}{\texttt{yxliao@mail.bnu.edu.cn}}}}
\date{November 14, 2026}
\begin{document}
\maketitle
\begin{abstract}
This is the second of two papers on isolated boundary singularities of
the scalar-flat boundary Yamabe equation. We construct positive
solutions that violate the critical upper bound, proving that the
dimension ranges in Part I are sharp. For general smooth metrics,
the conformal metric can be incomplete near the singularity in every
dimension $n\ge18$, and complete in every dimension $n\ge19$.
For Euclidean domains with smooth boundary, both types exist for
$n\ge23$. We also construct both types for metrics with totally geodesic
boundary in dimensions $n\ge35$.


\end{abstract}

\section{Introduction}
Let $g$ be a smooth Riemannian metric on the closed half-ball
$\overline{\mathcal B_1^+}\subset\R^n$, $n\ge3$. We study positive
solutions of the scalar-flat boundary Yamabe equation
\begin{equation}\label{eq:main}
 \begin{cases}
 -L_gu=0&\text{in }\mathcal B_1^+,\\
 B_gu=(n-2)u^{\frac{n}{n-2}}&\text{on }B_1\setminus\{0\}.
 \end{cases}
\end{equation}
Here $\mathcal B_r^+=\{x_n>0,\ |x|<r\}$ and
$B_r=\{x_n=0,\ |x'|<r\}$. The operators
\[
 L_g=\Delta_g-\frac{n-2}{4(n-1)}R_g,\qquad
 B_g=\partial_{\nu_g}+\frac{n-2}{2}h_g
\]
are the conformal Laplacian and conformal boundary operator.
Here $\Delta_g$ is the Laplace--Beltrami operator and $R_g$ is the scalar
curvature. The normal $\nu_g$ is outward, and $h_g$ is the average of the
principal curvatures, positive on a Euclidean sphere. Solutions are assumed smooth
up to the boundary away from zero.

The scalar-flat boundary Yamabe problem seeks a metric of zero scalar
curvature and constant boundary mean curvature in a prescribed conformal
class on a compact manifold with boundary \cite{Escobar,Almaraz}.
Equation~\eqref{eq:main} is the local form of this problem in the case
of positive boundary mean curvature. Allowing the conformal factor to
have isolated boundary singularities leads to the question of its growth
and asymptotic behavior near these points.

In Part~I~\cite[Theorem 1.1]{CompanionI}, we proved the critical upper bound
\[
 u(x)\le C|x|^{-\frac{n-2}{2}}\quad\text{near }0
\]
for general smooth metrics in dimensions $3\le n\le17$ and for
Euclidean domains with smooth boundary in dimensions $3\le n\le22$.
For general metrics, the range extends to $n=18$ when the conformal
metric is complete near the singularity.
The present paper establishes the sharpness of these ranges.
We say that $u^{4/(n-2)}g$ is complete near zero if every locally
rectifiable curve tending to zero has infinite length in this metric.

\begin{theorem}\label{thm:sharpness}
In each of the following cases, there exist a fixed smooth local
metric $g$ and a positive solution of \eqref{eq:main}, with an
isolated boundary singularity at zero, such that
\begin{equation*}
 \limsup_{\substack{x\to0\\x\in B_1}}
       |x|^{\frac{n-2}{2}} u(x)=+\infty.
\end{equation*}
\begin{enumerate}
\item For general smooth metrics, incomplete examples exist in every
dimension $n\ge18$, and complete examples exist in every dimension $n\ge19$.
\item For smooth Euclidean domains, both incomplete and complete examples
exist in every dimension $n\ge23$.
\end{enumerate}
In addition, both incomplete and complete examples exist in every
dimension $n\ge35$ for smooth metrics with totally geodesic boundary,
$\mathrm{II}_g\equiv0$.
All complete examples satisfy $u(x)\ge c |x|^{-{\frac{n-2}{2}}}$ near zero,
while all incomplete examples satisfy
$\liminf_{x\to0,\,x\in B_1}|x|^{\frac{n-2}{2}}u(x)=0$.
\end{theorem}

For flat Euclidean half-balls, Caffarelli et~al.~\cite{CJSX} proved
the critical upper bound in every dimension. Liao--Ma~\cite{LM}
proved the upper bound for smooth metrics in dimensions three and four.
For the interior Yamabe equation, Han et~al.~\cite{HJWXZ1,HJWXZ2}
proved the upper bound through dimension twenty-nine and constructed
counterexamples whose conformal metrics are complete near the
singularity in every dimension at least thirty.

The constructions use localized geometric perturbations and
finite-dimensional reduction, following Brendle~\cite{Brendle},
Brendle--Marques~\cite{BrendleMarques}, and the boundary construction
of Almaraz~\cite{AlmarazBlowup}.
We adapt the successive bubble insertion of Han et~al.~\cite{HJWXZ2},
using the flat singular models and linearized theory of
Jin--Liao--Xiong~\cite{JLXBoundary} and the local estimates and
geometric quadratic forms from Part~I~\cite{CompanionI}.
Starting from a smooth positive solution or the stationary singular
solution, we insert boundary bubbles at points tending to zero.
The perturbations are chosen so that the metrics converge in $C^\infty$
across zero and the solutions retain the previously inserted peaks.
A lower comparison with the stationary singular solution ensures
completeness. Dimension eighteen uses a separate
simultaneous insertion construction.

Section~\ref{sec:preliminaries-II} recalls the bubble, its linear
correction, the stationary singular model and the geometric quadratic
forms. Section~\ref{subsec:local-framework}
sets up local insertion and singular limits.
Section~\ref{subsec:successive-examples} treats the three geometric
classes, and Section~\ref{subsec:eighteen-examples} gives the
eighteenth-dimensional construction. The appendices provide local
analytic estimates, integration formulas, metric variations and
polynomial profile calculations.

\section{Preliminaries}\label{sec:preliminaries-II}
We use boundary Fermi coordinates, with $g_{nn}=1$ and $g_{in}=0$
for $i<n$. Write $\Gamma=\partial\R^n_+$,
$\Sigma_r=\{x_n>0,\ |x|=r\}$, and $\mathcal B_r=\{|x|<r\}$.
Tangential indices range from $1$ to $n-1$. Positive constants $c,C$
may change from line to line, and $0<\alpha<1$ is fixed in H\"older norms.
The estimates recalled below are established in \cite{CompanionI}.

Set $U(x)=\bigl(|x'|^2+(1+x_n)^2\bigr)^{-(n-2)/2}$ and $r=|x|$.
Its translates and dilates are
\begin{equation*}
 U_{\lambda,\xi}({x'},{x_n})=
 \left(\frac{\lambda}{|{x'}-\xi|^2+({x_n}+\lambda)^2}\right)^{\frac{n-2}{2}}.
\end{equation*}
For every $\lambda>0$ and $\xi\in\R^{n-1}$,
\[
 \begin{cases}
 -\Delta U_{\lambda,\xi}=0&\text{in }\R^n_+,\\
 \partial_\nu U_{\lambda,\xi}=(n-2)U_{\lambda,\xi}^{\frac n{n-2}}
                                      &\text{on }\Gamma.
 \end{cases}
\]

For a metric perturbation $\delta+\eta H$, we seek a correction
$U+\eta\psi_H$ satisfying the equation to first order in $\eta$.
For this flat linear problem, volume and boundary integrals are over
$\R^n_+$ and $\Gamma$, unless a domain is specified. Use the energy
\begin{equation}\label{eq:action}
 J_g(v)=\frac12\int(|\nabla v|_g^2+{\frac{n-2}{4(n-1)}}R_gv^2)\dd V_g
   +{\frac{n-2}{4}}\int h_gv^2\dd\sigma_g
   -{\frac{(n-2)^2}{2(n-1)}}\int(v_+)^{\frac{2(n-1)}{n-2}}\dd\sigma_g.
\end{equation}
Here $v_+=\max(v,0)$. Write $Q_g=J_g''(0)$ for the metric bilinear form
and $\mathcal R_g=J_g'$ for the weak residual. For compactly supported smooth tests,
\[
 \mathcal R_g(v)[\zeta]
 =-\int\zeta L_gv\,\dd V_g
  +\int\bigl(B_gv-(n-2)v_+^{n/(n-2)}\bigr)\zeta\,\dd\sigma_g.
\]
For a positive function $a$, denote the pointwise nonlinear remainder by
\[
 N_a(w)=(n-2)\left((a+w)_+^{\frac n{n-2}}-a^{\frac n{n-2}}\right)
                          -na^{\frac2{n-2}}w,
\]
so that
\[
 \bigl(\mathcal R_g(a+w)-\mathcal R_g(a)-D\mathcal R_g(a)w\bigr)[\zeta]
 =-\int N_a(w)\zeta\,\dd\sigma_g.
\]
The correction belongs to the energy space
\[
 \mathcal D^{1,2}(\R^n_+)
 =\overline{C_c^\infty(\overline{\R^n_+})}^{\|\nabla\cdot\|_{L^2}},
 \qquad \|v\|_{\mathcal D^{1,2}}=\|\nabla v\|_{L^2(\R^n_+)}.
\]
A Fermi-preserving variation satisfies $H_{in}=H_{ni}=0$
for $1\le i\le n$. For a compactly supported such variation, let
$J_1(H,f)=\partial_\eta J_{\delta+\eta H}(f)|_{\eta=0}$.
For polynomial $H$, $J_1$ denotes the cutoff limit whenever it exists.
The bubble Hessian is $Q_{\lambda,\xi}=J_\delta''(U_{\lambda,\xi})$, with
\[
 Q_{\lambda,\xi}(v,w)=\int\nabla v\cdot\nabla w
        -n\int_\Gamma U_{\lambda,\xi}^{\frac2{n-2}}vw.
\]
Write $Q=Q_{1,0}$ and $Q(v)=Q(v,v)$. Differentiating the equation for
$U+\eta\psi_H$ at $\eta=0$ gives
$Q(\psi_H,\zeta)=-D_fJ_1(H,U)[\zeta]$.

The bubble parameter derivatives lie in the kernel of $Q$,
\[
 \begin{aligned}
 Z_0&=\left.\partial_\lambda U_{\lambda,0}\right|_{\lambda=1}
       =-\frac{n-2}{2}U-x\cdot\nabla U,\\
 Z_i&=\left.\partial_{\xi_i}U_{1,\xi}\right|_{\xi=0}
       =-\partial_iU\quad\text{for }1\le i<n.
 \end{aligned}
\]
The Cayley map converts $Q$ to a ball Robin form,
\[
 T(x',x_n)=\left(\frac{2x'}{|x'|^2+(1+x_n)^2},
                 \frac{|x'|^2+x_n^2-1}{|x'|^2+(1+x_n)^2}\right).
\]
Since $T^*\delta=4U^{4/(n-2)}\delta$, substituting
$\phi=U(v\circ T)$ yields
\begin{equation}\label{eq:ball-robin-form}
 \begin{aligned}
 Q(\phi)&=2^{2-n}\left(\int_{\mathcal B_1}|\nabla v|^2
                         -\int_{\partial\mathcal B_1}v^2\right),\\
 \|\nabla\phi\|_2^2&=2^{2-n}\left(\int_{\mathcal B_1}|\nabla v|^2
                         +\frac{n-2}{2}\int_{\partial\mathcal B_1}v^2\right).
 \end{aligned}
\end{equation}
For a spherical harmonic $Y_\ell$ of degree $\ell$, its harmonic
extension is $r^\ell Y_\ell$ and its normal derivative on the sphere is
$\ell Y_\ell$. Thus the Steklov eigenvalues are $\ell=0,1,\ldots$,
and the Robin eigenvalues are $\ell-1$. Degree one corresponds to
$Z_0,\ldots,Z_{n-1}$. For the other degrees,
\[
 \frac{|\ell-1|}{\ell+(n-2)/2}\ge\frac2{n+2}\quad\text{for }\ell\ne1.
\]
For $v\in H^1_0(\mathcal B_1)$, $Q(\phi)=\|\nabla\phi\|_2^2=2^{2-n}\int_{\mathcal B_1}|\nabla v|^2$.
Consequently $\ker Q=\Span\{Z_0,\ldots,Z_{n-1}\}$, and $Q$ has a bounded inverse
on their Dirichlet-orthogonal complement, as in the flat case of
\cite[Proposition 2.7]{AlmarazBlowup}.

For compactly supported $H$, \cite[Section~2.1]{CompanionI} proves
$J_1(H,U_{\lambda,\xi})=0$.
Differentiating in $(\lambda,\xi)$ at $(1,0)$ yields the compatibility
conditions
\[
 D_fJ_1(H,U)[Z_\alpha]=0\quad\text{for }0\le\alpha<n.
\]

For a homogeneous Fermi-preserving $H$ of degree $d$, local
differentiation determines the volume and boundary sources
\begin{equation}\label{eq:metric-correction-source}
 \begin{aligned}
 F_H={}&\diver\bigl((\tfrac12\tr H\,I-H)\nabla U\bigr)
 -\frac{n-2}{4(n-1)}
       (\partial_i\partial_jH_{ij}-\Delta\tr H)U,\\
 F_H^{\partial}={}&\frac{n-2}{4(n-1)}(\partial_{x_n}\tr H)U.
 \end{aligned}
\end{equation}
The conormal and boundary-area variations cancel since $H_{in}=0$.
The first variation producing these sources is recorded in
Appendix~\ref{app:metric-variations}.
The linear problem is therefore
\[
 Q(\psi_H,\zeta)=\int F_H\zeta+\int_\Gamma F_H^\partial\zeta.
\]

\begin{lemma}\label{lem:jet-inverse}
Let $H$ be a homogeneous Fermi-preserving tensor of degree
$d\ge1$, and suppose $n>2d+2$. The sources satisfy
\[
 \int_{\R^n_+}F_HZ_\alpha+\int_\Gamma F_H^\partial Z_\alpha=0
 \quad\text{for }0\le\alpha<n.
\]
The correction equation
\begin{equation}\label{eq:jet-correction}
 \begin{cases}
 -\Delta\psi_H=F_H&\text{in }\R^n_+,\\
 (\partial_\nu-nU^{\frac2{n-2}})\psi_H=F_H^{\partial}&\text{on }\Gamma
 \end{cases}
\end{equation}
has a solution in $\mathcal D^{1,2}(\R^n_+)$, unique modulo
$\Span\{Z_0,\ldots,Z_{n-1}\}$. The normalization conditions
\begin{equation}\label{eq:corrector-normalization}
 \psi_H(0)=0,\qquad\nabla_{x'}\psi_H(0)=0
\end{equation}
determine $\psi_H$ uniquely.
For every integer $j\ge0$,
\begin{equation}\label{eq:correction-decay}
 |D^j\psi_H(x)|\le C_j|H|(1+|x|)^{d+2-n-j}.
\end{equation}
For a constant symmetric tangential tensor $A$, the solution of
\eqref{eq:jet-correction} with normalization
\eqref{eq:corrector-normalization} is
\begin{equation}\label{eq:constant-correction}
 \psi_A(x)=\tfrac12A_{ij}x_i\partial_jU(x),
\end{equation}
and satisfies \eqref{eq:correction-decay} with $d=0$.

For $s>0$ and $\xi\in\R^{n-1}$, let $\psi_{H;s,\xi}$ denote the
solution obtained by replacing $U$ with $U_{s,\xi}$ in
\eqref{eq:metric-correction-source} and \eqref{eq:jet-correction},
with normalization
\begin{equation}\label{eq:parameter-corrector-normalization}
 \psi_{H;s,\xi}(\xi,0)=0,\qquad
 \nabla_{x'}\psi_{H;s,\xi}(\xi,0)=0.
\end{equation}
The decay estimates hold
uniformly after any fixed number of derivatives in $(s,\xi)$
on compact subsets of $(0,\infty)\times\R^{n-1}$.
\end{lemma}
\begin{proof}
The source orthogonality, normalized solvability and spatial estimates
are proved in \cite[Lemma~2.2]{CompanionI}. We record the parameter
dependence needed for the reduced energy.
Extend $H\mapsto\psi_H$ by linearity to finite sums of homogeneous
tensors. For $\theta=(s,\xi)$, the polynomial
$H_\theta(y)=H(sy'+\xi,sy_n)$ has degrees zero through $d$.
Substitute $x=(sy'+\xi,sy_n)$ in the correction equation.
Under this substitution, \eqref{eq:corrector-normalization} becomes
\eqref{eq:parameter-corrector-normalization}.
Uniqueness in Lemma~\ref{lem:jet-inverse} implies
\begin{equation}\label{eq:corrector-covariance}
 \psi_{H;s,\xi}(x)=s^{-\frac{n-2}{2}}
 \psi_{H_\theta}
       \left(\frac{x'-\xi}{s},\frac{x_n}{s}\right).
\end{equation}
The homogeneous coefficients of $H_\theta$ depend smoothly on
$\theta$. Their derivatives are bounded by $C_{\Theta,\beta}|H|$
on every compact $\Theta\subset(0,\infty)\times\R^{n-1}$.
Differentiating \eqref{eq:corrector-covariance} and applying the
estimates for degrees $0,\ldots,d$ yields
\[
 |D_\theta^\beta D_x^j\psi_{H;s,\xi}(x)|
 \le C_{\Theta,\beta,j}|H|(1+|x|)^{d+2-n-j}.
\]
The same covariance formula and derivative bounds hold for the
Dirichlet-orthogonal correction $\psi_{H;s,\xi}^\perp$, determined by
\[
 \int_{\R^n_+}\nabla\psi_{H;s,\xi}^\perp\cdot\nabla(\partial_{\theta_\alpha}U_{s,\xi})=0
 \quad\text{for }0\le\alpha<n.
\]
\end{proof}

For a smooth compactly supported Fermi-preserving tensor $h$, take
$g_\eta=\exp(\eta h)$ and define $J_2(h)$ by
\[
 J_{g_\eta}(U)=J_\delta(U)+\eta^2J_2(h)+o(\eta^2).
\]
The first variation vanishes at $U$, so the same coefficient is
obtained along any smooth compactly supported Fermi metric path
with first variation $h$.

Let $\psi_h$ be the correction normalized by
\eqref{eq:corrector-normalization}. Its weak equation is
\[
 Q(\psi_h,\zeta)=-D_fJ_1(h,U)[\zeta].
\]
Testing with $\psi_h$ determines the quadratic contribution of the
correction. Thus
\[
 \begin{aligned}
 J_{g_\eta}(U+\eta\psi_h)
     &=J_\delta(U)+\eta^2\mathcal C(h)+o(\eta^2),\\
 \mathcal C(h)&=J_2(h)-\tfrac12Q(\psi_h).
 \end{aligned}
\]
Symmetric bilinear polarization gives
\begin{equation}\label{eq:reduced-form}
 \mathcal C(h,k)=J_2(h,k)-\frac12Q(\psi_h,\psi_k).
\end{equation}
For trace-free $h$, the direct quadratic coefficient is
\begin{equation}\label{eq:direct-form}
 \begin{aligned}
 J_2(h)={}&\frac14\int_{\R^n_+}|h\nabla U|^2
       -\frac{n-2}{32(n-1)}\int_{\R^n_+}|\partial h|^2U^2\\
 &+\frac{n-2}{16(n-1)}\int_{\R^n_+}|\diver h|^2U^2
 +\frac{n-2}{8(n-1)}\int_{\R^n_+}
                 h_{ik}(\partial_jh_{jk})\partial_i(U^2).
 \end{aligned}
\end{equation}
The metric and curvature expansions leading to this formula are
given in Appendix~\ref{app:metric-variations}.
For polynomial variations of maximal degree $d<(n-2)/2$ and a radial
cutoff $\chi_R(x)=\chi(x/R)$, homogeneity and \eqref{eq:correction-decay}
imply
\[
 \mathcal C(h,k)=\lim_{R\to\infty}\mathcal C(\chi_Rh,\chi_Rk),
 \qquad |\mathcal C(\chi_Rh,\chi_Rk)-\mathcal C(h,k)|=O(R^{2d+2-n}).
\]

The same construction at $U_{\lambda,\xi}$ defines
$\mathcal C(h,k;\lambda,\xi)$. We write $\mathcal C(h;\lambda,\xi)=\mathcal C(h,h;\lambda,\xi)$ for the associated quadratic value.
The notation without parameters refers to $U$.
For homogeneous tensors $H_i,H_j$ of degrees $i,j$,
\eqref{eq:corrector-covariance} and $x=\lambda y$ imply
\[
 \mathcal C(H_i,H_j;\lambda,0)
       =\lambda^{i+j}\mathcal C(H_i,H_j).
\]

Under $r=e^{-t}$, a flat singular solution has the form
$u(r\omega)=r^{-(n-2)/2}v(t,\omega)$, where
$S=\Sph^{n-1}_+$ and $\nu_S$ is its outward normal. Its equation is
\begin{equation}\label{eq:flat-cylinder}
 \begin{cases}
 v_{tt}+\Delta_Sv-\frac{(n-2)^2}{4}v=0&\text{in }\R\times S,\\
 \partial_{\nu_S}v=(n-2)v^{\frac n{n-2}}&\text{on }\R\times\partial S.
 \end{cases}
\end{equation}
By \cite{JLXBoundary}, \eqref{eq:flat-cylinder} admits a unique
stationary positive solution $v(t,\omega)=v_c(\omega)$, which defines
the singular model
\[
 u_c(r\omega)=r^{-\frac{n-2}{2}}v_c(\omega).
\]

For the Euclidean constructions, a boundary graph variation $f$ induces
the Fermi metric variation
\begin{equation}\label{eq:graph-metric-variation}
 (h_f)_{ij}=-2x_n\partial_{ij}f\quad\text{for }i,j<n,
 \qquad (h_f)_{in}=(h_f)_{ni}=0.
\end{equation}
Let $\Lambda$ be the Dirichlet-to-Neumann map of the unit ball
$\mathcal B_1$. On its boundary define
\begin{equation}\label{eq:shape-form}
 Q_\partial(f,g)
 =\langle f,\Lambda(-\Delta_{\Sph^{n-1}}-(n-1))g\rangle_{L^2(\Sph^{n-1})}.
\end{equation}
In stereographic coordinates $z\in\R^{n-1}$ on $\Sph^{n-1}$, write
\[
 \widehat f_{s,\xi}(z)=\frac{2}{1+|z|^2}s^{-1}f(\xi+sz),
 \qquad \widehat f=\widehat f_{1,0}.
\]
The geometric computation in \cite[Section 2.2]{CompanionI} yields
\begin{equation}\label{eq:shape-energy-coefficient}
 \mathcal C(h_f;s,\xi)
 =-\frac{(n-2)^2}{2^{n+1}(n-1)^2}
             Q_\partial(\widehat f_{s,\xi},\widehat f_{s,\xi}).
\end{equation}
On spherical harmonics $Y_\ell$ of degree $\ell$ on $\Sph^{n-1}$,
\[
 Q_\partial(Y_\ell,Y_\ell)=\ell(\ell-1)(\ell+n-1)\|Y_\ell\|_{L^2(\Sph^{n-1})}^2\ge0,
\]
so $Q_\partial\ge0$ and $\ker Q_\partial$ is spanned by spherical
harmonics of degrees $\ell=0,1$.

\section{Local insertion and singular limits}\label{subsec:local-framework}

Each new geometric perturbation is supported in a boundary half-ball
$D\subset\mathcal B_R^+$ with $\overline D\cap(\{0\}\cup\overline{\Sigma_R})=\varnothing$.
The accompanying solution correction is defined on $\mathcal B_R^+$.
We use the following two initial solutions.

\emph{Starting from a smooth solution.}
Take $u_0=U|_{\mathcal B_R^+}$ with metric $\delta$, and fix its
positive trace $u_0|_{\Sigma_R}$. The correction space is
\[
 X=\{v\in H^1(\mathcal B_R^+):v|_{\Sigma_R}=0\},
 \qquad \|v\|_X=\|\nabla v\|_2,\qquad Y=X^*.
\]
The trace bound $\int_{B_R}v^2\le C_nR\|v\|_X^2$ and
$U^{2/(n-2)}\le1$ yield
\[
 \int_{\mathcal B_R^+}|\nabla v|^2
       -n\int_{B_R}U^{2/(n-2)}v^2
 \ge(1-nC_nR)\|v\|_X^2\ge\tfrac12\|v\|_X^2
\]
for small $R$. Thus the initial mixed linearization $D\mathcal R_\delta(u_0):X\to Y$ is invertible.
To pass to $\mathcal B_1^+$, rescale by
$(u(y),g(y))\mapsto(R^{(n-2)/2}u(Ry),g(Ry))$.

\emph{Starting from the cylindrical solution.}
Set $R=1$ and start with metric $\delta$ and
\[
 u_c(r\omega)=r^{-(n-2)/2}v_c(\omega).
\]
Fix $0<\delta<1$ and write $D^c=\mathcal B_R^+\setminus D$.
For a solution difference $v$, define the exterior norms
\[
 \|v\|_{X(D^c)}^2=
 \begin{cases}
 \displaystyle\int_{D^c}|\nabla v|^2\,\dd x,
       &\text{when starting from }u_0,\\[6pt]
 \displaystyle\int_{D^c}r^{-2\delta}
                   (|\nabla v|^2+r^{-2}v^2)\,\dd x,
       &\text{when starting from }u_c.
 \end{cases}
\]
The cylindrical boundary conditions and weighted spaces $X_\delta,Y_\delta$
are specified below.

Assume that at every finite-stage pair $(g,u_*)$, each sufficiently
small boundary half-ball $D=\mathcal B_{8\rho}^+(p)$ centered at
$p\in B_R\setminus\{0\}$ with
\[
 \overline D\subset\overline{\mathcal B_R^+}\setminus
 \bigl(\{0\}\cup\overline{\Sigma_R}\cup\supp(g-\delta)\bigr)
\]
admits pairs $(\hat g,\hat u)$ satisfying the following conditions.
\begin{enumerate}
\item $\hat u>0$ solves \eqref{eq:main} for $\hat g$,
$\supp(\hat g-g)\Subset D$, and the full linearization
$D\mathcal R_{\hat g}(\hat u)$ is an isomorphism in $(X,Y)$ or $(X_\delta,Y_\delta)$.
For a scale $\lambda>0$, there is a local boundary maximum
$\xi_\lambda\in D\cap\Gamma$ with
\[
 |\xi_\lambda-p|\le C\lambda,\qquad
 \hat u(\xi_\lambda)\ge c\lambda^{-(n-2)/2},
\]
where $c,C$ may depend on the fixed stage and $D$.
\item For every $k\ge0$, every finite family of compact sets
$K_1,\ldots,K_q\subset\overline{\mathcal B_R^+}\setminus D$ on which
$u_*$ is smooth, and every $\eta>0$, such pairs exist at arbitrarily
small $\lambda>0$ with
\[
 \|\hat g-g\|_{C^k}
 +\sum_{a=1}^q\|\hat u-u_*\|_{C^k(K_a)}
 +\|\hat u-u_*\|_{X(D^c)}<\eta,
\]
and, when starting from $u_c$,
\[
 \sup_{D^c}\frac{|\hat u-u_*|}{u_c}<\eta.
\]
For Euclidean graph perturbations $g=\Phi^*\delta$ and $\hat g=\widehat\Phi^*\delta$,
$\|\hat g-g\|_{C^k}$ is replaced by $\|\widehat\Phi-\Phi\|_{C^k}$.
\end{enumerate}

\begin{proposition}\label{prop:insertion-limit}
For either initial pair under the insertion hypotheses above, there
exist a smooth metric $g_\infty$ and a positive punctured solution
$u_\infty$ of \eqref{eq:main} such that
\[
 \limsup_{\substack{x\to0\\x\in B_1}}
       |x|^{(n-2)/2}u_\infty(x)=\infty,\qquad
 D^\alpha(g_\infty-\delta)(0)=0\quad\text{for every }\alpha.
\]
If the iteration starts from $u_0$, the conformal metric
$u_\infty^{4/(n-2)}g_\infty$ is incomplete at zero and
\[
 \liminf_{\substack{x\to0\\x\in B_1}}|x|^{(n-2)/2}u_\infty(x)=0.
\]
If it starts from $u_c$, the conformal metric is complete at zero and
$u_\infty(x)\ge c|x|^{-(n-2)/2}$.
\end{proposition}

For the cylindrical construction, we first construct the inverse that
controls the decay of the correction relative to $u_c$.
Write the correction as $v(r\omega)=r^{-(n-2)/2}w(t,\omega)$ with $r=e^{-t}$.
Linearizing \eqref{eq:flat-cylinder} at $v_c$ gives
$\mathcal L_c=-\partial_t^2+A_c$, where
\[
 A_cf=-\Delta_Sf+\frac{(n-2)^2}{4}f\quad\text{in }S,\qquad
 \partial_{\nu_S}f=nv_c^{\frac{2}{n-2}}f\quad\text{on }\partial S.
\]
The boundary trace $v_c|_{\partial S}$ is constant, and $0<c\le v_c\le C$ on
$\overline S$. The bilinear form of $A_c$ is
\[
 Q_c^S(f,h)=\int_S\left(\nabla_Sf\cdot\nabla_Sh+
                         \frac{(n-2)^2}{4}fh\right)
                  -n\int_{\partial S}v_c^{2/(n-2)}fh.
\]
By \cite[Proposition 2.2]{JLXBoundary}, in the normalization \eqref{eq:main},
$0<n v_c^{2/(n-2)}\le n/2$ on $\partial S$, and the stationary equation gives
\[
 Q_c^S(v_c,v_c)=-2\int_{\partial S}v_c^{2(n-1)/(n-2)}<0.
\]
To separate this negative direction, let $\phi_c$ be the principal eigenfunction
from \cite[Section~4]{CompanionI}, satisfying
\[
 A_c\phi_c=-\kappa_c^2\phi_c,\qquad \|\phi_c\|_{L^2(S)}=1.
\]
Set $\Pi_-f=\langle f,\phi_c\rangle_{L^2(S)}\phi_c$ and
$\Pi_+=I-\Pi_-$, and decompose $w(t,\omega)=w_0(t)\phi_c(\omega)+w_+(t,\omega)$.

All following Sobolev norms are on $(0,\infty)\times S$.
The domain and test spaces are
\[
 \begin{aligned}
 X_\delta&=\{w:e^{\delta t}w\in H^1,\ \Pi_+w(0)=0\},\\
 T_{-\delta}&=\{\eta:e^{-\delta t}\eta\in H^1,\ \eta(0)=0\},
 \end{aligned}
\]
with norms
\[
 \|w\|_{X_\delta}=\|e^{\delta t}w\|_{H^1},
 \qquad
 \|\eta\|_{T_{-\delta}}=\|e^{-\delta t}\eta\|_{H^1},
 \qquad Y_\delta=T_{-\delta}^*.
\]
The weak operator $\mathcal L_c:X_\delta\to Y_\delta$ is defined by
\[
 (\mathcal L_cw)[\eta]
 =\int_0^\infty
       \left(\int_S w_t\eta_t+Q_c^S(w,\eta)\right)\dd t.
\]
Since $\eta(0)=0$, integration by parts in $t$ produces no boundary term at $t=0$.
On $\operatorname{ran}\Pi_-$, the equation is $-w_0''-\kappa_c^2w_0=f_0$ on $(0,\infty)$
with $e^{\delta t}w_0\in H^1(0,\infty)$, which determines $w_0(0)$, while
$w_+(0)=0$. For $v=r^{-(n-2)/2}w$, change of variables gives
\begin{equation*}
 \|w\|_{X_\delta}^2\asymp
 \int_{\mathcal B_1^+}r^{-2\delta}
                  (|\nabla v|^2+r^{-2}v^2)\,\dd x.
\end{equation*}

\begin{lemma}\label{lem:complete-cylinder-inverse}
For $0<\delta<1$, the weak map
$\mathcal L_c:X_\delta\to Y_\delta$, with test space $T_{-\delta}$,
is an isomorphism.
\end{lemma}
\begin{proof}
We first establish the angular spectral gap
\[
 \operatorname{spec}(A_c)\cap(-\infty,1)=\{-\kappa_c^2\},
 \qquad\dim\ker(A_c+\kappa_c^2)=1.
\]
Let $\lambda_{\ell,j}$ be the Robin eigenvalues of $A_c$ in tangential
spherical-harmonic degree $\ell$. By \cite[Lemma 2.4]{JLXBoundary},
\[
 \lambda_{0,1}=-\kappa_c^2<0\text{ is simple},\qquad
 \lambda_{0,2}\ge \frac{n^2}{4}.
\]
With $\theta\in(0,\pi/2]$ the polar angle from the inward normal, tangential
translation of $u_c$ yields the degree-one radial eigenfunction
\[
 f_1(\theta)=\frac{n-2}{2}\sin\theta\,v_c(\theta)-\cos\theta\,v_c'(\theta)>0
 \quad\text{for }0<\theta\le\frac\pi2,\qquad\lambda_{1,1}=1,
\]
whose positivity follows from \cite[equation (2.40)]{JLXBoundary}.
Since the centrifugal potential satisfies
\[
 \frac{\ell(\ell+n-3)}{\sin^2\theta}\ge\frac{n-2}{\sin^2\theta}
 \quad\text{for }\ell\ge1,
\]
the Rayleigh quotient gives $\lambda_{\ell,1}\ge\lambda_{1,1}=1$ for all $\ell\ge1$.
This proves the angular gap.

On $\operatorname{ran}\Pi_-$, every functional $f_0$ in the dual of
$H^1_{0,-\delta}(0,\infty)=\{\eta_0:e^{-\delta t}\eta_0\in H^1(0,\infty),\ \eta_0(0)=0\}$
admits a Riesz representation
\[
 f_0[\eta_0]=\int_0^\infty\bigl(f(s)\eta_0(s)-h(s)\eta_0'(s)\bigr)\dd s,
 \qquad \|e^{\delta t}f\|_2+\|e^{\delta t}h\|_2\le C\|f_0\|.
\]
For $t\ge0$, the function $\zeta_t(s)=-\kappa_c^{-1}\mathbf1_{\{s\ge t\}}\sin(\kappa_c(s-t))$
satisfies $\zeta_t(0)=0$ and $\zeta_t\in H^1_{0,-\delta}(0,\infty)$.
Setting $w_0(t)=f_0[\zeta_t]$ gives the variation-of-constants formulas
\begin{align*}
 w_0(t)&=-\int_t^\infty\frac{\sin(\kappa_c(s-t))}{\kappa_c}f(s)\,\dd s
                         +\int_t^\infty\cos(\kappa_c(s-t))h(s)\,\dd s,\\
 w_0'(t)&=\int_t^\infty\cos(\kappa_c(s-t))f(s)\,\dd s-h(t)
                         +\kappa_c\int_t^\infty\sin(\kappa_c(s-t))h(s)\,\dd s.
\end{align*}
Since $e^{-\delta\tau}|\sin(\kappa_c\tau)|+e^{-\delta\tau}|\cos(\kappa_c\tau)|\in L^1(0,\infty)$
for every $\delta>0$, Young's convolution inequality yields
\[
 \|e^{\delta t}w_0\|_{H^1(0,\infty)}
 \le C_{\delta,\kappa_c}\bigl(\|e^{\delta t}f\|_2+\|e^{\delta t}h\|_2\bigr)
 \le C_{\delta,\kappa_c}\|f_0\|.
\]
Since $e^{\delta t}\sin(\kappa_ct),e^{\delta t}\cos(\kappa_ct)\notin L^2(0,\infty)$,
the homogeneous equation $-w_0''-\kappa_c^2w_0=0$ has no nonzero solution
with $e^{\delta t}w_0\in H^1(0,\infty)$.

On $\operatorname{ran}\Pi_+$, write $\widehat w=e^{\delta t}w_+$. The conjugated operator is
\[
 -\widehat w_{tt}+2\delta \widehat w_t+(A_c-\delta^2)\widehat w,\qquad \widehat w(0)=0.
\]
Since $\int_0^\infty\int_S\widehat w_t\widehat w\,\dd\sigma_S\dd t=0$ and
$Q_c^S(\phi,\phi)\ge\|\phi\|_{L^2(S)}^2$ for $\phi\in\operatorname{ran}\Pi_+$,
\[
 \int_0^\infty\bigl(\|\widehat w_t\|_{L^2(S)}^2+Q_c^S(\widehat w,\widehat w)-\delta^2\|\widehat w\|_{L^2(S)}^2\bigr)\dd t
 \ge\int_0^\infty\bigl(\|\widehat w_t\|_{L^2(S)}^2+(1-\delta^2)Q_c^S(\widehat w,\widehat w)\bigr)\dd t.
\]
Moreover, the spectral gap and trace inequality on $S$ give
\[
 \|\phi\|_{H^1(S)}^2
 \le C\bigl(Q_c^S(\phi,\phi)+\|\phi\|_{L^2(S)}^2\bigr)
 \le C Q_c^S(\phi,\phi)\quad\text{for }\phi\in\operatorname{ran}\Pi_+.
\]
Lax--Milgram therefore provides the inverse on $\operatorname{ran}\Pi_+$.
Combining both components gives
\[
 \|w\|_{X_\delta}\le C\|\mathcal L_cw\|_{Y_\delta},\qquad
 \mathcal L_c(X_\delta)=Y_\delta.
\qedhere
\]
\end{proof}

Let $D$ be a pulled-back half-ball of radius comparable to $\rho$,
and let $\varphi$ solve
\[
 \begin{cases}
 -L_g\varphi=0&\text{in }D,\\
 B_g\varphi=0&\text{on }\partial D\cap\Gamma,\\
 \varphi=1&\text{on }\partial D\setminus\Gamma.
 \end{cases}
\]
\begin{lemma}\label{lem:core-barrier}
If $\sup_D(|g-I|+\rho|Dg|+\rho^2|D^2g|)$ is sufficiently small,
then $1/2\le\varphi\le2$ on $D$. Let $v>0$ solve \eqref{eq:main} on $D$, with
$m=\inf_{\partial D\setminus\Gamma}v>0$. Then
\begin{equation}\label{eq:complete-core-barrier}
 v\ge m\varphi\ge \frac{m}{2}\quad\text{on }D.
\end{equation}
\end{lemma}
\begin{proof}
Rescaling \eqref{eq:positive-gauge-bound} to $D$ via
\eqref{eq:small-mixed-coercivity} gives
\[
 \|\varphi-1\|_{L^\infty(D)}
 \le C(\rho^2\|R_g\|_{L^\infty(D)}+\rho\|h_g\|_{L^\infty(\partial D\cap\Gamma)})\le\frac12,
\]
so $1/2\le\varphi\le2$.
Set $f=v-m\varphi$. Then $f\ge0$ on $\partial D\setminus\Gamma$ and
\[
 \begin{cases}
 -L_gf=0&\text{in }D,\\
 B_gf=(n-2)v^{n/(n-2)}\ge0&\text{on }\partial D\cap\Gamma.
 \end{cases}
\]
Thus $f_-=\max(-f,0)\ge0$ vanishes on $\partial D\setminus\Gamma$ and satisfies
$ff_-=-f_-^2$ and $\langle\nabla f,\nabla f_-\rangle_g=-|\nabla f_-|_g^2$ almost everywhere.
Testing with $f_-$ and using \eqref{eq:small-mixed-coercivity} yields
\[
 c\|f_-\|_{H^1(D)}^2
 \le Q_g(f_-,f_-)=-Q_g(f,f_-)
 =-(n-2)\int_{\partial D\cap\Gamma}v^{n/(n-2)}f_-\dd\sigma_g\le0.
\]
Hence $f_-=0$, which proves \eqref{eq:complete-core-barrier}.
\end{proof}

\begin{proof}[Proof of Proposition~\ref{prop:insertion-limit}]
Set $p_j=r_je_1$, $r_j=R4^{-j-3}$, $\rho_j=r_j/64$, and
$D_j=\mathcal B_{8\rho_j}^+(p_j)$. Define the rays
\[
 \gamma_n=\{te_n:0<t\le R/2\},\qquad
 \gamma_-=\{-te_1:0<t\le R/2\}.
\]
The closures $\overline{D_j}$ are pairwise disjoint and disjoint from
$\{0\}\cup\gamma_n\cup\gamma_-\cup\overline{\Sigma_R}$.
At stage $j$, choose the insertion pair $(g_j,u_j)$ on $D_j$ at scale
$\lambda_j>0$ with local boundary maximum $\xi_j\in D_j\cap\Gamma$ such that
\[
 \|g_j-g_{j-1}\|_{C^{j+3}}\le2^{-j}r_j^j,\qquad
 \frac{\lambda_j}{\rho_j}<2^{-j},\qquad
 u_j(\xi_j)\ge4^jr_j^{-(n-2)/2},\qquad |\xi_j-p_j|=o(r_j).
\]
For Euclidean graph deformations $g_j=\Phi_j^*\delta$, impose the $C^{j+3}$ bound
on $\Phi_j-\Phi_{j-1}$ with $\sum_{j=1}^\infty\|\Phi_j-\Phi_{j-1}\|_{C^1}$ small.
Set
\[
 K_j=\bigl\{x\in\overline{\mathcal B_{R/2}^+}:|x|\ge2r_j\bigr\}
     \cup\{\xi_1,\ldots,\xi_{j-1}\}.
\]

When starting from $u_0$, finite energy and \eqref{eq:critical-robin-test}
make $u_j$ smooth across $0$. Put $\widetilde K_j=K_j\cup\overline{\gamma_n}\cup\overline{\gamma_-}$
and choose $\eps_j>0$ with $\sum_{j=1}^\infty\eps_j<\frac12\min_{\overline{\mathcal B_{R/2}^+}}u_0$
such that
\[
 \|u_j-u_{j-1}\|_{C^j(\widetilde K_j)}+\|u_j-u_{j-1}\|_{X(D_j^c)}\le\eps_j.
\]
For every $k\ge0$,
\[
 \sum_{j\ge k}\|g_j-g_{j-1}\|_{C^k}
       \le\sum_{j\ge k}2^{-j}r_j^j<\infty,
\]
and $\supp(g_j-g_0)\subset\{|x|\ge r_j/2\}$, so $g_j\to g_\infty$ in $C^\infty(\overline{\mathcal B_R^+})$
with $D^\alpha(g_\infty-g_0)(0)=0$ for all $\alpha$.
For every compact $K\subset\overline{\mathcal B_{R/2}^+}\setminus\{0\}$ and $l>j$ with $2r_j<\dist(0,K)$,
\[
 \|u_l-u_j\|_{C^k(K)}\le\sum_{i=j+1}^l\eps_i
       \longrightarrow0\quad\text{as }j\to\infty.
\]
Thus $u_j\to u_\infty$ in $C^\infty_{\mathrm{loc}}(\overline{\mathcal B_{R/2}^+}\setminus\{0\})$.
On $\gamma_n\cup\gamma_-$, $u_\infty\ge u_0-\sum_{j=1}^\infty\eps_j>0$, so the strong maximum
and boundary-point principles give $u_\infty>0$. At the inserted peaks,
\[
 u_\infty(\xi_j)\ge u_j(\xi_j)-\sum_{i>j}\eps_i
 \ge 4^jr_j^{-(n-2)/2}-\sum_{i>j}\eps_i,
\]
which implies $\limsup_{x\to0,\,x\in B_1}|x|^{(n-2)/2}u_\infty(x)=\infty$.
Since $\sup_{\gamma_-\cup\gamma_n}u_\infty<\infty$, we have on $\gamma_-$ and $\gamma_n$,
\[
 \lim_{t\to0^+}t^{\frac{n-2}{2}}u_\infty(-te_1)=0,\qquad
 \int_{\gamma_n}u_\infty^{\frac{2}{n-2}}\dd s_{g_\infty}
 \le C\int_0^{R/2}\dd t<\infty,
\]
so $u_\infty^{4/(n-2)}g_\infty$ is incomplete at zero.

When starting from $u_c$, Lemma~\ref{lem:complete-cylinder-inverse} provides the initial inverse.
Impose at stage $j$
\[
 \|u_j-u_{j-1}\|_{C^j(K_j)}+\|u_j-u_{j-1}\|_{X(D_j^c)}
 +\sup_{D_j^c}\frac{|u_j-u_{j-1}|}{u_c}\le\eps_j.
\]
First require $\sum_{j=1}^\infty\eps_j<1/2$. Since $\partial D_j\setminus\Gamma\subset\bigcap_{i\le j}D_i^c$,
\[
 u_j\ge\left(1-\sum_{i=1}^j\eps_i\right)u_c\ge\frac12u_c
 \quad\text{on }\partial D_j\setminus\Gamma.
\]
Since $|x|\asymp r_j$ on $D_j$ and $0<c\le v_c\le C$,
Lemma~\ref{lem:core-barrier} yields $u_j\ge c_0u_c$ on $D_j$ with $c_0>0$ independent of $j$.
Taking $\sum_{j=1}^\infty\eps_j<\frac12\min(1,c_0)$ gives, for all $N\ge j$,
\[
 u_N\ge\left(c_0-\sum_{i=j+1}^N\eps_i\right)u_c
                    \ge\frac{c_0}{2}u_c
 \quad\text{on }D_j,
\]
while $u_N\ge\frac12u_c$ on $\mathcal B_R^+\setminus\bigcup_{j=1}^N D_j$.
Hence
\[
 u_N\ge\frac12\min(1,c_0)u_c\quad\text{on }\mathcal B_R^+\setminus\{0\}.
\]
Passing to the limit $N\to\infty$ gives $g_\infty,u_\infty$ with
$u_\infty\ge cu_c\ge c'|x|^{-(n-2)/2}$. Every curve $\gamma$ tending to zero satisfies
\[
 L_{u_\infty^{4/(n-2)}g_\infty}(\gamma)
                 \ge c\int_\gamma\frac{|\dd r|}{r}=\infty,
\]
so $u_\infty^{4/(n-2)}g_\infty$ is complete at zero.
For Euclidean domains, $D\Phi_\infty(0)=I$, so these properties transfer to
$\Phi_\infty(\mathcal B_R^+)$. Rescaling to a smaller half-ball removes $\Sigma_R$.
\end{proof}

At a fixed stage, we construct the pair required in
Proposition~\ref{prop:insertion-limit} by the local projection of
\cite[Propositions 2.7--2.9]{AlmarazBlowup}.
At a finite stage of the cylindrical construction, write
$u_*=r^{-\frac{n-2}{2}}(v_c+w)$ with $w\in X_\delta$, where
$\supp(g-\delta)\Subset\overline{\mathcal B_1^+}\setminus(\{0\}\cup\overline{\Sigma_1})$.
Subtracting the exact equation for $u_c$ leaves the nonlinear boundary remainder
$-N_{v_c}(w)$. The critical trace inequality and
\[
 e^{\delta t}|w|^{\frac{n}{n-2}}
       =e^{-\frac{2}{n-2}\delta t}|e^{\delta t}w|^{\frac{n}{n-2}}
\]
imply, for $\|w\|_{X_\delta}\le1$,
\[
 \begin{aligned}
 \|N_{v_c}(w)\|_{Y_\delta}
 &\le C\|w\|_{X_\delta}^{\min\{2,n/(n-2)\}},\\
 \|DN_{v_c}(w)\|_{X_\delta\to Y_\delta}
 &\le C\|w\|_{X_\delta}^{\min\{1,2/(n-2)\}}.
 \end{aligned}
\]
Thus $w\mapsto N_{v_c}(w)$ is $C^1$ on $X_\delta$ with $DN_{v_c}(0)=0$.
The metric perturbation satisfies
\[
 \|(\mathcal L_{g,u_c}-\mathcal L_{\delta,u_c})\|_{X_\delta\to Y_\delta}
 \le C\sup_{\mathcal B_1^+}\bigl(|g-I|+r|\partial g|+r^2|\partial^2g|\bigr).
\]
On unit cylinder strips $[T,T+1]\times\overline S$ with large $T$,
absorption of the critical boundary term in \eqref{eq:critical-robin-test}
and \eqref{eq:weighted-strip-regularity} give
\[
 \|w\|_{C^k([T,T+1]\times\overline S)}
 \le C_ke^{-\delta T}\|w\|_{X_\delta},
 \qquad v_c+w\ge\tfrac12\min_{\overline S} v_c.
\]
Since the cylindrical boundary derivative carries the factor $r=e^{-t}$,
the linearized boundary coefficient difference satisfies
\[
 rn\bigl(u_*^{2/(n-2)}-u_c^{2/(n-2)}\bigr)(r\omega)
 =n\bigl((v_c+w)^{2/(n-2)}-v_c^{2/(n-2)}\bigr)
 =O(e^{-\delta t})\qquad\text{as }t\to\infty.
\]
Local trace compactness on $[0,T]\times\partial S$ and the bound
$O(e^{-\delta T})$ on $[T,\infty)\times\partial S$ make multiplication by
this difference compact from $X_\delta$ to $Y_\delta$.
By Lemma~\ref{lem:complete-cylinder-inverse}, the small metric perturbation
and compact boundary-potential perturbation preserve Fredholm index zero at
every finite stage.

In what follows, $(X,Y)$ denotes either the unweighted spaces on $\mathcal B_R^+$
or the weighted spaces $(X_\delta,Y_\delta)$.
For $W>0$, let $\mathcal L_{g,W}=D\mathcal R_g(W)$, with bilinear form
\[
 Q_{g,W}(v,w)=(\mathcal L_{g,W}v)[w]
 =Q_g(v,w)-n\int W^{\frac2{n-2}}vw\,\dd\sigma_g.
\]
In particular, $Q_{\delta,U_{s,\xi}}=Q_{s,\xi}$.

Fix a boundary Fermi chart $D$ centered at $p\ne0$, with
$\overline D\cap(\{0\}\cup\overline{\Sigma_R})=\varnothing$ and $g(p)=I$.
Let $\Theta\Subset(0,\infty)\times\R^{n-1}$ be compact.
Write $\theta=(s,\xi)$ and $\partial_\alpha=\partial_{\theta_\alpha}$ for $0\le\alpha<n$.
Fix cutoffs $\chi_D,\widetilde\chi_D\in C_c^\infty(D\cup(D\cap\Gamma))$ equal to $1$
near $p$ with $\widetilde\chi_D=1$ on $\supp\chi_D$, and set
\[
 U_{\lambda,\theta}(x)=\chi_D(x)\lambda^{-\frac{n-2}{2}}U_{s,\xi}\Bigl(\frac{x-p}{\lambda}\Bigr).
\]
For a positive solution $u_{*,\lambda}$ of \eqref{eq:main} for
$g_\lambda$, independent of $\theta\in\Theta$, set
\[
 W_\theta=u_{*,\lambda}+U_{\lambda,\theta}>0,\qquad
 Z_\alpha=\partial_\alpha U_{\lambda,\theta},
\]
\[
 \ell_\alpha(\eta)=\int_D\nabla Z_\alpha\cdot\nabla(\widetilde\chi_D\eta)\,\dd x
 =\int_D\nabla Z_\alpha\cdot\nabla\eta\,\dd x.
\]
Write $\ell(v)=(\ell_\alpha(v))_\alpha$,
$G_{\alpha\beta}=\ell_\beta(Z_\alpha)$ and
$\mathcal L_{\lambda,\theta}=\mathcal L_{g_\lambda,W_\theta}$.
Define the augmented operator $\mathcal A_{\lambda,\theta}:X\times\R^n\to Y\times\R^n$ by
\begin{equation}\label{eq:insertion-augmented}
 \mathcal A_{\lambda,\theta}(v,{c})=
 \left(\mathcal L_{\lambda,\theta}v-\sum_\alpha c_\alpha\ell_\alpha,
                             \ell(v)\right).
\end{equation}
For a fixed metric $g$ and solution $u_*$, write
$\mathcal L_*=\mathcal L_{g,u_*}:X\to Y$, and identify
$U_{\lambda,\theta}^{\frac{2}{n-2}}$ with the boundary multiplication operator
$v\mapsto U_{\lambda,\theta}^{\frac{2}{n-2}}v|_{\Gamma}\in Y$.

\begin{lemma}\label{lem:insertion-inverse}
Assume $\mathcal L_*:X\to Y$ is invertible and $g_\lambda\to g$ in $C^\infty$
on the fixed charts as $\lambda\to0$. Suppose, uniformly for $\theta\in\Theta$,
\[
 \mathcal L_{\lambda,\theta}
 =\mathcal L_*-nU_{\lambda,\theta}^{\frac2{n-2}}+o_{X\to Y}(1).
\]
For all sufficiently small $\lambda$, the augmented operator
$\mathcal A_{\lambda,\theta}:X\times\R^n\to Y\times\R^n$ in
\eqref{eq:insertion-augmented} is invertible, and
\[
 \sup_{\theta\in\Theta}
 \|\mathcal A_{\lambda,\theta}^{-1}\|_{Y\times\R^n\to X\times\R^n}\le C,
\]
with $C$ depending only on $g$ and $u_*$.
\end{lemma}
\begin{proof}
For fixed $\lambda>0$, $\supp U_{\lambda,\theta}\Subset D\cup(D\cap\Gamma)$,
so local trace compactness on $D\cap\Gamma$ makes
$v\mapsto nU_{\lambda,\theta}^{2/(n-2)}v$ compact from $X$ to $Y$.
Since the $o_{X\to Y}(1)$ term is absorbed by $\mathcal L_*^{-1}$ and
the finite-rank augmentation
\[
\mathcal A_{\lambda,\theta}
=\begin{pmatrix}
    \mathcal L_{\lambda,\theta}&-\ell\\
    \ell&0
\end{pmatrix}
\]
preserves the Fredholm index, we have
\[
\operatorname{ind}\mathcal A_{\lambda,\theta}
=\operatorname{ind}\mathcal L_*=0.
\]

If the uniform bound failed, there would exist $\lambda_j\to0$,
$\theta_j\to\theta_*$ and $(v_j,c_j)\in X\times\R^n$ with $\|v_j\|_X+|c_j|=1$ such that
\[
 \mathcal L_jv_j-\ell_jc_j=h_j\longrightarrow0\quad\text{in }Y,
 \qquad \ell_j(v_j)\longrightarrow0,
\]
where the subscript $j$ denotes evaluation at $(\lambda_j,\theta_j)$.
Since $\|\mathcal L_jZ_{j,\alpha}\|_Y=o(1)$, testing with $Z_{j,\alpha}$ gives
\[
 (G_jc_j)_\alpha
 =\langle \mathcal L_jZ_{j,\alpha},v_j\rangle-h_j[Z_{j,\alpha}]=o(1).
\]
Since $G_j\to(\int_{\R^n_+}\nabla Z_\alpha\cdot\nabla Z_\beta)>0$, $G_j^{-1}$ is uniformly
bounded and $c_j\to0$.
Passing to a subsequence, $v_j\rightharpoonup v$ in $X$.
For any $K>1$,
\[
 \|U_{\lambda_j,\theta_j}^{2/(n-2)}\|_{L^{n-1}(|x-p|>K\lambda_j)}
                              \le \frac{C}{K}.
\]
For any $\zeta\in C_c^\infty$, splitting the boundary integral at $|x-p|=K\lambda_j$ gives
\[
 \begin{aligned}
 \left|\int_\Gamma U_{\lambda_j,\theta_j}^{\frac{2}{n-2}}v_j\zeta\,\dd x'\right|
 &\le C\|\zeta\|_{L^\infty}(K\lambda_j)^{\frac{n-1}{2}}\|v_j\|_X
     +\|U_{\lambda_j,\theta_j}^{\frac{2}{n-2}}\|_{L^{n-1}(|x-p|>K\lambda_j)}
      \|v_j\|_X\|\zeta\|_X\\
 &\le o_j(1)+CK^{-1}\|\zeta\|_X.
 \end{aligned}
\]
Letting $j\to\infty$ and then $K\to\infty$ yields $\mathcal L_*v=0$, so $v=0$.

At the bubble scale, set
\[
 \widehat v_j(y)=\lambda_j^{\frac{n-2}{2}}(\widetilde\chi_Dv_j)(p+\lambda_jy).
\]
Then $\widehat v_j\rightharpoonup\widehat v$ in $\mathcal D^{1,2}(\R^n_+)$, where
$Q_{\theta_*}(\widehat v,\zeta)=0$ for all $\zeta\in\mathcal D^{1,2}(\R^n_+)$.
For the flat Jacobi fields $Z_\alpha$ at $\theta_*$,
\[
 \left|\int_{|y|>K}\nabla\widehat v_j\cdot\nabla Z_\alpha\,\dd y\right|
 \le C\|\nabla Z_\alpha\|_{L^2(|y|>K)}\longrightarrow0\quad\text{as }K\to\infty.
\]
Combined with $\ell_j(v_j)\to0$, this yields
$\int_{\R^n_+}\nabla\widehat v\cdot\nabla Z_\alpha=0$ for $0\le\alpha<n$,
so the kernel characterization following \eqref{eq:ball-robin-form} implies $\widehat v=0$.

By local trace compactness on $B_K\subset\Gamma$,
$\|\widehat v_j\|_{L^{2(n-1)/n}(B_K)}\to0$ for each fixed $K$.
Dual trace H\"older and the exterior $L^{n-1}$ bound on $|x-p|>K\lambda_j$ give
\[
 \|U_{\lambda_j,\theta_j}^{\frac{2}{n-2}}v_j\|_Y
 \le C\|\widehat v_j\|_{L^{\frac{2(n-1)}{n}}(B_K)}
     +\|U_{\lambda_j,\theta_j}^{\frac{2}{n-2}}\|_{L^{n-1}(|x-p|>K\lambda_j)}\|v_j\|_X
 \le o_j(1)+CK^{-1}\|v_j\|_X.
\]
Since
\[
 \mathcal L_*v_j
 =nU_{\lambda_j,\theta_j}^{\frac{2}{n-2}}v_j+\ell_jc_j+h_j+o_Y(1),
\]
applying $\mathcal L_*^{-1}$ yields
\[
 \|v_j\|_X\le C o_j(1)+CK^{-1}\|v_j\|_X.
\]
Choosing $K>2C$ and $j$ large gives $\|v_j\|_X\to0$, contradicting $\|v_j\|_X+|c_j|=1$.
Since $\operatorname{ind}\mathcal A_{\lambda,\theta}=0$, this proves invertibility and the uniform bound.
\end{proof}

To determine the bubble parameters after projection, let $J_{g_\lambda,D}$
be the energy \eqref{eq:action} restricted to $D$ and $D\cap\Gamma$, and define
\[
 I_\lambda(\theta)=J_{g_\lambda,D}(W_\theta)
                         -J_{g_\lambda,D}(u_{*,\lambda}).
\]
Let $\mu\to0$ be the perturbation amplitude, $\mathcal R_\theta^{(1)}$
its first-order residual on $D$, and $\psi_\theta$ the corresponding
approximate first correction on $D$, with rescaled form
\[
 \widehat\psi_{\lambda,\theta}(y)=\lambda^{\frac{n-2}{2}}\psi_\theta(p+\lambda y).
\]
In rescaled half-space coordinates, write $\mathcal L_\theta$ for the weak operator of
$Q_\theta=Q_{s,\xi}$, $\psi_\theta^0$ for the flat correction solving
$\mathcal L_\theta\psi_\theta^0=-\mathcal R_\theta^{(1)}$, $J_2(\theta)$ for the
direct quadratic energy coefficient, and $\mathcal C(\theta)$ for the completed
coefficient in \eqref{eq:reduced-form}.
For a parameter-dependent correction $w=w_\theta$, define the projected residual
\[
 \mathcal F_{\lambda,\alpha}(\theta)
      =\mathcal R_{g_\lambda}(W_\theta+w_\theta)[Z_\alpha],
 \qquad \mathcal F_\lambda=(\mathcal F_{\lambda,\alpha})_\alpha.
\]

\begin{lemma}\label{lem:local-insertion-reduction}
Let $n\ge5$, and suppose $\|\mathcal A_{\lambda,\theta}^{-1}\|$ and $\|G^{-1}\|$
are uniformly bounded on $\Theta$. Let $\mu=\mu(\lambda)\to0$ as
$\lambda\to0$. Assume that the residual has the expansion
\[
 \mathcal R_{g_\lambda}(W_\theta)
       =\mu\mathcal R_\theta^{(1)}+o_{C^2(\Theta;Y)}(\mu),
\]
and that the approximate correction satisfies
\[
 \ell(\psi_\theta)=0,\qquad
 \mathcal L_{\lambda,\theta}\psi_\theta+\mathcal R_\theta^{(1)}
       =o_{C^1(\Theta;Y)}(1).
\]
Suppose also that
\[
 I_\lambda(\theta)
       =I_{\lambda,0}+\mu^2J_2(\theta)+o_{C^2(\Theta)}(\mu^2),
\]
where $I_{\lambda,0}$ is independent of $\theta$. Finally, assume that
the rescaled corrections satisfy
\[
 \widehat\psi_{\lambda,\theta}\longrightarrow\psi_\theta^0
 \quad\text{in }C^1\bigl(\Theta;H^1_{\mathrm{loc}}(\overline{\R^n_+})\bigr),
\]
with uniformly vanishing exterior Dirichlet integrals,
\[
 \lim_{K\to\infty}\sup_{\lambda,\theta}\sum_{|\beta|\le1}\|\nabla D_\theta^\beta\widehat\psi_{\lambda,\theta}\|_{L^2(|y|>K)}=0.
\]
Then there exist a unique $w=w_\theta\in X$ with $\ell(w)=0$ and $\|w\|_X\le C\mu$,
and multipliers $c=(c_\alpha)$, satisfying
\begin{equation}\label{eq:insertion-correction}
 \mathcal R_{g_\lambda}(W_\theta+w_\theta)=\sum_\alpha c_\alpha\ell_\alpha.
\end{equation}
Moreover, $w_\theta=\mu\psi_\theta+o_{C^1(\Theta;X)}(\mu)$ and
\begin{equation}\label{eq:insertion-force}
 \mathcal F_\lambda(\theta)=\mu^2D\mathcal C(\theta)+o_{C^1(\Theta)}(\mu^2).
\end{equation}
Every nondegenerate critical point of $\mathcal C$ yields an exact solution
with invertible full linearization.
\end{lemma}
\begin{proof}
Fix $\lambda$ and suppress $\theta$ in $W,\psi,\mathcal R^{(1)}$.
Write $\mathcal L=\mathcal L_{\lambda,\theta}$, $\mathcal A=\mathcal A_{\lambda,\theta}$
and $\mathcal R=\mathcal R_{g_\lambda}$, and set $p=n/(n-2)\in(1,2)$.
The nonlinear residual is
\[
 \mathcal N_W(w)[\eta]
 =\bigl(\mathcal R(W+w)-\mathcal R(W)-\mathcal L w\bigr)[\eta]
 =-\int N_W(w)\eta\,\dd\sigma_{g_\lambda}.
\]
The critical trace inequality gives
\begin{align*}
 \|\mathcal N_W(w)\|_Y&\le C\|w\|_X^p,\\
 \|\mathcal N_W(w)-\mathcal N_W(\widetilde w)\|_Y
 &\le C(\|w\|_X+\|\widetilde w\|_X)^{p-1}\|w-\widetilde w\|_X.
\end{align*}
On $\{w\in X:\|w\|_X\le C\mu\}$, the Lipschitz constant is
$O(\mu^{p-1})<1/2$ for small $\mu$, so the contraction mapping theorem solves
\[
 (w,c)=\mathcal A^{-1}
 \bigl(-\mathcal R(W)-\mathcal N_W(w),0\bigr),
 \qquad \|w\|_X\le C\mu.
\]
Subtracting the approximate first correction gives
\[
 \|w-\mu\psi\|_X
 \le C\bigl(\|\mathcal R(W)-\mu\mathcal R^{(1)}\|_Y
       +\mu\|\mathcal L\psi+\mathcal R^{(1)}\|_Y+\|w\|_X^p\bigr)=o(\mu).
\]
Since $|D_\theta W|+|D_\theta Z|\le CW$ and
\[
 W\left|(W+w)_+^{p-1}-W^{p-1}-(p-1)W^{p-2}w\right|
 \le C|w|^p,
\]
differentiating the fixed-point equation in $\theta$ yields
$w=\mu\psi+o_{C^1(\Theta;X)}(\mu)$.

Since $\partial_\alpha W=Z_\alpha$ is supported in $D$, symmetry of
the energy Hessian on these compactly supported variations yields
\[
 \mathcal F_{\lambda,\alpha}
 =\mathcal R(W)[Z_\alpha]
  +\partial_\alpha\mathcal R(W)[w]
  +\mathcal N_W(w)[Z_\alpha].
\]
Because $\supp Z_\alpha\Subset D\cup(D\cap\Gamma)$, the expansion of $I_\lambda$ gives
\[
 \mathcal R(W)[Z_\alpha]=\partial_\alpha I_\lambda
        =\mu^2\partial_\alpha J_2+o_{C^1}(\mu^2),
\]
and the residual expansion together with $w=\mu\psi+o_{C^1}(\mu)$ gives
\[
 \partial_\alpha\mathcal R(W)[w]
 =\mu^2\langle\partial_\alpha\mathcal R^{(1)},\psi\rangle
                                      +o_{C^1}(\mu^2).
\]
For the third term, Taylor's formula on $|w|\le W/2$ and direct bounds on $|w|>W/2$ give
\[
 W\left|(W+w)_+^p-W^p-pW^{p-1}w
                     -\frac{p(p-1)}2W^{p-2}w^2\right|
       \le C|w|^{p+1},
\]
and similarly after one $\theta$-derivative, so trace H\"older yields
\[
 \begin{aligned}
 &\left\|\mathcal N_W(w)[Z_\alpha]
 -\frac12(\partial_\alpha Q_{g_\lambda,W})(w,w)\right\|_{C^1(\Theta)}\\
 &\qquad\le C\sup_{\Theta}\int_{D\cap\Gamma}
       \bigl(|w|^p|D_\theta w|+|w|^{p+1}\bigr)\dd\sigma_{g_\lambda}
 =O(\mu^{p+1})=o(\mu^2).
 \end{aligned}
\]
Rescaled local convergence $\widehat\psi_{\lambda,\theta}\to\psi_\theta^0$ and
the uniform exterior Dirichlet integral bound therefore imply
\[
 \mathcal N_W(w)[Z_\alpha]
 =\frac{\mu^2}{2}(\partial_\alpha Q_\theta)(\psi_\theta^0,\psi_\theta^0)
                              +o_{C^1}(\mu^2).
\]
Differentiating $\mathcal L_\theta\psi_\theta^0=-\mathcal R^{(1)}$ gives
\[
 \partial_\alpha \mathcal C=\partial_\alpha J_2
       +\langle\partial_\alpha\mathcal R^{(1)},\psi_\theta^0\rangle
       +\frac12(\partial_\alpha Q_\theta)(\psi_\theta^0,\psi_\theta^0),
\]
which proves \eqref{eq:insertion-force}.

By \eqref{eq:insertion-force}, every nondegenerate critical point
of $\mathcal C$ produces a nearby zero $\theta_\lambda$ of $\mathcal F_\lambda$.
Since $\mathcal F_\lambda=Gc$ and $G$ is invertible, $c(\theta_\lambda)=0$, giving an exact solution $u=W+w$.
Since
\[
 \|D\mathcal R(W+w)-D\mathcal R(W)\|_{X\to Y}
          \le C\|w\|_X^{p-1}=o(1),
\]
the augmented operator $\mathcal A_u$ remains invertible.
Put $\mathcal L_u=D\mathcal R(u)$. Differentiating \eqref{eq:insertion-correction}
and $\ell(w)=0$ at $c=0$ gives
\[
 \mathcal L_uD_\theta u=\ell D_\theta c,
 \qquad \ell(D_\theta u)=G+O(\mu),
\]
and
\[
 D_\theta c=G^{-1}D_\theta \mathcal F_\lambda
       =G^{-1}\bigl(\mu^2D^2\mathcal C+o(\mu^2)\bigr),
\]
which is invertible by nondegeneracy of $D^2\mathcal C$.
For any $h\in Y$, if $(v,d)=\mathcal A_u^{-1}(h,0)$, then
\[
 \mathcal L_u\bigl(v-D_\theta u(D_\theta c)^{-1}d\bigr)=h.
\]
Since $\operatorname{ind}\mathcal L_u=0$, this right inverse is the two-sided inverse
$\mathcal L_u^{-1}:Y\to X$, with $\|\mathcal L_u^{-1}\|_{Y\to X}=O(\mu^{-2})$.
Finally, strong rescaled energy convergence and \eqref{eq:critical-robin-test}
yield smooth convergence on fixed half-balls and exterior compact sets, while
\eqref{eq:weighted-strip-regularity} gives weighted convergence on the cylindrical end.
\end{proof}

In both cases, positivity follows when $Q_g$ is coercive on energy
functions with zero trace on the outer spherical boundary.
When starting from $u_0$, $u|_{\Sigma_R}=u_0|_{\Sigma_R}>0$.
When starting from $u_c$, the boundary condition $\Pi_+w(0)=0$ gives
\[
 u|_{\Sigma_1}=v_c+w_0(0)\phi_c>0
\]
provided the total correction is sufficiently small. In this case,
\eqref{eq:weighted-strip-regularity} also gives $u>0$ near zero.
Thus $u_-=\max(-u,0)$ has zero outer spherical trace and belongs to
the energy space. In the singular case it vanishes near zero.
Testing the positive-part equation with $u_-$ yields
\begin{equation}\label{eq:insertion-positivity}
 Q_g(u_-,u_-)
 =-(n-2)\int_{B_R}u_+^{\frac{n}{n-2}}u_-\dd\sigma_g=0.
\end{equation}
Coercivity implies $u_-=0$, and the strong maximum and boundary-point
principles then give $u>0$.

\section{Successive insertion for the three geometric classes}\label{subsec:successive-examples}

The geometric condition in Lemma~\ref{lem:local-insertion-reduction}
is a nondegenerate critical point of $\mathcal C$ in
\eqref{eq:insertion-force}.
We first choose metric polynomials for which $\mathcal C(H;s,\xi)$
has a nondegenerate critical point with strictly negative value.

Use rescaled half-space coordinates $y=(y',y_n)$.
Let $A_1=\diag(1,-1,0,\ldots,0)$ and $A_2=\diag(0,0,1,-1,0,\ldots,0)$
on $\R^{n-1}$, and consider
\begin{equation}\label{eq:metric-profile-family}
 H=\sum_{\alpha=1}^2
       \bigl\{y_nP_\alpha(|y'|^2)+y_n^2Q_\alpha(|y'|^2)\bigr\}A_\alpha,
\end{equation}
where
\[
 \deg P_\alpha\le1,\qquad
 \begin{cases}
 \deg Q_\alpha\le1,\quad\deg H\le4&\text{for }n=19,\\
 Q_\alpha\equiv0,\quad\deg H\le3&\text{for }n\ge20.
 \end{cases}
\]

For these trace-free tangential tensors, use the metric path
\begin{equation}\label{eq:profile-reduced-form}
 g_\eta^{-1}=\exp(\eta H),\qquad g_\eta=\exp(-\eta H),\qquad
 \left.\partial_\eta g_\eta\right|_{\eta=0}=-H.
\end{equation}
Let $\psi_{H;s,\xi}^{\perp}$ be the Dirichlet-orthogonal correction
from Lemma~\ref{lem:jet-inverse}, which solves
\begin{equation}\label{eq:inverse-profile-source}
 \begin{cases}
 -\Delta\psi_{-H;s,\xi}^{\perp}=F_{-H;s,\xi}
       &\text{in }\R^n_+,\\
 (\partial_\nu-nU_{s,\xi}^{2/(n-2)})\psi_{-H;s,\xi}^{\perp}=0
       &\text{on }\Gamma.
 \end{cases}
\end{equation}
Since $\psi_{H;s,\xi}^{\perp}-\psi_{H;s,\xi}\in\ker Q_{s,\xi}$, we have
$Q_{s,\xi}(\psi_{H;s,\xi}^{\perp})=Q_{s,\xi}(\psi_{H;s,\xi})$.
Linearity and \eqref{eq:reduced-form} yield
\[
 \psi_{-H;s,\xi}^{\perp}=-\psi_{H;s,\xi}^{\perp},
 \qquad \mathcal C(-H;s,\xi)=\mathcal C(H;s,\xi).
\]
Consequently, the reduced energy coefficient of
\eqref{eq:profile-reduced-form} is $\mathcal C(H;s,\xi)$.

We record a coefficient-selection lemma which will also be used for
the boundary graphs and totally geodesic metrics.
Let $G$ be a positive definite $k\times k$ matrix, $k\ge2$, and
let $0<d_0<\cdots<d_{k-1}$ be integers.
For $c=(c_1,c_2)\in\R^k\oplus\R^k$, with
$c_\alpha=(c_{\alpha0},\ldots,c_{\alpha,k-1})$, define
\[
 E_c(s)=\sum_{\alpha=1}^2\sum_{i,j=0}^{k-1}
              G_{ij}c_{\alpha i}c_{\alpha j}s^{d_i+d_j},
 \qquad q(c)=E_c'(1).
\]
In the metric application, $G$ is the Gram matrix of
$-\mathcal C$ on the profile basis of homogeneous degrees $d_i$, so
\[
 E_c(s)=-\mathcal C(H;s,0),\qquad
 h_i(c)=-\partial_{\xi_i\xi_i}\mathcal C(H;1,0).
\]

\begin{lemma}\label{lem:profile-selection}
Suppose $\inf_{c\ne0}q(c)<0$. Let $h_1,\ldots,h_N$ be quadratic
forms on $\R^k\oplus\R^k$ with $h_j\notin\R q$ for $1\le j\le N$.
There exists $c=(c_1,c_2)$ with $c_{10}\ne0$ and $c_{1,k-1}^2+c_{2,k-1}^2>0$ such that
\[
 q(c)=0,\qquad \partial_{\log s}^2E_c(1)<0,\qquad
 h_j(c)\ne0\quad\text{for }1\le j\le N.
\]
\end{lemma}
\begin{proof}
The matrix of $q$ is $\diag(M,M)$ with $M_{ij}=(d_i+d_j)G_{ij}$.
Since $M_{ii}=2d_iG_{ii}>0$ and $\inf_{x\ne0}x^tMx<0$, each block $M$
has at least one positive and one negative eigenvalue, so
$\operatorname{rank}q\ge4$ and $q$ is an irreducible indefinite quadratic form.
Write
\[
 \partial_{\log s}E_c(s)=s^{2d_0}p_c(s),
 \qquad \deg p_c\le2(d_{k-1}-d_0).
\]
For $c_1=(1,0,\ldots,0)$ and $c_2=(0,\ldots,0,1)$,
\[
 p_c(s)=2d_0G_{00}
       +2d_{k-1}G_{k-1,k-1}s^{2(d_{k-1}-d_0)}
\]
has only simple roots, so $\operatorname{Disc}(p_c)\not\equiv0$.
Since $\{q<0\}$ is open and nonempty, we may choose $c$ with
\[
 q(c)<0,\qquad
 \operatorname{Disc}(p_c)
 (c_{10}^2+c_{20}^2)(c_{1,k-1}^2+c_{2,k-1}^2)\ne0.
\]
Then $p_c(0)>0$, $\lim_{s\to\infty}p_c(s)=\infty$ and $p_c(1)=q(c)<0$, so the first positive root
$s_*\in(0,1)$ of $p_c$ satisfies $p_c'(s_*)<0$.
Replacing $c_{\alpha i}$ by $s_*^{d_i}c_{\alpha i}$ gives
\[
 q(c)=0,\qquad
 Dq(c)[(d_ic_{\alpha i})_{\alpha,i}]
             =\partial_{\log s}^2E_c(1)<0.
\]
Thus $\Sigma=\{c:q(c)=0,\ \partial_{\log s}^2E_c(1)<0\}$ is a nonempty smooth
hypersurface patch.

If a quadratic form $h$ vanishes on a nonempty open subset $U\subset\Sigma$,
choose coordinates in which $q(x)=x_1^2-Q(x')$. Subtracting a multiple of $q$
gives $h(x)\equiv x_1L(x')+R(x')\pmod{q(x)}$, so $x_1L(x')+R(x')=0$ on $U$.
Since the projection of $U$ onto $x'$ is open, squaring yields $R^2=QL^2$ as polynomials.
Since $\operatorname{rank}Q=\operatorname{rank}q-1\ge3$, $Q$ is irreducible, so
$L=R=0$ and $h\in\R q$.
Consequently, the zero sets of $h_1,\ldots,h_N$, $c_{10}^2$ and
$c_{1,k-1}^2+c_{2,k-1}^2$ have empty interior in $\Sigma$, and there exists
$c\in\Sigma$ outside their union.
\end{proof}

\begin{lemma}\label{lem:metric-wells}
For every $n\ge19$, the family \eqref{eq:metric-profile-family}
contains a polynomial $H_n$ such that $(s,\xi)=(1,0)$ is a
nondegenerate critical point of $\mathcal C(H_n;s,\xi)$ and
\[
 \mathcal C(H_n;1,0)<0,
 \qquad \deg H_n=\begin{cases}4,&n=19,\\3,&n\ge20.\end{cases}
\]
\end{lemma}
\begin{proof}
\emph{Step 1. The positive form and the two copies.}
Since $\deg H\le4<(n-2)/2$, \cite[equation~(2.26) and Lemma~2.5]{CompanionI}
imply $\mathcal C(H)\le0$.
By \cite[Lemma~2.6 and equation~(2.31)]{CompanionI}, $\mathcal C(H)=0$ implies
\[
 H=\mathcal L_Y\delta+f\delta,\qquad Y_n=0\quad\text{on }\Gamma
\]
near zero, and hence
\[
 (D\mathrm{II}_\delta[H])^\circ=0,\qquad DW_\delta[H]=0.
\]
For $H$ of the form \eqref{eq:metric-profile-family},
\[
 (D\mathrm{II}_\delta[H])^\circ
 =-\tfrac12\partial_nH|_\Gamma
 =-\frac12\sum_{\alpha=1}^2P_\alpha(|y'|^2)A_\alpha=0
 \quad\Longrightarrow\quad P_1=P_2=0.
\]
At $n=19$, $H=y_n^2\sum_{\alpha=1}^2 Q_\alpha(|y'|^2)A_\alpha$ then satisfies
\[
 DW_\delta[H]_{ninj}\big|_\Gamma
 =-\frac{n-3}{2(n-2)}\partial_{nn}H_{ij}\big|_\Gamma
 =-\frac{n-3}{n-2}
       \sum_{\alpha=1}^2Q_\alpha(|y'|^2)(A_\alpha)_{ij}=0
 \quad\Longrightarrow\quad Q_1=Q_2=0.
\]
Thus $\mathcal C(H)=0\implies H=0$, so $-\mathcal C>0$ on the coefficient space.

Under $(y_1,y_2)\mapsto(y_2,y_1)$, $(A_1,A_2)\mapsto(-A_1,A_2)$, while under
$(y_1,y_2)\leftrightarrow(y_3,y_4)$, $A_1\leftrightarrow A_2$.
Invariance of $\mathcal C$ at $(s,\xi)=(1,0)$ gives, for radial polynomials $P,Q$,
\[
 \mathcal C(PA_1,QA_2)=0,\qquad
 \mathcal C(PA_1,QA_1)=\mathcal C(PA_2,QA_2).
\]
In the ordered scalar bases
\[
 \renewcommand{\arraystretch}{1.3}
 \begin{array}{c|c|c}
 n&\text{basis}&\text{degrees }d\\ \hline
 19&(y_n,y_n^2,y_n|y'|^2,y_n^2|y'|^2)&(1,2,3,4)\\
 n\ge20&(y_n,y_n|y'|^2)&(1,3)
 \end{array}
\]
let $G>0$ be the Gram matrix of $-\mathcal C$ on the $A_1$ copy.
For $c=(c_1,c_2)$, homogeneity gives
\[
 -\mathcal C(H;s,0)=E_c(s),\qquad q(c)=-\partial_s\mathcal C(H;1,0).
\]

\emph{Step 2. A negative value of the scale derivative.}
At $n=19$, taking $c_1=(3,-3,-1,1)$ and $c_2=0$ gives
$H=y_n(1-y_n)(3-|y'|^2)A_1$, for which \eqref{eq:nineteen-negative-dilation} yields
$q(c)<-1/6<0$.
For $n\ge20$, \eqref{eq:cubic-metric-dilation} gives $4G_{01}^2-3G_{00}G_{11}>0$, so
taking $c_1=(2G_{01}/G_{00},-1)$ and $c_2=0$ yields
\[
 q(c)=-\frac{2(4G_{01}^2-3G_{00}G_{11})}{G_{00}}<0.
\]

\emph{Step 3. The center derivatives and coefficient selection.}
For $\theta=(s,\xi)$ and $H_\theta(y)=H(sy'+\xi,sy_n)$,
\eqref{eq:corrector-covariance} and \eqref{eq:reduced-form} give
$\mathcal C(H;s,\xi)=\mathcal C(H_\theta)$, and Lemma~\ref{lem:jet-inverse} yields
\begin{equation}\label{eq:metric-profile-center-derivatives}
 \partial_{\xi_i\xi_j}\mathcal C(H;\theta)
 =2\mathcal C(\partial_{\xi_i}H_\theta,
                         \partial_{\xi_j}H_\theta)
       +2\mathcal C(H_\theta,
                         \partial_{\xi_i\xi_j}H_\theta).
\end{equation}
Since $H$ is invariant under each coordinate reflection $y_i\mapsto-y_i$,
$\mathcal C(H;s,\xi)$ is even in each $\xi_i$, so
\[
 \nabla_\xi\mathcal C(H;s,0)=0,\qquad
 \partial_{s\xi_i}\mathcal C(H;s,0)=0,\qquad
 \partial_{\xi_i\xi_j}\mathcal C(H;s,0)=0\quad\text{for }i\ne j.
\]
Set $h_i(c)=-\partial_{\xi_i\xi_i}\mathcal C(H;1,0)$.
For $H=y_nA_1$, $\partial_{\xi_i}H_\theta=0$, so
\[
 h_i(c)=0,\qquad q(c)=2G_{00}>0.
\]
For $H=y_n|y'|^2A_1$, \eqref{eq:metric-profile-center-derivatives},
$\mathcal C\le0$ and \eqref{eq:metric-positive-coupling} yield
\[
 h_i(c)
 =-2\mathcal C(2y_ny_iA_1)
                   -4\mathcal C(y_n|y'|^2A_1,y_nA_1)
 \ge-4\mathcal C(y_n|y'|^2A_1,y_nA_1)>0.
\]
Hence $h_i\notin\R q$ for $1\le i\le n-1$.
Lemma~\ref{lem:profile-selection} therefore provides $c\ne0$ with
$E_c'(1)=0$, $E_c''(1)<0$ and $h_i(c)\ne0$ for all $i$.
The resulting polynomial $H_n$ satisfies $D_{(s,\xi)}\mathcal C(H_n;1,0)=0$,
\[
 D^2_{(s,\xi)}\mathcal C(H_n;1,0)
 =\diag\bigl(-E_c''(1),-h_1(c),\ldots,-h_{n-1}(c)\bigr),
 \qquad \det D^2_{(s,\xi)}\mathcal C(H_n;1,0)\ne0,
\]
and $\mathcal C(H_n;1,0)=-E_c(1)<0$.
\end{proof}

Fix a tangential trace-free polynomial $H$ of degree $d_*\ge1$
with $n>4d_*+2$ such that $\theta_*$ is a nondegenerate critical point
of $\mathcal C(H;\cdot)$ with $\mathcal C(H;\theta_*)<0$.
Choose a small closed box $\Theta$ around $\theta_*$ containing no other
critical point, so that
\begin{equation}\label{eq:profile-margin}
 \min_{\theta\in\partial\Theta}|D_\theta\mathcal C(H;\theta)|=m_*>0,
 \qquad \deg(D_\theta\mathcal C(H;\cdot),\Theta,0)
 =\operatorname{sgn}\det D_\theta^2\mathcal C(H;\theta_*)\in\{-1,1\}.
\end{equation}

For either choice of $(X,Y)$ in Section~\ref{subsec:local-framework},
let $g$ be a smooth metric with $Q_g$ coercive on energy functions
with zero outer spherical trace, and let $u_*>0$
solve \eqref{eq:main} with $\mathcal L_{g,u_*}:X\xrightarrow{\cong}Y$.
Let $\mathcal U=\mathcal B_{8\rho}^+(p)$ be a boundary half-ball centered at
$p\ne0$ with $\overline{\mathcal U}\cap(\{0\}\cup\overline{\Sigma_R}\cup\supp(g-\delta))=\varnothing$.

\begin{proposition}\label{prop:complete-insertion}
For all sufficiently small $\lambda>0$, there exist a smooth metric
$g_\lambda$ with $\supp(g_\lambda-g)\Subset\mathcal U$ and a positive
solution $u_\lambda$ of \eqref{eq:main} such that
$\mathcal L_{g_\lambda,u_\lambda}:X\to Y$ is invertible, with a local
boundary maximum $\xi_\lambda\in\mathcal U\cap\Gamma$ satisfying
\[
 |\xi_\lambda-p|\le C\lambda,\qquad
 u_\lambda(\xi_\lambda)\ge c\lambda^{-\frac{n-2}{2}}.
\]
For every fixed $k\ge0$ and compact
$K\subset\overline{\mathcal B_R^+}\setminus
(\overline{\mathcal U}\cup\overline{\Sigma_R})$
on which $u_*$ is smooth, including zero when $u_*$ is smooth there,
\[
 \|g_\lambda-g\|_{C^k}
 +\|u_\lambda-u_*\|_{C^k(K)}
 +\|u_\lambda-u_*\|_{X(\mathcal U^c)}\longrightarrow0.
\]
\end{proposition}
\begin{proof}
Choose $d_*<\beta<(n-2)/4$ and set $\mu=\lambda^\beta$ and $L=\rho/\lambda$, so that
\[
 \mu L^{d_*}=\rho^{d_*}\lambda^{\beta-d_*}\to0,\qquad
 \frac{\lambda^{(n-2)/2}}{\mu^2}
       =\lambda^{(n-2)/2-2\beta}\to0.
\]
Fix a smooth radial cutoff $\chi\in C_c^\infty(\mathcal B_2)$ with $\chi=1$ on $\mathcal B_1$
and $0\le\chi\le1$. Define $g_\lambda=g$ outside $\mathcal U$ and, on $\mathcal U$,
\[
 g_\lambda^{-1}(x)
 =\exp\left(\mu\chi\Bigl(\frac{x-p}{\rho}\Bigr)H\Bigl(\frac{x-p}{\lambda}\Bigr)\right).
\]
Writing $H=\sum_{d\le d_*}H_d$ into homogeneous components gives
\[
 \mu\chi\Bigl(\frac{x-p}{\rho}\Bigr)H\Bigl(\frac{x-p}{\lambda}\Bigr)
 =\chi\Bigl(\frac{x-p}{\rho}\Bigr)\sum_{d\le d_*}\lambda^{\beta-d}H_d(x-p),
\]
so for every fixed $k\ge0$,
\[
 \|g_\lambda-g\|_{C^k(\mathcal U)}
 \le C_{\rho,k}\lambda^{\beta-d_*}\longrightarrow0.
\]
On $\mathcal U$, the conditions $H_{in}=0$ and $\tr H=0$ imply
\[
 \begin{gathered}
 (g_\lambda)_{nn}=1,\qquad (g_\lambda)_{in}=0\quad\text{for }i<n,\qquad
 \det g_\lambda=1,\\
 \dd\sigma_{g_\lambda}=\dd x',\qquad h_{g_\lambda}=0,
 \end{gathered}
\]
and $(g_\lambda)_T|_\Gamma=\delta$ whenever $H|_\Gamma=0$.

By the implicit function theorem and \eqref{eq:insertion-positivity},
there is a unique continuation $u_{*,\lambda}>0$ solving
\[
 \mathcal R_{g_\lambda}(u_{*,\lambda})=0,
 \qquad \|u_{*,\lambda}-u_*\|_X\longrightarrow0,
\]
which converges smoothly to $u_*$ on $\mathcal U$.
Take $D=\mathcal U$ in the local framework, with $\chi_D=1$ on
$\mathcal B_{2\rho}^+(p)$.
Set $W_\theta=u_{*,\lambda}+U_{\lambda,\theta}$ for $\theta\in\Theta$.

The following residual and energy estimates hold after up to two
$\theta$-derivatives, since the corresponding bubble derivatives
satisfy the same spatial decay bounds.
In the rescaled variable $y=(x-p)/\lambda$ with $r=|y|$, the linear metric source is
$\mu O((1+r)^{d_*-n})$. Its $L^{\frac{2n}{n+2}}$ norm on $r>L$ and the
quadratic volume source norm satisfy
\begin{align}
 \mu\left(\int_L^\infty
 r^{n-1+(d_*-n)\frac{2n}{n+2}}\dd r\right)^{\frac{n+2}{2n}}
 &\le C\mu L^{d_*-\frac{n-2}{2}},
 \label{eq:metric-linear-source}\\
 \mu^2\left(\int_1^\infty
 r^{n-1+(2d_*-n)\frac{2n}{n+2}}\dd r\right)^{\frac{n+2}{2n}}
 &\le C\mu^2\qquad\text{for }n>4d_*+2.
 \label{eq:metric-quadratic-source}
\end{align}
On $\mathcal B_{2\rho}(p)\cap\Gamma$, $B_{g_\lambda}U_{\lambda,\theta}-(n-2)U_{\lambda,\theta}^{n/(n-2)}=0$.
Higher-order metric sources carry an extra factor $\mu L^{d_*}\to0$.
The bubble cutoff error satisfies
\begin{equation}\label{eq:metric-bubble-cutoff}
 \|\mathcal R_\delta(U_{\lambda,\theta})\|_{C^2(\Theta;Y)}
 \le C\lambda^{\frac{n-2}{2}}.
\end{equation}
Since $0<c\le u_{*,\lambda}\le C$ on $\mathcal U$, the boundary interaction satisfies,
for $0\le j\le2$,
\[
 \bigl|D_\theta^j\bigl((u_{*,\lambda}+U_{\lambda,\theta})^{\frac{n}{n-2}}
 -u_{*,\lambda}^{\frac{n}{n-2}}-U_{\lambda,\theta}^{\frac{n}{n-2}}\bigr)\bigr|
 \le C\min\bigl(U_{\lambda,\theta}^{\frac{2}{n-2}},U_{\lambda,\theta}\bigr).
\]
Writing $r=|x'-p|$ for the physical boundary radius and splitting
at $r=\sqrt\lambda$, we obtain, for $n\ge5$,
\begin{equation}\label{eq:metric-boundary-interaction}
 \begin{gathered}
 \int_{r<\sqrt\lambda}
      \left(\frac{\lambda}{\lambda^2+r^2}\right)^{\frac{2(n-1)}n} r^{n-2}\dd r
      \le C\lambda^{\frac{n-1}{2}},\\
 \int_{\sqrt\lambda<r<\rho}
      \lambda^{\frac{(n-1)(n-2)}n}r^{-\frac{(n-2)^2}n}\dd r
      \le C\lambda^{\frac{n-1}{2}},
 \end{gathered}
\end{equation}
so its $L^{\frac{2(n-1)}{n}}(\Gamma)$ norm is $O(\lambda^{n/4})=o(\mu)$.
With $\mathcal R_\theta^{(1)}$ denoting the first-order geometric residual on $\mathcal U$,
\begin{equation}\label{eq:metric-insertion-residual}
 \norm{\mathcal R_{g_\lambda}(W_\theta)-\mu\mathcal R_\theta^{(1)}}
                         _{C^2(\Theta;Y)}
 \le C\bigl(\mu L^{d_*-\frac{n-2}{2}}+\mu^2+\lambda^{\frac{n-2}{2}}+\lambda^{n/4}\bigr)
                                      =o(\mu).
\end{equation}
The four terms in \eqref{eq:metric-insertion-residual} come respectively
from truncation of the linear geometric source
\eqref{eq:metric-linear-source}, the nonlinear metric remainder
\eqref{eq:metric-quadratic-source}, the bubble cutoff
\eqref{eq:metric-bubble-cutoff}, and the boundary interaction
\eqref{eq:metric-boundary-interaction}.

Concavity bounds the absolute value of the potential difference by
$C\min(1,U_{\lambda,\theta}^{2/(n-2)})$.
The same split then gives, uniformly for $\theta\in\Theta$,
\begin{equation}\label{eq:insertion-potential-interaction}
 \left\|(u_{*,\lambda}+U_{\lambda,\theta})^{\frac2{n-2}}
       -u_{*,\lambda}^{\frac2{n-2}}
       -U_{\lambda,\theta}^{\frac2{n-2}}\right\|_{L^{n-1}(D\cap\Gamma)}
 \le C\sqrt\lambda.
\end{equation}
Together with the convergence of $g_\lambda$ and $u_{*,\lambda}$,
this implies
\[
 \mathcal L_{\lambda,\theta}
 =\mathcal L_{g,u_*}-nU_{\lambda,\theta}^{\frac2{n-2}}
       +o_{X\to Y}(1)
\]
uniformly on $\Theta$. Lemma~\ref{lem:insertion-inverse} therefore
applies to $\mathcal A_{\lambda,\theta}$.

Let $\psi_{-H;\theta}^\perp$ be the flat correction from Lemma~\ref{lem:jet-inverse}
(using \eqref{eq:constant-correction} for degree zero), and define the cut-off,
moment-adjusted correction on $\mathcal U$ by
\[
 \psi_\theta=\widetilde\psi_\theta-\sum_{\alpha,\beta}(G^{-1})_{\alpha\beta}\ell_\beta(\widetilde\psi_\theta)Z_\alpha,
 \qquad \widetilde\psi_\theta(x)=\chi_D(x)\lambda^{-\frac{n-2}{2}}\psi_{-H;\theta}^\perp\Bigl(\frac{x-p}{\lambda}\Bigr).
\]
By \eqref{eq:correction-decay},
\[
 \sum_{|\gamma|\le1}\int_{|y|>L}\bigl|\nabla D_\theta^\gamma\psi_{-H;\theta}^\perp(y)\bigr|^2\dd y
 \le CL^{2d_*+2-n}\longrightarrow0,
\]
and $\|\mathcal L_{\lambda,\theta}\psi_\theta+\mathcal R_\theta^{(1)}\|_{C^1(\Theta;Y)}=o(1)$.
Since $J_1(H,U_\theta)=0$ and the cubic metric and interaction energy
errors have $C^2(\Theta)$ norm bounded by
\[
 C\mu^3L^{d_*}+C\lambda^{\frac{n-2}{2}}u_*(p)=o(\mu^2),
\]
Lemma~\ref{lem:local-insertion-reduction} applies with $\mathcal C(\theta)=\mathcal C(H;\theta)$
and yields
\begin{equation*}
 \mathcal F_\lambda(\theta)=\mu^2D_\theta\mathcal C(H;\theta)+o_{C^1(\Theta)}(\mu^2).
\end{equation*}
Nondegeneracy at $\theta_*$ gives an exact solution $u_\lambda=W_{\theta_\lambda}+w_{\theta_\lambda}$
with $\theta_\lambda\to\theta_*$ and invertible linearization $\mathcal L_{g_\lambda,u_\lambda}:X\to Y$.
By \eqref{eq:insertion-positivity}, $u_\lambda>0$.
Since $\lambda^{(n-2)/2}u_\lambda(p+\lambda y)\to U_{\theta_*}(y)$ in $C^2_{\mathrm{loc}}(\overline{\R^n_+})$,
$u_\lambda$ has a local boundary maximum $\xi_\lambda=p+O(\lambda)$ with
$u_\lambda(\xi_\lambda)\ge c\lambda^{-(n-2)/2}$, and satisfies the stated exterior estimates.
\end{proof}

\begin{proof}[The metric cases $n\ge19$ of Theorem~\ref{thm:sharpness}]
Take $H_n$ from Lemma~\ref{lem:metric-wells}, which has degree
\[
 d_*=\begin{cases}4,&n=19,\\3,&n\ge20,\end{cases}
 \qquad n>4d_*+2.
\]
Applying Proposition~\ref{prop:complete-insertion} and Proposition~\ref{prop:insertion-limit}
starting from $u_0$ and $u_c$ yields singular solutions with incomplete and complete
conformal metrics at zero, respectively.
At every stage, $H_n|_\Gamma=0$ and $\det g=1$, so $g_T|_\Gamma=I$, $h_g=0$,
and $D^\alpha(g_\infty-\delta)(0)=0$ for all $\alpha$.
\end{proof}

For the Euclidean case, we choose a boundary graph polynomial $f$ whose
completed coefficient \eqref{eq:shape-energy-coefficient} has a nondegenerate
critical point with negative value.
Set $Q_1(z)=z_1^2-z_2^2$ and $Q_2(z)=z_3^2-z_4^2$, and consider
\begin{equation}\label{eq:degree-six-graph}
 f(z)=\sum_{\alpha=1}^2\sum_{j=0}^2c_{\alpha j}Q_\alpha(z)|z|^{2j}.
\end{equation}

\begin{lemma}\label{lem:euclidean-well}
For every $n\ge23$, there is a degree-six polynomial $f$ of the
form \eqref{eq:degree-six-graph}, with $c_{10}=1$, such that
$(s,\xi)=(1,0)$ is a nondegenerate critical point of
$\mathcal C(h_f;s,\xi)$ and
\[
 \mathcal C(h_f;1,0)<0.
\]
\end{lemma}
\begin{proof}
\emph{Step 1. The positive shape form and the two copies.}
The spherical harmonic decomposition of \eqref{eq:shape-form}
shows that $Q_\partial\ge0$ and $\ker Q_\partial$
is spanned by spherical harmonics of degrees $\ell=0,1$ on $\Sph^{n-1}$.
Since $\widehat{Q_1|z|^{2i}}$ ($0\le i\le2$) have tangential harmonic degree $2$,
they are orthogonal to $\ker Q_\partial$, and
\[
 Q_1(z)(c_0+c_1|z|^2+c_2|z|^4)\equiv0
       \quad\Longrightarrow\quad c_0=c_1=c_2=0.
\]
Thus their Gram matrix $G$, defined by
\[
 G_{ij}=Q_\partial\!\left(\widehat{Q_1|z|^{2i}},
                         \widehat{Q_1|z|^{2j}}\right)
             \qquad\text{for }0\le i,j\le2,
\]
is positive definite.
Invariance of $Q_\partial$ under $(z_1,z_2)\mapsto(z_2,z_1)$ and
$(z_1,z_2)\leftrightarrow(z_3,z_4)$ yields
\[
 Q_\partial\!\left(\widehat{Q_1|z|^{2i}},
                   \widehat{Q_2|z|^{2j}}\right)=0,\qquad
 Q_\partial\!\left(\widehat{Q_2|z|^{2i}},
                   \widehat{Q_2|z|^{2j}}\right)=G_{ij}.
\]
For $c=(c_1,c_2)$, set $E_c(s,\xi)=Q_\partial(\widehat f_{s,\xi},\widehat f_{s,\xi})$.
Since $\deg(Q_\alpha|z|^{2j})=2j+2$ induces scale degree $d_j=2j+1\in\{1,3,5\}$,
$E_c(s,0)=E_c(s)$ has the form in Lemma~\ref{lem:profile-selection}.

\emph{Step 2. A negative value of the scale derivative.}
By \eqref{eq:euclidean-negative-dilation}, the trial profile
$f(z)=Q_1(z)(1-|z|^2/4)^2$ satisfies
\[
 \frac{q(c)}{G_{00}}<-\frac1{128}
 \qquad\text{for }n\ge23.
\]

\emph{Step 3. The center derivatives and coefficient selection.}
Reflection invariance in each coordinate implies that $E_c(s,\xi)$ is even in each $\xi_i$, so
\[
 \nabla_\xi E_c(s,0)=0,\qquad
 \partial_{s\xi_i}E_c(s,0)=0,\qquad
 \partial_{\xi_i\xi_j}E_c(s,0)=0\quad\text{for }i\ne j.
\]
Set $h_i(c)=\partial_{\xi_i\xi_i}E_c(1,0)$.
For $f=Q_1$,
\[
 Q_1(z+\xi)-Q_1(z)
 =2\xi_1z_1-2\xi_2z_2+Q_1(\xi)\in\Span\{1,z_1,\ldots,z_{n-1}\},
\]
whose compactification lies in $\ker Q_\partial$, so
\[
 h_i(c)=0,\qquad q(c)=2G_{00}>0.
\]
For $f=Q_1|z|^2$, differentiating $E_c(1,\xi)$ gives
\[
 h_i(c)
 =2Q_\partial(\widehat{\partial_i f},\widehat{\partial_i f})
       +2Q_\partial(\widehat f,\widehat{\partial_{ii}f}),
 \qquad \partial_{ii}(Q_1|z|^2)=2(\partial_{ii}Q_1)|z|^2+4z_i\partial_iQ_1+2Q_1.
\]
Projecting $\partial_{ii}(Q_1|z|^2)$ onto the tangential harmonic $Q_1$ gives
$6Q_1$ for $i=1,2$ and $2Q_1$ for $3\le i\le n-1$.
Since $Q_\partial\ge0$, \eqref{eq:euclidean-positive-coupling} yields
\[
 h_i(c)\ge 2Q_\partial(\widehat f,\widehat{\partial_{ii}f})
 =\begin{cases}
 12G_{01},&i=1,2,\\
 4G_{01},&3\le i\le n-1
 \end{cases}
 \ge4G_{01}>0.
\]
Hence $h_i\notin\R q$ for $1\le i\le n-1$.
Lemma~\ref{lem:profile-selection} provides $c$ with $c_{10}\ne0$, $c_{12}^2+c_{22}^2>0$,
\[
 \partial_sE_c(1,0)=0,\qquad
 \partial_{ss}E_c(1,0)<0,\qquad
 h_i(c)\ne0\quad\text{for }1\le i\le n-1,
\]
and
\[
 D^2_{(s,\xi)}E_c(1,0)
 =\diag\bigl(E_c''(1),h_1(c),\ldots,h_{n-1}(c)\bigr).
\]
Normalizing by $c\mapsto c/c_{10}$ gives $c_{10}=1$ and $\deg f=6$.
By \eqref{eq:shape-energy-coefficient}, $\mathcal C(h_f;s,\xi)$ is a negative
constant multiple of $E_c(s,\xi)$, so $(1,0)$ is a nondegenerate critical point
with $\mathcal C(h_f;1,0)<0$.
\end{proof}

Fix $n\ge23$, the polynomial $f$ from Lemma~\ref{lem:euclidean-well},
$\theta_*=(1,0)$, and a small closed box $\Theta$ around $\theta_*$ containing
no other critical point.
Let $u_*>0$ be an exact solution in a Euclidean domain $\Omega$ with
$Q_\delta$ coercive on energy functions with zero outer spherical
trace and $\mathcal L_{\delta,u_*}:X\xrightarrow{\cong}Y$.
Fix a planar boundary half-ball $\mathcal U=\mathcal B_{8\rho}^+(p)$ centered at
$p\ne0$ avoiding $\{0\}$, the earlier graph supports and $\overline{\Sigma_R}$.

\begin{proposition}\label{prop:euclidean-insertion}
For all sufficiently small $\lambda>0$, there exist a smooth graph
deformation $T_\lambda$, equal to the identity outside $\mathcal U$,
and a positive solution $u_\lambda$ of \eqref{eq:main} in $T_\lambda(\Omega)$
such that
\[
 \mathcal L_{T_\lambda^*\delta,\,u_\lambda\circ T_\lambda}:
                  X\longrightarrow Y
\]
is invertible, with a local boundary maximum
$\xi_\lambda\in T_\lambda(\mathcal U\cap\Gamma)$ satisfying
\[
 |\xi_\lambda-p|\le C\lambda,\qquad
 u_\lambda(\xi_\lambda)\ge c\lambda^{-\frac{n-2}{2}}.
\]
For each fixed $k\ge0$ and exterior compact set $K$,
\[
 \|T_\lambda-I\|_{C^k}
 +\|u_\lambda\circ T_\lambda-u_*\|_{C^k(K)}
 +\|u_\lambda\circ T_\lambda-u_*\|_{X(\mathcal U^c)}
 \longrightarrow0.
\]
\end{proposition}
\begin{proof}
\emph{Step 1. The graph deformation and continuation of $u_*$.}
Put $\beta=41/8$, $\mu=\lambda^\beta$ and $L=\rho/\lambda$.
Fix smooth cutoffs $\chi\in C_c^\infty(B_2)$ with $\chi=1$ on $B_1$ and
$\eta\in C_c^\infty([0,2))$ with $\eta=1$ on $[0,1]$, $0\le\chi,\eta\le1$.
Define
\begin{equation}\label{eq:euclidean-amplitude}
 f_\lambda(x')=\mu\lambda\chi\Bigl(\frac{x'-p}{\rho}\Bigr)f\Bigl(\frac{x'-p}{\lambda}\Bigr),
\end{equation}
and the vertical extension
\begin{equation}\label{eq:graph-extension}
 T_\lambda(x',x_n)=\bigl(x',x_n+f_\lambda(x')\eta(x_n/\rho)\bigr).
\end{equation}
Then $\supp f_\lambda\subset B_{2\rho}(p)$ and $\supp(T_\lambda-I)\subset\mathcal B_{3\rho}^+(p)$.
By \eqref{eq:degree-six-graph} and \eqref{eq:euclidean-amplitude},
\[
 f_\lambda(x')=\chi\Bigl(\frac{x'-p}{\rho}\Bigr)
 \sum_{\alpha=1}^2\sum_{j=0}^2c_{\alpha j}\lambda^{\beta-1-2j}
                Q_\alpha(x'-p)|x'-p|^{2j}.
\]
Since $\min_{0\le j\le2}(\beta-1-2j)=\beta-5=1/8$, we have for every $k\ge0$,
\[
 \|T_\lambda-I\|_{C^k}\le C_k\lambda^{1/8},
\]
so $T_\lambda$ is a diffeomorphism for small $\lambda$.
The invertible linearization at $u_*$ and the implicit function theorem
continue $u_*$ to an exact solution on $T_\lambda(\Omega)$.
We denote its pullback to $\Omega$ by $u_{*,\lambda}$, so that
$\|u_{*,\lambda}-u_*\|_X\to0$.
When starting from $u_0$, the outer spherical trace is fixed and positive.
When starting from $u_c$, the continuation preserves $\Pi_+w(0)=0$.
The small correction keeps the spherical trace positive, and
\eqref{eq:weighted-strip-regularity} ensures positivity near zero.
The small metric deformation also preserves coercivity on energy
functions with zero outer spherical trace. The negative-part test
\eqref{eq:insertion-positivity} therefore gives $u_{*,\lambda}>0$.

In rescaled coordinates $y=(x-p)/\lambda$ with $r=|y|$, write
$f_L(y')=\chi(y'/L)f(y')$ and
\[
 v=\mu\eta(y_n/L)\nabla f_L(y'),\qquad
 \sigma=\mu L^{-1}\eta'(y_n/L)f_L(y').
\]
Then $DT_\lambda=\left(\begin{smallmatrix}I&0\\v^t&1+\sigma\end{smallmatrix}\right)$,
and the pulled-back Dirichlet matrix
$A_\lambda=(\det DT_\lambda)(DT_\lambda)^{-1}(DT_\lambda)^{-t}$ satisfies
\begin{equation*}
 A_\lambda=
 \begin{pmatrix}(1+\sigma)I&-v\\
                -v^t&(1+|v|^2)/(1+\sigma)\end{pmatrix},
\end{equation*}
so
\[
 A_\lambda-I=
 \begin{pmatrix}\sigma I&-v\\-v^t&-\sigma\end{pmatrix}
 +\begin{pmatrix}0&0\\0&\sigma^2+|v|^2\end{pmatrix}
 -\frac{\sigma(\sigma^2+|v|^2)}{1+\sigma}e_n\otimes e_n.
\]
Here $|v|+|\sigma|\le C\mu(1+r)^5$ and $|Dv|+|D\sigma|\le C\mu(1+r)^4$,
with $\|v\|_\infty+\|\sigma\|_\infty\le C\mu L^5\to0$.
On $\Gamma$, the area element and mean-curvature density are
\[
 J_\Gamma=(1+|\nabla f_\lambda|^2)^{1/2},\qquad
 hJ_\Gamma=\frac1{n-1}
  \left(\delta_{ij}-\frac{(f_\lambda)_i(f_\lambda)_j}
                              {1+|\nabla f_\lambda|^2}\right)
                              (f_\lambda)_{ij}.
\]

\emph{Step 2. The residual expansion.}
Take the cutoff bubble $U_{\lambda,\theta}$ equal to the uncut bubble on
$\mathcal B_{3\rho}^+(p)$ and zero outside $\mathcal B_{4\rho}^+(p)$,
and set $W_\theta=u_{*,\lambda}+U_{\lambda,\theta}$ on $D=\mathcal U=\mathcal B_{8\rho}^+(p)$.
From $|\nabla f_L|\le C(1+r)^5$, $|\nabla^2f_L|\le C(1+r)^4$ and
$|D^jU_\theta|\le C_j(1+r)^{2-n-j}$, the linear volume and face sources are
$\mu O((1+r)^{5-n})$ and $\mu O((1+|y'|)^{6-n})$, whose
$L^{\frac{2n}{n+2}}(\{r>L\})$ and $L^{\frac{2(n-1)}{n}}(\{|y'|>L\})$ norms
are $O(\mu L^{5-(n-2)/2})$.
The quadratic volume density is $\mu^2O((1+r)^{10-n})$, and for $n>22$,
\[
 \int_1^\infty r^{n-1+(10-n)\frac{2n}{n+2}}\dd r
 =\frac{n+2}{n(n-22)}<\infty,
\]
giving norm $O(\mu^2)$. On $\Gamma$, $\partial_{y_n}U_\theta=-(n-2)U_\theta^{n/(n-2)}$,
so the quadratic normal-flux and area densities are $\mu^2O((1+|y'|)^{10-n})$,
and for $n>20$,
\[
 \int_1^\infty r^{n-2+(10-n)\frac{2(n-1)}{n}}\dd r
 =\frac{n}{(n-1)(n-20)}<\infty.
\]
The cubic and higher terms have dual norm bounded by
\[
 C\mu^3\begin{cases}
 L^{16-n/2},&23\le n<32,\\
 1+\log L,&n=32,\\
 1,&n>32
 \end{cases}
 \le C\mu^2\lambda^{5/8}(1+|\log\lambda|)=o(\mu^2).
\]
The bubble cutoff error is $O(\lambda^{\frac{n-2}{2}})$, and the interaction term
\[
 \sum_{j=0}^2\bigl|D_\theta^j\bigl((u_{*,\lambda}+U_{\lambda,\theta})^{\frac{n}{n-2}}
                  -u_{*,\lambda}^{\frac{n}{n-2}}
                  -U_{\lambda,\theta}^{\frac{n}{n-2}}\bigr)\bigr|
 \le C\min\bigl(U_{\lambda,\theta}^{\frac{2}{n-2}},U_{\lambda,\theta}\bigr)
\]
has $L^{\frac{2(n-1)}{n}}(\Gamma)$ norm $O(\lambda^{n/4})$, as in Proposition~\ref{prop:complete-insertion}.
Consequently, with $\mathcal R_\theta^{(1)}$ the first-order geometric residual on $D$,
\begin{equation}\label{eq:euclidean-residual}
 \norm{\mathcal R_{T_\lambda^*\delta}(W_\theta)-\mu\mathcal R_\theta^{(1)}}_{C^2(\Theta;Y)}
 \le C\bigl(\mu L^{5-\frac{n-2}{2}}+\mu^2
                         +\lambda^{\frac{n-2}{2}}+\lambda^{n/4}\bigr)=o(\mu).
\end{equation}
These estimates are local to the fixed half-ball $D$, where $u_*$ is
smooth. Since $p\ne0$, the weights defining $X_\delta$ and $Y_\delta$
are smooth and bounded above and below there. Thus the same residual
estimate holds in the weighted spaces when starting from $u_c$.

\emph{Step 3. The direct energy and first correction.}
Integration by parts reduces the first graph variation to a face
integral, including the terms arising from the normal cutoff in
$T_\lambda$. Its density vanishes by the identity
\[
 \tfrac12\bigl((\partial_{y_n}U)^2-|\nabla_{y'}U|^2\bigr)
                    +\frac{n-2}{4(n-1)}\Delta_{y'}(U^2)=0,
\]
and differentiating in $\theta$ gives
\[
 \mathcal R_\theta^{(1)}[Z_\alpha]=0\quad\text{for }0\le\alpha<n.
\]
Using $\mathcal R_{T_\lambda^*\delta}(u_{*,\lambda})=0$, we obtain the local energy expansion
\begin{equation}\label{eq:euclidean-raw-action}
 J_{T_\lambda^*\delta,D}(W_\theta)-J_{T_\lambda^*\delta,D}(u_{*,\lambda})
 =J_\delta(U)+\mu^2J_2(f;U_\theta)+o_{C^2(\Theta)}(\mu^2),
\end{equation}
where the remainder is bounded by
$C(\lambda^{n-2}+\mu^2L^{12-n}+\mu^3+\lambda^{(n-2)/2})=o(\mu^2)$ since
\begin{equation}\label{eq:euclidean-scale-margins}
 5<\beta=\frac{41}{8}<\frac{n-2}{4},\qquad
 \frac{\lambda^{(n-2)/2}}{\mu^2}\le\lambda^{1/4},\qquad
 \frac{\lambda^{n/4}}{\mu}\le\lambda^{5/8}.
\end{equation}

The graph deformation tends to the identity in $C^2$, so its
contribution to the linearized operator tends to zero in operator norm.
Estimate~\eqref{eq:insertion-potential-interaction} controls the boundary
interaction with $u_{*,\lambda}$. The new bubble potential has the exterior $L^{n-1}$
bound $C/K$ on $|x-p|>K\lambda$, as in the proof of
Lemma~\ref{lem:insertion-inverse}. These facts verify the operator
hypotheses of that lemma and give the uniform inverse for
$\mathcal A_{\lambda,\theta}$.
Let $\psi_\theta^0$ be the flat Dirichlet-orthogonal correction for the
vertical graph variation, which differs from the Fermi graph variation $h_f$
by a boundary-tangent coordinate variation. By diffeomorphism invariance
\cite[Section~2.2]{CompanionI}, \eqref{eq:reduced-form} and
\eqref{eq:shape-energy-coefficient},
\[
 \mathcal C(h_f;\theta)
 =J_2(f;U_\theta)-\tfrac12Q_\theta(\psi_\theta^0,\psi_\theta^0).
\]
Rescaling $\psi_\theta^0$ into $D$, cutting off and adjusting its moments
defines $\psi_{\lambda,\theta}$ with $\ell(\psi_{\lambda,\theta})=0$.
Local coefficient convergence and the decay of $\psi_\theta^0$ imply
\[
 \|\mathcal L_{\lambda,\theta}\psi_{\lambda,\theta}
             +\mathcal R_\theta^{(1)}\|_{C^1(\Theta;Y)}=o(1).
\]
The same decay controls the omitted exterior Dirichlet integrals,
\[
 \sum_{|\gamma|\le1}\int_{|y|>L}\bigl|\nabla D_\theta^\gamma\psi_\theta^0(y)\bigr|^2\dd y
 =O(L^{12-n})\longrightarrow0.
\]
Equations~\eqref{eq:euclidean-residual} and \eqref{eq:euclidean-raw-action}
verify all hypotheses of Lemma~\ref{lem:local-insertion-reduction}.

\emph{Step 4. The exact solution and exterior convergence.}
Lemma~\ref{lem:local-insertion-reduction} gives
\[
 \mathcal F_\lambda(\theta)=\mu^2D_\theta\mathcal C(h_f;\theta)+o_{C^1(\Theta)}(\mu^2).
\]
By Lemma~\ref{lem:euclidean-well}, $\mathcal F_\lambda$ has a zero
$\theta_\lambda\to\theta_*$, yielding an exact pulled-back solution
$u_\lambda\circ T_\lambda=W_{\theta_\lambda}+w_{\theta_\lambda}$ with invertible full
linearization. Its spherical trace is positive, and in the cylindrical
construction it is positive near zero. Coercivity and
\eqref{eq:insertion-positivity} then imply $u_\lambda>0$.
Strong rescaled energy convergence and the local Robin estimates give
a boundary maximum at $p+O(\lambda)$ with height at least
$c\lambda^{-(n-2)/2}$.

On exterior compact sets where $u_*$ is smooth, the same estimates
give smooth convergence to $u_*$, including at zero when starting
from $u_0$. The energy estimates give convergence in $X(\mathcal U^c)$.
When starting from $u_c$, local estimates on fixed annuli and
\eqref{eq:weighted-strip-regularity} also yield
\[
 \sup_{\mathcal U^c}
 \frac{|u_\lambda\circ T_\lambda-u_*|}{u_c}\longrightarrow0.
\]
This is the relative exterior convergence required in
Proposition~\ref{prop:insertion-limit}.
\end{proof}

\begin{proof}[The Euclidean cases of Theorem~\ref{thm:sharpness}]
Use the disjoint half-balls $D_j=\mathcal B_{8\rho_j}^+(p_j)$ of
Proposition~\ref{prop:insertion-limit}. At stage $j$, the graph
\eqref{eq:euclidean-amplitude} is
\[
 f_j(p_j+w)=\lambda_j^{\beta+1}\chi(w/\rho_j)f(w/\lambda_j),
 \qquad\beta=\frac{41}{8},
\]
which satisfies $\|f_j\|_{C^k}\le C_{j,k}\lambda_j^{\beta-5}\to0$ as $\lambda_j\to0$.
Choosing $\lambda_j$ small enough to satisfy the bounds of
Proposition~\ref{prop:insertion-limit}, the compositions $T^{(N)}=T_{\lambda_N}\circ\cdots\circ T_{\lambda_1}$
converge in $C^\infty$ to a diffeomorphism $T^{(\infty)}$ with
\[
 D^\alpha(T^{(\infty)}-I)(0)=0\quad\text{for all }\alpha,\qquad
 g_\infty=(T^{(\infty)})^*\delta.
\]
Then $u_\infty\circ(T^{(\infty)})^{-1}$ solves \eqref{eq:main} on the Euclidean domain
$T^{(\infty)}(\mathcal B_R^+)$, whose boundary has second fundamental form and mean curvature at $p_j$ given by
\[
 \mathrm{II}(p_j)=\lambda_j^{\beta-1}D^2f(0)\ne0,\qquad
 h(p_j)=\frac{1}{n-1}\tr D^2f(0)=0.
\]

When starting from $u_0$, Proposition~\ref{prop:euclidean-insertion} and
Proposition~\ref{prop:insertion-limit} give the incomplete Euclidean example.
When starting from $u_c$, choose the $C^2$ graph changes summably small so that
$h_-=\max(-h,0)$ satisfies, with $C$ the uniform Sobolev trace constant,
\[
 \frac{n-2}{2}C\|h_-\|_{L^{n-1}}< \frac12
\]
at every stage. Then for every $v\in X$, H\"older's and the trace inequalities yield
\[
 \int|\nabla v|^2+\frac{n-2}{2}\int hv^2
 \ge\left(1-\frac{n-2}{2}C\|h_-\|_{L^{n-1}}\right)\int|\nabla v|^2
 \ge\frac12\int|\nabla v|^2.
\]
Thus $Q_g$ remains coercive at every stage, and Lemma~\ref{lem:core-barrier}
gives $u_j\ge c_0u_c$ on $D_j$. Proposition~\ref{prop:insertion-limit}
then yields the complete Euclidean example.
\end{proof}

For totally geodesic boundaries, we use metric polynomials independent of $x_n$,
cut off evenly in $x_n$ so that $\mathrm{II}_{ij}=-\frac12\partial_{x_n}g_{ij}(x',0)=0$.
Take two nonzero algebraic Weyl tensors $W^1,W^2$ supported in disjoint
four-dimensional tangential subspaces, normalized by
\[
 (\Omega_\alpha)_{pq}=\sum_{i,j,l}
   (W^\alpha_{ipjl}+W^\alpha_{iljp})
   (W^\alpha_{iqjl}+W^\alpha_{iljq}),
 \qquad \tr\Omega_1=\tr\Omega_2=1,\qquad \Omega_1\Omega_2=0.
\]
Fix an orthonormal basis of $\R^{n-1}$ diagonalizing both $\Omega_1$ and $\Omega_2$.
Put $(K_\alpha)_{ij}(x')=W^\alpha_{ikjl}x'_kx'_l$ and consider
\[
 H(x')=\sum_{\alpha=1}^2P_{c_\alpha}(|x'|^2)K_\alpha(x'),
 \qquad P_{c_\alpha}(v)=\sum_{j=0}^2c_{\alpha j}v^j.
\]

\begin{lemma}\label{prop:geodesic-profile}
For every $n\ge35$, there is a degree-six polynomial $H$ in the
two-copy Weyl family above such that $(s,\xi)=(1,0)$ is a nondegenerate critical
point of $\mathcal C(H;s,\xi)$ and
\[
 \mathcal C(H;1,0)<0.
\]
\end{lemma}
\begin{proof}
The algebraic Weyl symmetries imply
\[
 \tr H=0,\qquad \diver H=0,\qquad Hx'=0,\qquad \partial_{x_n}H=0.
\]
Since $\nabla_{x'}U_{s,0}$ is parallel to $x'$, $H\nabla U_{s,0}=0$, so
\eqref{eq:metric-correction-source} gives $F_{-H;s,0}=0$ and $F^\partial_{-H;s,0}=0$.
By Lemma~\ref{lem:jet-inverse}, $\psi_{-H;s,0}=0$.
Substituting $H\nabla U_{s,0}=0$, $\diver H=0$ and $\psi_{-H;s,0}=0$ into
\eqref{eq:direct-form} gives
\[
 -\mathcal C(H;s,0)
 =\frac{n-2}{32(n-1)}\int_{\R^n_+}|\partial H|^2U_{s,0}^2\,\dd x>0
 \qquad\text{for }H\ne0,
\]
where the integrand is $O((1+r)^{10-2(n-2)})$, which is integrable since $n>14$.

Disjointness of the supports of $W^1,W^2$ and $\tr\Omega_1=\tr\Omega_2=1$
imply that the two copies are orthogonal and share the same positive definite
Gram matrix
\[
 G_{ij}=-\mathcal C\bigl(|x'|^{2i}K_1,|x'|^{2j}K_1\bigr)
                  \qquad\text{for }0\le i,j\le2.
\]
With homogeneous degrees $d=(2,4,6)$, we have $-\mathcal C(H;s,0)=E_c(s)$
as in Lemma~\ref{lem:profile-selection}.

Taking $P_{c_1}(v)=(1-v/3)^2$ and $P_{c_2}=0$ yields $q(c)<0$
for every $n\ge35$ by \eqref{eq:geodesic-negative-dilation}.

Since $H(-x')=H(x')$, we have $D_\xi\mathcal C(H;s,0)=0$ and $D_{s\xi}\mathcal C(H;s,0)=0$.
For $U_{s,\xi}$, $Hx'=0$ implies $H_{ij}(x')(x_i-\xi_i)(x_j-\xi_j)=H_{ij}(x')\xi_i\xi_j$,
so \eqref{eq:metric-correction-source} gives
\[
 F_{-H;s,\xi}(x)
 =\frac{n(n-2)s^{(n-2)/2}H_{ij}(x')\xi_i\xi_j}
 {(|x'-\xi|^2+(x_n+s)^2)^{(n+2)/2}},
 \qquad F^\partial_{-H;s,\xi}=0.
\]
Hence $\|F_{-H;s,\xi}\|_{L^{\frac{2n}{n+2}}(\R^n_+)}=O(|\xi|^2)$, and Lemma~\ref{lem:jet-inverse} yields
\[
 \|\psi_{-H;s,\xi}^{\perp}\|_{\mathcal D^{1,2}}
       \le C|\xi|^2,\qquad
 Q_{s,\xi}(\psi_{-H;s,\xi}^{\perp})=O(|\xi|^4),\qquad
 \left.D^2_{\xi\xi}Q_{s,\xi}(\psi_{-H;s,\xi}^{\perp})\right|_{\xi=0}=0.
\]
Thus $D^2_{\xi\xi}\mathcal C(H;1,0)=D^2_{\xi\xi}J_2(H;U_{1,0})$, which by
\eqref{eq:geodesic-center-hessian} is a linear combination of $\Omega_1,\Omega_2,I$
and hence diagonal in the chosen orthonormal basis.
Set $h_i(c)=-\partial_{\xi_i\xi_i}\mathcal C(H;1,0)$.
For $P_{c_1}=1$ and $P_{c_2}=0$, \eqref{eq:geodesic-center-hessian} gives
\[
 h_i(c)=0,\qquad q(c)=4G_{00}>0.
\]
For $P_{c_1}(v)=v$ and $P_{c_2}=0$, \eqref{eq:geodesic-center-hessian} and the
identities following it give, with $c_*>0$,
\[
 D^2_{\xi\xi}\mathcal C(H;1,0)
 =c_*\Bigl(-\Omega_1-\frac{n+6}{6(n+5)}I\Bigr)<0
 \quad\Longrightarrow\quad h_i(c)>0\quad\text{for }1\le i\le n-1.
\]
Thus $h_i\notin\R q$ for $1\le i\le n-1$.
Lemma~\ref{lem:profile-selection} provides $c$ with $c_{12}^2+c_{22}^2>0$,
$E_c'(1)=0$, $E_c''(1)<0$ and $h_i(c)\ne0$ for all $i$, so $\deg H=6$,
$D_{(s,\xi)}\mathcal C(H;1,0)=0$,
\[
 D^2_{(s,\xi)}\mathcal C(H;1,0)
 =\diag\bigl(-E_c''(1),-h_1(c),\ldots,-h_{n-1}(c)\bigr)
\]
is invertible, and $\mathcal C(H;1,0)=-E_c(1)<0$.
\end{proof}

\begin{proof}[The totally geodesic cases of Theorem~\ref{thm:sharpness}]
The polynomial $H$ in Lemma~\ref{prop:geodesic-profile} has $d_*=6$, so
$n>4d_*+2$ for $n\ge35$. Applying Proposition~\ref{prop:complete-insertion}
with the radial cutoff $\chi((x-p)/\rho)$ (which is even in $x_n$) gives
\[
 (g_\lambda)_{nn}=1,\qquad (g_\lambda)_{in}=0\quad\text{for }i<n,\qquad
 \mathrm{II}_{g_\lambda,ij}=-\frac12\partial_{x_n}(g_\lambda)_{ij}\Big|_\Gamma=0
\]
at every stage. Proposition~\ref{prop:insertion-limit} starting from $u_0$ and $u_c$
then yields the incomplete and complete examples with $\mathrm{II}_{g_\infty}\equiv0$
and $D^\alpha(g_\infty-\delta)(0)=0$ for all $\alpha$.
\end{proof}

\section{The eighteenth-dimensional construction}\label{subsec:eighteen-examples}

The successive-insertion argument in
Proposition~\ref{prop:complete-insertion} requires
$d_*<\beta<(n-2)/4$. For $n=18$ and $d_*=4$, this interval is empty.
We instead give separated bubbles a common fourth-order metric jet,
starting with the local polynomial below.

Write $y=(Y,Z,y_n)\in\R^{16}\times\R\times[0,\infty)$, with
$y_n$ the normal variable.
Take $A=\diag(1,-1,0,\ldots,0)$ on $\R^{17}$, with its two nonzero
entries in $Y$, and set
\begin{equation}\label{eq:profile18}
 h_{18}(Y,Z,y_n)=y_n(1-y_n)(3-|Y|^2)-\tfrac1{10}y_nZ^2,
 \qquad H_{18}=h_{18}A.
\end{equation}
Write $\mathcal C_{18}(\theta)=\mathcal C(H_{18};\theta)$ for
$\theta=(s,\xi)\in(0,\infty)\times\R^{17}$.
The quartic term $y_n^2|Y|^2A$ is invariant under $Z$-translations
and will be used as a common metric term, while the cubic term
$-\frac1{10}y_nZ^2A$ supplies the $Z$-center dependence.

\begin{lemma}\label{lem:metric-well-eighteen}
There is $s_{18}\in(15/16,1)$ such that $(s_{18},0)$ is a
nondegenerate strict local minimum of $\mathcal C_{18}$ and
$\mathcal C_{18}(s_{18},0)<0$.
\end{lemma}
\begin{proof}
The endpoint and Bernstein bounds in
Appendix~\ref{app:metric-wells} give, on $I=[15/16,1]$,
\[
 \partial_s\mathcal C_{18}(15/16,0)<0<\partial_s\mathcal C_{18}(1,0),\qquad
 \partial_{ss}\mathcal C_{18}(s,0)>0\quad\text{for }s\in I,
\]
so $s\mapsto\mathcal C_{18}(s,0)$ has a unique minimizer $s_{18}\in(15/16,1)$.
Reflection invariance in each coordinate of $(Y,Z)$ gives
\[
 D_\xi\mathcal C_{18}(s,0)=0,\qquad
 D_{s\xi}\mathcal C_{18}(s,0)=0,\qquad
 \partial_{\xi_i\xi_j}\mathcal C_{18}(s,0)=0\quad\text{for }i\ne j.
\]
The Bernstein lower bounds for $\partial_{\xi_i\xi_i}\mathcal C_{18}(s,0)$
on $I$ ($i=1,3,17$) and \eqref{eq:18-profile-bounds} yield
\[
 D^2_{(s,\xi)}\mathcal C_{18}(s_{18},0)>0,\qquad
 \mathcal C_{18}(s_{18},0)\le\mathcal C_{18}(1,0)<0.
\qedhere
\]
\end{proof}

Choose a small closed box $\Theta$ around $(s_{18},0)$ with boundary
gradient margin $m_*>0$ as in \eqref{eq:profile-margin}, now for
$\mathcal C_{18}$.
Write $x=(Y,Z,x_n)\in\R^{16}\times\R\times[0,\infty)$.
Set
\[
 p_j=2^{-j}e_Z,\qquad \rho_j=2^{-j-6},\qquad
 \lambda_j=\exp(-128M8^j).
\]
Fix disjoint rescaled cutoffs $\chi_j\in C_c^\infty(\mathcal B_{2\rho_j}(p_j))$
with $\chi_j=1$ on $\mathcal B_{\rho_j}(p_j)$, a cutoff $\chi_0\in C_c^\infty(\mathcal B_2)$
with $\chi_0=1$ on $\mathcal B_1$, and a small constant $a_0>0$.
Since all centers lie on the $Z$ axis, the quartic part of every
rescaled profile $\lambda_j^4h_{18}((x-p_j)/\lambda_j)$ is the same
polynomial $x_n^2|Y|^2$. We prescribe this common fourth-order term
through $f_0$ and use $f_j$ to insert the remaining terms near $p_j$.
Extending $A$ by zero in the normal direction, set
\begin{equation}\label{eq:metric18-blocks}
 \begin{aligned}
 f_0(x)&=a_0\chi_0(x)x_n^2|Y|^2,\\
 f_j(x)&=a_0\chi_j(x)\left\{
 \lambda_j^4h_{18}\Bigl(\frac{x-p_j}{\lambda_j}\Bigr)-x_n^2|Y|^2\right\},\\
 g_N^{-1}&=\exp\Bigl(\Bigl(f_0+\sum_{j=1}^Nf_j\Bigr)A\Bigr).
 \end{aligned}
\end{equation}
The quartic terms cancel in $f_j$, leaving
\[
 f_j(x)=a_0\chi_j(x)\left\{
 3\lambda_j^3x_n-3\lambda_j^2x_n^2
 -\lambda_jx_n|Y|^2
 -\frac{\lambda_j}{10}x_n(Z-2^{-j})^2\right\}.
\]
On $\mathcal B_{\rho_j}^+(p_j)$, $\chi_0=\chi_j=1$ and all other
local cutoffs vanish. Consequently, for $1\le j\le N$, the rescaled
metric is
\[
 g_N^{-1}(p_j+\lambda_jy)
 =\exp\bigl(a_0\lambda_j^4H_{18}(y)\bigr)
 \quad\text{for }y\in\mathcal B_{\rho_j/\lambda_j}^+.
\]
The polynomial multiplying $\chi_j$ has degree at most three, and
each term contains a positive power of $\lambda_j$. Differentiating
the rescaled cutoff therefore gives, for every fixed $k,m\ge0$,
\[
 \|D^kf_j\|_\infty\le C_ka_0\lambda_j\rho_j^{3-k}
                              =o(|p_j|^m)\qquad\text{as }j\to\infty.
\]
Hence $g_N\to g_\infty$ in $C^\infty$, with
$g_\infty^{-1}=\exp((f_0+\sum_{j=1}^\infty f_j)A)$ having fourth-order jet
$a_0x_n^2|Y|^2A$ at zero.
The exponential defines a positive definite metric.
Since $A_{in}=0$, $\tr A=0$ and $f_0=f_j=0$ on $\Gamma$, it satisfies
\[
 \begin{gathered}
 (g_N)_{nn}=1,\qquad (g_N)_{in}=0\quad\text{for }i<18,\\
 \det g_N=1,\qquad (g_N)_T|_\Gamma=I,\qquad h_{g_N}=0.
 \end{gathered}
\]

For $\theta=(\theta_1,\ldots,\theta_N)\in\Theta^N$, with
$\theta_j=(s_j,\xi_j)$, write
\[
 \begin{aligned}
 U_j(x)&=\lambda_j^{-8}
       U_{s_j,\xi_j}\Bigl(\frac{x-p_j}{\lambda_j}\Bigr),
       \qquad \mu_j=a_0\lambda_j^4,\\
 W_N(x)&=\sum_{j=1}^NU_j(x).
 \end{aligned}
\]

Fix $a_0>0$ small and $M$ large. Since
$\|g_N-I\|_\infty+\|R_{g_N}\|_{L^9(\R^{18}_+)}\le Ca_0$,
Sobolev absorption gives, uniformly in $N$,
\[
 Q_{g_N}(v,v)\ge c\|\nabla v\|_2^2
                   \quad\text{for all }v\in\mathcal D^{1,2}(\R^{18}_+),
\]
and the Neumann Green kernel of $-L_{g_N}$ satisfies
\[
 0<G_{g_N}(x,y)\le C|x-y|^{-16}.
\]
The kernel comparison is established in the proof of
Lemma~\ref{lem:bridge-potentials}.

For $g=g_N$, let $G_g$ and $G_{g,\Gamma}$ be the volume and boundary Green operators.

Write $d_{ij}=|p_i-p_j|$ and define
\begin{equation}\label{eq:interaction-rows}
 e_{ij}=\left(\frac{\lambda_i\lambda_j}{d_{ij}^2}\right)^{8},
 \qquad
 e_j=\sum_{i\ne j}e_{ij}+\lambda_j^{16}\rho_j^{-16}.
\end{equation}
Since
\[
 \frac{e_j}{\mu_j^2}\le\frac C{a_0^2}
 \left(\sum_{i<j}\lambda_i^8|p_i|^{-16}
       +|p_j|^{-16}\sum_{i>j}\lambda_i^8
       +\Bigl(\frac{\lambda_j}{\rho_j}\Bigr)^8\rho_j^{-8}\right),
\]
we have
\begin{equation}\label{eq:interaction-small}
 \lim_{M\to\infty}\sup_{j\ge1}\frac{e_j}{\mu_j^2}=0.
\end{equation}

Set
\[
 \delta_M=\sup_{i\ge1}\sum_{j\ne i}\frac{\sqrt{\lambda_i\lambda_j}}{d_{ij}}
 \longrightarrow0\quad\text{as }M\to\infty.
\]

Let $Z_{j\alpha}$ ($0\le\alpha<18$) be the
cut-off Jacobi fields supported in $\mathcal B_{2\rho_j}^+(p_j)$ and equal to
the normalized scale and center derivatives of $U_j$ on $\mathcal B_{\rho_j}^+(p_j)$.
Define
\[
 \mathcal X_N=\bigl\{w\in\mathcal D^{1,2}(\R^{18}_+):
       Q_{g_N}(w,Z_{j\alpha})=0
                         \text{ for }1\le j\le N,\ 0\le\alpha<18\bigr\}.
\]
Since $\supp Z_{i\alpha}\cap\supp Z_{j\beta}=\varnothing$ for $i\ne j$,
$Q_{g_N}(Z_{i\alpha},Z_{j\beta})=0$ for $i\ne j$, and the diagonal blocks
$Q_{g_N}(Z_{j\alpha},Z_{j\beta})$ have uniformly bounded inverses.

For the source norm, put
\[
 d_N(x)=\min_{1\le j\le N}(\lambda_j+|x-p_j|).
\]

For energy-dual volume and boundary sources $(f,f^{\partial})$, define
\[
 \norm{(f,f^{\partial})}_{\mathcal Y_N}
 =\sup_{\R^{18}_+}\frac{G_g(|f|)+G_{g,\Gamma}(|f^{\partial}|)}{W_N}
       +\sup_{\R^{18}_+}\frac{d_N^2|f|}{W_N}
       +\sup_\Gamma\frac{d_N|f^{\partial}|}{W_N}.
\]

\begin{lemma}\label{lem:relative-many-inverse}
Let $W\asymp W_N$ be positive, with uniform comparison constants.
Suppose that, for every fixed $R>0$,
\[
 \lim_{M\to\infty}\sup_{N\ge1}\sup_{\theta\in\Theta^N}
 \max_{1\le j\le N}
 \|\lambda_j^8W(p_j+\lambda_j\,\cdot)-U_{s_j,\xi_j}\|
                         _{L^\infty(\mathcal B_R^+)}=0.
\]
For sufficiently large $M$, there exist a unique
$w\in\mathcal X_N$ and a unique multiplier family $\eta=(\eta_{j\alpha})$ satisfying
\[
 \mathcal L_{g,W}w=(f,f^{\partial})+\sum_{j,\alpha}\eta_{j\alpha}
                          Q_g(Z_{j\alpha},\cdot),
\]
and, with $C$ independent of $N$,
\begin{equation}\label{eq:relative-many-inverse}
 \|w/W_N\|_\infty+\sup_j|\eta_j|
                    \le C\|(f,f^{\partial})\|_{\mathcal Y_N}.
\end{equation}
\end{lemma}
\begin{proof}
For each fixed $N$, $W^{1/8}\in L^{17}(\Gamma)$, so $Q_{g,W}$ is a compact
perturbation of the coercive form $Q_g$ on $\mathcal X_N\times\mathcal X_N$
and the projected operator has Fredholm index zero.
The annihilator of $\mathcal X_N$ is spanned by
$Q_g(Z_{j\alpha},\cdot)$, which supplies the multipliers in the
stated equation. We prove the estimate for any energy solution
of this equation.
Since $W_N(x)\asymp(1+|x|)^{-16}$ for $|x|\ge R$ and
\[
 \sup_{\R^{18}_+\setminus\mathcal B_R^+}
 \frac{G_{g,\Gamma}(18W^{1/8}\mathbf1_{\Gamma\setminus B_R}|w|)}{W_N}
 \le CR^{-1}\sup_{\R^{18}_+\setminus\mathcal B_R^+}\frac{|w|}{W_N},
\]
the exterior Neumann series, for large $R$, constructs a solution
with finite relative norm from the spherical boundary values and
the source potential. The small exterior $L^{17}$ norm of $W^{1/8}$
makes the exterior form coercive, so this solution agrees with the
given energy solution. Local Robin estimates on the remaining
half-ball then imply $\|w/W_N\|_\infty<\infty$.

Testing the equation with $Z_{j\alpha}$ gives
\[
 \sum_\beta Q_g(Z_{j\beta},Z_{j\alpha})\eta_{j\beta}
 =\langle\mathcal L_{g,W}Z_{j\alpha},w\rangle
            -\langle(f,f^\partial),Z_{j\alpha}\rangle.
\]
On each fixed rescaled half-ball,
$\lambda_jW^{1/8}(p_j+\lambda_jy)\to U_{s_j,\xi_j}^{1/8}(y)$ and
$g(p_j+\lambda_j\,\cdot)\to\delta$ in $C^2$, uniformly in $j,N$.
The flat Jacobi equation cancels the limiting contribution there.
For the complementary region, the self-bubble moment satisfies
\[
 \int_{\Gamma\setminus B_{R\lambda_j}(p_j)}
                  U_j^{9/8}|Z_{j\alpha}|\,\dd x'\le CR^{-17},
\]
while the other bubbles contribute $C\sum_{i\ne j}e_{ij}$ by
Lemma~\ref{lem:bridge-potentials}. Combining these bounds with the
metric and cutoff estimates gives
\[
 \bigl|\langle\mathcal L_{g,W}Z_{j\alpha},w\rangle\bigr|
 \le\Bigl(o_M(1)+CR^{-17}+C\sum_{i\ne j}e_{ij}+Ca_0\lambda_j\rho_j^3\Bigr)
    \Bigl\|\frac{w}{W_N}\Bigr\|_\infty
 =o(1)\Bigl\|\frac{w}{W_N}\Bigr\|_\infty.
\]
Here $R$ is chosen large before $M$, with all bounds uniform in $j,N$.
The absolute metric sources and cutoff commutators of $Z_{j\alpha}$
have uniformly bounded $W_N$ moments. Green symmetry therefore gives
\[
 |\langle(f,f^\partial),Z_{j\alpha}\rangle|
 \le C\|(f,f^\partial)\|_{\mathcal Y_N}.
\]
Inverting the Gram block $Q_g(Z_{j\alpha},Z_{j\beta})$ yields
\[
 \sup_j|\eta_j|
 \le C\|(f,f^\partial)\|_{\mathcal Y_N}
                    +o(1)\|w/W_N\|_\infty
 \qquad\text{as }M\to\infty.
\]

Suppose the bound on $\|w/W_N\|_\infty$ failed. Then there would exist $M_k\to\infty$
and solutions with
\[
 \|w_k/W_{N_k}\|_\infty=1,
 \qquad
 \|(f_k,f_k^\partial)\|_{\mathcal Y_{N_k}}\longrightarrow0,
 \qquad
 \sup_j|\eta_{k,j}|\longrightarrow0.
\]
For any sequence $j_k\le N_k$,
\[
 \lambda_{j_k}^8|w_k(p_{j_k}+\lambda_{j_k}y)|
 \le\lambda_{j_k}^8W_{N_k}(p_{j_k}+\lambda_{j_k}y)
 \le C_R U(y)
 \qquad\text{for }|y|\le R,
\]
and the pointwise terms in the $\mathcal Y_{N_k}$ norm make the
rescaled sources tend to zero. Local Robin compactness then gives,
after passing to a subsequence,
$\lambda_{j_k}^8w_k(p_{j_k}+\lambda_{j_k}\,\cdot)\to w_\infty$
locally, where $w_\infty$ solves the flat Jacobi equation at $U_{\theta_*}$ with
$|w_\infty|\le CU_{\theta_*}$. Caccioppoli's inequality gives
\[
 \int_{\{R<|y|<2R\}}|\nabla w_\infty|^2\,\dd y\le CR^{-16}
 \qquad\text{for }R\ge1,
\]
so $w_\infty\in\mathcal D^{1,2}(\R^{18}_+)$.
To pass $Q_g(w_k,Z_{j_k,\alpha})=0$ to the limit, first cut off the
Jacobi field at a fixed rescaled radius $R$. Local weak energy
convergence applies on this half-ball. On its complement, integrate
by parts onto the Jacobi field and use the preceding moment and
cutoff estimates. Letting $k\to\infty$ and then $R\to\infty$
preserves the Dirichlet orthogonality. The flat kernel characterization
therefore yields
\[
 w_\infty=\sum_\alpha a_\alpha Z_\alpha,
 \qquad
 \int_{\R^{18}_+}\nabla w_\infty\cdot\nabla Z_\beta\,\dd y=0\quad\text{for all }\beta,
\]
so the positive flat Gram matrix implies $w_\infty=0$.
Since $j_k$ was arbitrary, any failure of uniform convergence over
the bubble indices would produce a nonzero rescaled limit.
Consequently, for every fixed $R>0$,
\[
 \lim_{k\to\infty}\max_{1\le j\le N_k}\sup_{\mathcal B_{R\lambda_j}^+(p_j)}
 \frac{|w_k|}{W_{N_k}}=0.
\]
The Green representation is
\[
 w_k=G_g(f_k)+G_{g,\Gamma}(f_k^\partial)
 +18G_{g,\Gamma}(W^{1/8}w_k)
 +\sum_{j,\alpha}\eta_{k,j\alpha}Z_{j\alpha}.
\]
The source potentials are controlled by the $\mathcal Y_{N_k}$ norm,
and $|Z_{j\alpha}|\le CU_j$ controls the multiplier sum by
$C\sup_j|\eta_{k,j}|W_{N_k}$.
Split the remaining boundary potential over
$\bigcup_j B_{R\lambda_j}(p_j)$ and its complement.
On the disks, the uniform convergence just proved and the first
bound in \eqref{eq:bridge-estimates} give $o_k(1)W_{N_k}$.
On the complement, $|w_k|\le W_{N_k}$ and the second bound there
give $CR^{-1}W_{N_k}$. Dividing by $W_{N_k}$ and taking the supremum yields
\[
 1=\|w_k/W_{N_k}\|_\infty
 \le o_k(1)+CR^{-1}+
 C\bigl(\|(f_k,f_k^\partial)\|_{\mathcal Y_{N_k}}+\sup_j|\eta_{k,j}|\bigr),
\]
a contradiction for $R>2C$ and $k$ large.
For $(f,f^\partial)=0$, this estimate gives $w=0$ and $\eta=0$, so
Fredholm index zero yields existence and uniqueness.
The argument uses only the uniform separation and local profile
bounds, so it also applies to finite subfamilies with their
corresponding metrics.
\end{proof}

To compare each bubble's energy with $\mathcal C_{18}$, let
$g^{(j)}$ be the metric containing only $f_0$ and $f_j$:
\[
 (g^{(j)})^{-1}=\exp\bigl((f_0+f_j)A\bigr),\qquad
 g_0^{-1}=\exp(f_0A).
\]
Set
\[
 \mathcal X^{(j)}=\bigl\{w\in\mathcal D^{1,2}(\R^{18}_+):
                  Q_{g^{(j)}}(w,Z_{j\alpha})=0\text{ for }0\le\alpha<18\bigr\}.
\]
The single-bubble projected equation is
\begin{equation}\label{eq:single-projected-equation}
 J_{g^{(j)}}'(\widetilde U_j)
   =\sum_\alpha \hat c_{j\alpha}
                Q_{g^{(j)}}(Z_{j\alpha},\cdot),
 \qquad \widetilde U_j=U_j+\psi_j,\quad \psi_j\in\mathcal X^{(j)}.
\end{equation}
Since $|N_{U_j}(v)-N_{U_j}(\widetilde v)|\le C\kappa U_j^{1/8}|v-\widetilde v|$
for $|v|,|\widetilde v|\le\kappa U_j$, Lemma~\ref{lem:relative-many-inverse}
with $N=1$ and the flat energy inverse yield a unique solution $\psi_j$ satisfying,
with $\kappa\le Ca_0$,
\begin{equation*}
 |\psi_j|+|D_{\theta_j}\psi_j|\le\kappa U_j,\qquad
 \norm{\psi_j}_{\mathcal D^{1,2}}
       +\norm{D_{\theta_j}\psi_j}_{\mathcal D^{1,2}}\le C\mu_j.
\end{equation*}

Set $R_j=\rho_j/\lambda_j$ and define the energy remainder by
\begin{equation}\label{eq:single-bubble-energy}
 J_{g^{(j)}}(\widetilde U_j)
 =J_\delta(U)+\mu_j^2\bigl(\mathcal C_{18}(\theta_j)+r_j(\theta_j)\bigr).
\end{equation}
On fixed rescaled half-balls, the metric has the expansion
\[
 g^{(j)}(p_j+\lambda_jy)=I-\mu_jH_{18}(y)+O(\mu_j^2).
\]
Thus the flat first correction is $\psi_{-H_{18};\theta_j}^{\perp}$.
Writing $\widehat\psi_j(y)=\lambda_j^8\psi_j(p_j+\lambda_jy)$,
the comparison \eqref{eq:endpoint-response-convergence} gives
\[
 \widehat\psi_j/\mu_j\longrightarrow\psi_{-H_{18};\theta_j}^{\perp}
 \quad\text{in }C^1\bigl(\Theta;\mathcal D^{1,2}(\R^{18}_+)\bigr),
\]
uniformly in $j$ as $M\to\infty$.
The direct quadratic term and this correction combine through
\eqref{eq:reduced-form} to produce $\mathcal C_{18}$.
Lemma~\ref{lem:single-bubble-localization} controls the remaining
terms in \eqref{eq:single-bubble-energy}. In particular,
\[
 \lim_{M\to\infty}\sup_{j\ge1}\|r_j\|_{C^1(\Theta)}=0,
 \qquad |\hat c_j|\le C\mu_j^2.
\]

Now sum the corrected bubbles and write their nonlinear boundary
interaction as
\[
 \widetilde W_N=\sum_{j=1}^N\widetilde U_j,\qquad
 E_N^\partial=-16\left(\widetilde W_N^{9/8}
                     -\sum_{j=1}^N\widetilde U_j^{9/8}\right).
\]

Note that $a_0\sum_{j=1}^\infty\lambda_j\rho_j^3\to0$ as $M\to\infty$.
The operator difference $\mathcal D_j=L_{g_0}-L_{g^{(j)}}$ is supported in
$\mathcal B_{2\rho_j}^+(p_j)$, where its second-, first-, and zeroth-order
coefficients $c_j^{ik},b_j^i,v_j$ satisfy, with $d_j(x)=\lambda_j+|x-p_j|$,
\begin{equation}\label{eq:patch-operator-envelope}
 |c_j^{ik}|+d_j|b_j^i|+d_j^2|v_j|
                              \le C(a_0\lambda_j)d_j^{3}.
\end{equation}
Disjointness of the patches $\mathcal B_{2\rho_j}^+(p_j)$ gives the exact identity
\[
 -L_{g_N}=-L_{g_0}+\sum_{j=1}^N\mathcal D_j,
\]
and $g_N=g^{(j)}$ on $\supp Z_{j\alpha}$.
Summing \eqref{eq:single-projected-equation} over $1\le j\le N$ yields
\begin{equation}\label{eq:exact-finite-residual}
 J_{g_N}'(\widetilde W_N)
 =\sum_{j,\alpha}\hat c_{j\alpha}
                Q_{g_N}(Z_{j\alpha},\cdot)
       +\mathcal R_N,
 \qquad \mathcal R_N=(E_N,E_N^\partial),\qquad
 E_N=\sum_{k\ne i}\mathcal D_k\widetilde U_i.
\end{equation}

The pointwise bound \eqref{eq:cross-patch-pointwise} and the
potential estimates of Lemma~\ref{lem:bridge-potentials} imply
\[
 \|\mathcal R_N\|_{\mathcal Y_N}
 \le C\left(\delta_M+a_0\sum_{j=1}^N\lambda_j\rho_j^3\right).
\]

Fix $C_0>0$ large and set
\[
 \eps_M=C_0\biggl(\delta_M+a_0\sum_{j=1}^\infty\lambda_j\rho_j^3
 +\sup_{j\ge1}\bigl\{R_j^{-4}+\mu_j^{1/8}
             +\mu_j(1+|\log\lambda_j|)^{5/9}\bigr\}\biggr)
 \longrightarrow0\quad\text{as }M\to\infty.
\]

The single-bubble convergence and the separation of the centers
verify the local convergence hypothesis of
Lemma~\ref{lem:relative-many-inverse} for $\widetilde W_N$.
The projected nonlinear equation for $w_N\in\mathcal X_N$ is
\begin{equation}\label{eq:many-projected-equation}
 \mathcal L_{g_N,\widetilde W_N}w_N
 =-\mathcal R_N+\bigl(0,N_{\widetilde W_N}(w_N)\bigr)
       +\sum_{j,\alpha}\eta_{j\alpha}
                       Q_{g_N}(Z_{j\alpha},\cdot).
\end{equation}
For $|w|,|\widetilde w|\le C\eps_MW_N$, we have
$|N_{\widetilde W_N}(w)-N_{\widetilde W_N}(\widetilde w)|
       \le C\eps_MW_N^{1/8}|w-\widetilde w|$,
$d_NW_N^{1/8}\le C$ and $G_{g,\Gamma}(W_N^{9/8})\le CW_N$.
Lemma~\ref{lem:relative-many-inverse} and the contraction mapping theorem
therefore yield a unique solution $w_N\in\mathcal X_N$ ($C^1$ in $\theta\in\Theta^N$) satisfying
\begin{equation}\label{eq:coarse-many-correction}
 |w_N|\le\eps_MW_N,\qquad
 \sup_{1\le j\le N}|\eta_j|\le\eps_M.
\end{equation}
Rescaled elliptic and Schauder estimates on balls of radius $cd_k(x)$ in
$\mathcal B_{2\rho_k}^+(p_k)$ also give
\begin{equation}\label{eq:many-spatial-derivatives}
 |w_N|+d_k|\nabla w_N|+d_k^2|D^2w_N|
                         \le C\eps_MW_N\quad\text{on }\mathcal B_{2\rho_k}^+(p_k).
\end{equation}

For $\theta\in\Theta^N$, write
\[
 u_N(\theta)=\widetilde W_N+w_N,\qquad
 \mathcal J_N(\theta)=J_{g_N}(u_N(\theta)).
\]

\begin{proof}[The metric case $n=18$ of Theorem~\ref{thm:sharpness}]
Fix $a_0>0$ sufficiently small. The preceding construction supplies
a $C^1$ family $u_N(\theta)$ for large $M$. We choose its parameters
so that all multipliers vanish.

\emph{Step 1. The reduced gradient.}
The interaction estimates of Lemma~\ref{lem:interaction-moments}
control the metric residual and the weighted trace of $w_N$.
Differentiating the orthogonality condition $Q_g(w_N,Z_{i\alpha})=0$
gives
\[
 Q_g(D_{\theta_j}w_N,Z_{i\alpha})
 =-\delta_{ij}Q_g(w_N,D_{\theta_j}Z_{j\alpha}).
\]
Together with \eqref{eq:exact-finite-residual} and
\eqref{eq:many-projected-equation}, this identity implies
\begin{equation}\label{eq:exact-force-identity}
\begin{aligned}
 &D_{\theta_j}\bigl(J_{g_N}(\widetilde W_N+w_N)-J_{g_N}(\widetilde W_N)\bigr)\\
 &\quad=\langle J_{g_N}'(\widetilde W_N+w_N)-J_{g_N}'(\widetilde W_N),
                       D_{\theta_j}\widetilde W_N\rangle
 +\langle J_{g_N}'(\widetilde W_N+w_N),D_{\theta_j}w_N\rangle\\
 &\quad=\langle D_{\theta_j}\mathcal R_N,w_N\rangle
       -\int_\Gamma N_{\widetilde W_N}(w_N)D_{\theta_j}\widetilde U_j\,\dd x'
 -\sum_\alpha\eta_{j\alpha}
                  Q_g(w_N,D_{\theta_j}Z_{j\alpha}).
\end{aligned}
\end{equation}
By \eqref{eq:bridge-moments}, \eqref{eq:cross-metric-moments} and $|w_N|\le\eps_MW_N$,
$|\langle D_{\theta_j}\mathcal R_N,w_N\rangle|\le C\eps_Me_j$.
By \eqref{eq:local-trace-moment},
\[
 \left|\int_\Gamma N_{\widetilde W_N}(w_N)D_{\theta_j}\widetilde U_j\,\dd x'\right|
 \le C\eps_M\int_\Gamma \widetilde W_N^{1/8}U_j|w_N|\,\dd x'
             \le C\eps_Me_j,
\]
and $|Q_g(w_N,D_{\theta_j}Z_{j\alpha})|\le C\eps_M$, so
$|\sum_\alpha\eta_{j\alpha}Q_g(w_N,D_{\theta_j}Z_{j\alpha})|\le C\eps_Me_j$.
Hence
\[
 \bigl|D_{\theta_j}\bigl(J_{g_N}(\widetilde W_N+w_N)-J_{g_N}(\widetilde W_N)\bigr)\bigr|
                                      \le C\eps_Me_j.
\]
Moreover, \eqref{eq:exact-finite-residual} gives
\[
 D_{\theta_j}J_{g_N}(\widetilde W_N)
       -D_{\theta_j}J_{g^{(j)}}(\widetilde U_j)
 =\sum_{i\ne j,\alpha}\hat c_{i\alpha}
                  Q_g(Z_{i\alpha},D_{\theta_j}\widetilde U_j)
                    +\langle\mathcal R_N,D_{\theta_j}\widetilde U_j\rangle,
\]
which is bounded by $C\sum_{i\ne j}\mu_i^2e_{ij}+Ce_j\le Ce_j$.
Combining this with \eqref{eq:single-bubble-energy} and
Lemma~\ref{lem:single-bubble-localization} yields
\begin{equation}\label{eq:many-bubble-force}
 |D_{\theta_j}\mathcal J_N(\theta)-\mu_j^2D\mathcal C_{18}(\theta_j)|
       \le \mu_j^2\eps_M+Ce_j.
\end{equation}

\emph{Step 2. Choice of parameters and the exact equation.}
For $0\le t\le1$, consider the homotopy with components
\[
 (\mathcal H_t(\theta))_j
 =(1-t)\mu_j^{-2}D_{\theta_j}\mathcal J_N(\theta)
                         +tD\mathcal C_{18}(\theta_j).
\]
By \eqref{eq:many-bubble-force} and \eqref{eq:interaction-small},
we may choose $M$ independently of $N$ so that
\[
 |(\mathcal H_t(\theta))_j|\ge\frac{m_*}{2}
 \quad\text{whenever }\theta_j\in\partial\Theta.
\]
Thus the homotopy has no zero on $\partial(\Theta^N)$.
Multiplication of each gradient block by $\mu_j^{-2}>0$ preserves
degree, and consequently
\[
 \deg(D\mathcal J_N,\Theta^N,0)
       =\deg(D\mathcal C_{18},\Theta,0)^N=1.
\]
For every $N$, there is a critical point in
$\operatorname{int}(\Theta^N)$.

Fix such a critical point $\theta$.
The projected equation gives
\[
 J_{g_N}'(\widetilde W_N+w_N)=\sum_{i,\alpha}c_{i\alpha}
 Q_{g_N}(Z_{i\alpha},\cdot),\qquad
 c_{i\alpha}=\hat c_{i\alpha}+\eta_{i\alpha}.
\]
Let $G_N$ be the block-diagonal flat Jacobi Gram matrix, and define $R_N$ by
\[
 (G_N+R_N)_{i\alpha,j\beta}
 =Q_g\!\left(Z_{i\alpha},
       \partial_{\theta_j^\beta}(\widetilde W_N+w_N)\right).
\]
Since $Q_g(Z_{i\alpha},D_{\theta_j}w_N)=-\delta_{ij}Q_g(D_{\theta_j}Z_{j\alpha},w_N)$,
we have $|(R_N)_{jj}|\le C(a_0+\eps_M)$ and $|(R_N)_{ij}|\le Ce_{ij}$ for $i\ne j$.
Thus, for small $a_0$ and large $M$,
\[
 \|G_N^{-1}R_N\|_{\ell^\infty\to\ell^\infty}
 \le C\left(a_0+\eps_M+\sup_i\sum_{j\ne i}e_{ij}\right)<1,
\]
so $G_N+R_N$ is invertible. Since $0=D\mathcal J_N(\theta)=(G_N+R_N)^tc$,
we obtain $c=0$. Finally, $|\psi_j|\le\kappa U_j$ and $|w_N|\le\eps_MW_N$ with
$\kappa+\eps_M<1/2$ yield
\[
 \tfrac12 W_N\le(1-\kappa-\eps_M)W_N\le u_N(\theta)=\widetilde W_N+w_N
                  \le(1+\kappa+\eps_M)W_N\le\tfrac32 W_N.
\]

Choose one critical point for each $N$ and denote the resulting
exact solution by $u_N$. The comparison constants are independent
of $N$.

\emph{Step 3. The singular limit.}

For every compact $K\subset\overline{\R^{18}_+}\setminus\{0\}$ and $k\ge0$,
\[
 \sum_{j\gg1}\|D^kU_j\|_{L^\infty(K)}
       \le C_{K,k}\sum_{j\gg1}\lambda_j^8<\infty,
\]
so Appendix~\ref{app:local-robin} gives $\sup_N\|u_N\|_{C^k(K)}<\infty$.
Passing to a diagonal subsequence yields a smooth solution $u_\infty>0$ for $g_\infty$
on $\overline{\R^{18}_+}\setminus\{0\}$.
Since $u_N(p_j)\ge c\lambda_j^{-8}$ for all $N\ge j$,
\[
 |p_j|^8u_\infty(p_j)
       \ge c\Bigl(\frac{|p_j|}{\lambda_j}\Bigr)^8\longrightarrow\infty
 \quad\text{as }j\to\infty.
\]
On the rays $\{-re_Z:0<r<r_0\}$ and $\{re_n:0<r<r_0\}$,
\[
 \sup_{0<r<r_0}
 \bigl(u_\infty(-re_Z)+u_\infty(re_n)\bigr)
       \le C\sum_{j=1}^\infty\lambda_j^8|p_j|^{-16}<\infty,
\]
so $\lim_{r\to0^+}r^8u_\infty(-re_Z)=0$ and $\int_0^{r_0}u_\infty(re_n)^{1/8}\,\dd r<\infty$.
Restricting to $\mathcal B_1^+$ completes the $n=18$ case, and hence the proof of
Theorem~\ref{thm:sharpness}.
\end{proof}

\appendix

\section{Local analytic estimates}
\subsection{Mixed Green kernels and positive gauges}\label{app:green-kernels}
For \(v|_{\Sigma_r}=0\), scaling Poincaré and trace inequalities give
\[
\int_{\mathcal B_r^+}v^2\le Cr^2\int_{\mathcal B_r^+}|\nabla v|^2,
\qquad \int_{B_r}v^2\le Cr\int_{\mathcal B_r^+}|\nabla v|^2.
\]
Consequently the mixed form satisfies
\begin{equation}\label{eq:small-mixed-coercivity}
    Q_g(v,v)\ge(c-Cr-Cr^2)\|\nabla v\|_2^2.
\end{equation}
All coefficient bounds below are uniform after pulling
back to one fixed coordinate neighborhood. Thus
\(g_\rho(x)=g(\rho x)\) on \(\mathcal B_R^+\), with \(\rho R\le\delta_0\),
has ellipticity constants independent of \(R\) and \(\rho\).
The constants below are allowed to depend on the fixed metric.

Let \(G_0\) be the Green kernel for \(-\Delta_{g_\rho}\) on a
half-ball or half-annulus, with Neumann condition on the flat
face and Dirichlet condition on the spherical boundaries.
Write \(C^{ij}=\sqrt{\det g_\rho}\,(g_\rho)^{ij}\). Reflect the
tangential block and the normal coefficient \(C^{nn}\) evenly, and
the normal-tangential coefficients \(C^{in}\) oddly. The Fermi
conditions make the reflected divergence-form coefficients
Lipschitz and uniformly elliptic. The reflected domain has
Dirichlet condition on its spherical boundaries. The Dirichlet
Green estimates \cite{LSW,GW} give
\begin{equation}\label{eq:mixed-green-upper}
    0<G_0(x,y)\le C|x-y|^{2-n}.
\end{equation}
For the two-sided estimates, take the half-ball \(\mathcal B_R^+\).
In the reflected ball put \(d_R(x)=\dist(x,\partial \mathcal B_R)\) and
\[
K_R(x,y)=|x-y|^{2-n}
\left(1\wedge\frac{d_R(x)}{|x-y|}\right)
\left(1\wedge\frac{d_R(y)}{|x-y|}\right).
\]
The interior Green lower bound of \cite{LSW,GW}, followed by Harnack
chains at distance comparable to \(|x-y|\) from the spherical boundary,
gives the lower estimate away from the spherical boundary. To reach
the spherical boundary \(\Sigma_R\), let \(\ell\) be the distance to
\(\Sigma_R\) in a patch of radius \(r\le R/8\) containing no pole.
The coefficients are Lipschitz, \(|\nabla\ell|=1\), and
\(|D^2\ell|\le C/R\). For fixed large \(K\), on a smaller patch one has
\[
-\Delta_{g_\rho}(\ell+K\ell^2/r)<0
<-\Delta_{g_\rho}(\ell-K\ell^2/r).
\]
Indeed, \(|\Delta_{g_\rho}\ell|\le C(R^{-1}+\rho)\) and
\(|\nabla\ell|_{g_\rho}^2\ge c\). On \(\ell\le r/(4K)\), the quadratic
term of size \(K/r\) dominates when \(K\) is large and
\(\rho R\le\delta_0\). Comparison, interior Harnack and scaling give
\[
\frac{w(x)}{w(x_r)}\asymp\frac{d_R(x)}r
\]
for a positive harmonic function vanishing on the spherical boundary.
Here \(x_r\) is an interior point at distance comparable to \(r\).
Apply this first in \(x\) and then in \(y\), using symmetry. Thus
\(G_0\asymp K_R\), with constants independent of \(\rho,R\) when
\(\rho R\le\delta_0\).

By even reflection of the Neumann condition, the half-ball kernel is
the sum of the reflected Dirichlet kernels with poles at \(y\) and its
reflection. For a face pole the two terms coincide. The
boundary-distance comparison gives
\begin{subequations}\label{eq:green-three}
    \begin{equation}
        G_0(0,y)\asymp|y|^{2-n}-R^{2-n}.
    \end{equation}
    For three points, it gives
    \begin{equation}
        \frac{G_0(x,z)G_0(z,y)}{G_0(x,y)}
        \le C\bigl(|x-z|^{2-n}+|z-y|^{2-n}\bigr).
    \end{equation}
\end{subequations}
The boundary pole is the limit of interior poles after reflection.
The second inequality follows from the two-sided Dirichlet estimate
by splitting according to which of \(|x-z|,|z-y|\) is at least
\(|x-y|/2\) and retaining the boundary-distance factors on that side.

The same estimates hold for the metric pair \((-L_{g_\rho},B_{g_\rho})\)
on a sufficiently small coordinate neighborhood before rescaling. In fact
\(|{\frac{n-2}{4(n-1)}}R_{g_\rho}|\le C\rho^2\) and
\(|{\frac{n-2}{2}}h_{g_\rho}|\le C\rho\).
By \eqref{eq:green-three},
\begin{align*}
    &\sup_{x,y}\frac{
        \int G_0(x,z)|{\frac{n-2}{4(n-1)}}R_{g_\rho}(z)|G_0(z,y)\dd z
        +\int_\Gamma G_0(x,z)|{\frac{n-2}{2}}h_{g_\rho}(z)|G_0(z,y)\dd z}
    {G_0(x,y)}\\
    &\qquad\le C\sup_x\left(
    \rho^2\int_{\mathcal B_R^+}|x-z|^{2-n}\dd z
    +\rho\int_{B_R}|x-z|^{2-n}\dd z\right)
    \le C((\rho R)^2+\rho R).
\end{align*}
Choose \(\delta_0\) so that the last bound is less than \(1/4\).
Iterating this bound in the signed volume/face resolvent series gives
\[
 |G-G_0|\le\sum_{m=1}^\infty4^{-m}G_0=\tfrac13G_0.
\]
Thus \(0<G\asymp G_0\). The same smallness makes the mixed form
coercive. This argument allows the scalar curvature and boundary
mean curvature potentials to have either sign.
On a half-annulus, coercivity again gives a positive mixed kernel.
Comparison on the annulus with the half-ball kernel gives
\eqref{eq:mixed-green-upper} for this metric kernel as well.

On the original half-ball, choose \(r\) so that
\eqref{eq:small-mixed-coercivity} is coercive. The Lax--Milgram
theorem gives a unique weak solution of
\[
\begin{cases}
    -L_g\varphi=0&\text{in }\mathcal B_r^+,\\
    B_g\varphi=0&\text{on }B_r,\\
    \varphi=1&\text{on }\Sigma_r.
\end{cases}
\]
Testing with \(\varphi_-\) and using coercivity gives \(\varphi_-=0\),
hence \(\varphi\ge0\). The strong maximum principle and boundary
point lemma give \(\varphi>0\). Write \(G[f]\) and
\(G_\Gamma[f^{\partial}]\) for the volume and face Green
potentials. The Green representation with the kernel just constructed gives
\begin{equation}\label{eq:positive-gauge-bound}
    \varphi-1=-G\left[\frac{n-2}{4(n-1)}R_g\right]
    -G_\Gamma\left[\frac{n-2}{2}h_g\right].
\end{equation}
The kernel upper bound gives
\(\int G\le Cr^2\) and \(\int G_\Gamma\le Cr\). Hence
\[
\|\varphi-1\|_\infty\le C(r^2\|R_g\|_\infty+r\|h_g\|_\infty).
\]
After a further shrinking, \(1/2\le\varphi\le3/2\).

\subsection{Local Robin and critical trace estimates}\label{app:local-robin}

Let $D'\Subset D\cup\Gamma_D$ be nested half-balls or cylinder strips,
where $\Gamma_D$ is the smooth boundary part. Their separation from
$\partial D\setminus\Gamma_D$ is fixed and positive, and $D$ avoids
the edge. Assume uniform ellipticity and fixed smooth coefficient
bounds. The constants below depend on these bounds and the separation.

For the linear problem
\[
 \begin{cases}
 Aw=f&\text{in }D,\\
 \partial_\nu w+\kappa w=f^{\partial}&\text{on }\Gamma_D,
 \end{cases}
\]
with smooth bounded $\kappa$, the local oblique estimates of
\cite{ADN,Lieberman} give, for every $D'\Subset D\cup\Gamma_D$,
\begin{equation}\label{eq:buffered-robin-estimate}
 \|w\|_{C^{k+2,\alpha}(D')}
 \le C\left(\|w\|_{L^2(D)}+\|f\|_{C^{k,\alpha}(D)}
                    +\|f^{\partial}\|_{C^{k+1,\alpha}(\Gamma_D)}\right).
\end{equation}
The obliqueness constant has a uniform positive lower bound.
For $q\ge1$, the local linear maximum estimate and interpolation
on nested patches also give
\begin{equation*}
 \|w\|_{L^\infty(D')}
 \le C_q\bigl(\|w\|_{L^q(D)}+\|f\|_{L^\infty(D)}
                         +\|f^{\partial}\|_{L^\infty(\Gamma_D)}\bigr).
\end{equation*}
This estimate applies to signed $w$. Together with
\eqref{eq:buffered-robin-estimate}, it gives local $C^1$ compactness
from local $L^1$ bounds and smooth bounded coefficients and data.
For positive solutions of the homogeneous equation with bounded
Robin coefficient, the boundary Harnack argument in
\cite[proof of Proposition 2.4]{LM} applies.
In particular, on local half-annuli,
\begin{equation}\label{eq:buffered-robin-harnack}
 \sup_{D'}w\le C\inf_{D'}w.
\end{equation}
The constants depend on a bound for the scaled Robin potential
and are uniform under multiplication of $w$ by a positive constant.

For a bounded positive semilinear solution, the initial boundary
H\"older estimate can be obtained from its local Neumann
representation. For a bounded face density $f$, split the kernel difference at
$|x-y|=2r$. The local kernel and outer gradient bounds give
\[
 |G_\Gamma f(x)-G_\Gamma f(x')|
 \le C\|f\|_\infty\left(r+r\int_{2r}^1\frac{\dd s}{s}\right)
 \le C\|f\|_\infty r(1+|\log r|)
\]
whenever $|x-x'|\le r$. Volume potentials satisfy a stronger bound.
Hence $u\in C^\alpha$
on smaller patches for every $\alpha<1$.
If $0<c\le u\le C$, then $u^{n/(n-2)}\in C^\alpha$.
The boundary Schauder estimate first gives $u\in C^{1,\alpha}$.
Iteration of \eqref{eq:buffered-robin-estimate} then gives uniform
$C^k$ bounds after Harnack normalization.

Let $u_j$ be weak solutions converging strongly in $H^1$ to
a smooth positive solution on a fixed patch. To obtain a uniform local
supremum bound, first consider a positive solution $u$ and test
with bounded truncations of
$\eta^2u^{2\beta-1}$, $\beta>1$. Ellipticity, Young's inequality,
and the trace Sobolev inequality give
\begin{equation}\label{eq:critical-robin-test}
 \begin{aligned}
 &\|\nabla(\eta u^\beta)\|_{L^2(D)}^2
       +\|\eta u^\beta\|_{L^{\frac{2(n-1)}{n-2}}(\Gamma)}^2
       +\|\eta u^\beta\|_{L^{\frac{2n}{n-2}}(D)}^2\\
 &\quad\le C\beta^2\left\{
 \int_D(|\nabla\eta|^2+\eta^2)u^{2\beta}
 +\int_\Gamma(1+u^{\frac{2}{n-2}})(\eta u^\beta)^2\right\}.
 \end{aligned}
\end{equation}
For the positive-part equation before positivity has been proved,
use the signed test $\eta^2u|u|^{2\beta-2}$, with bounded
truncations. The same inequality holds with $u$ replaced by
$|u|$, because the nonlinear boundary contribution is bounded
by $C|u|^{\frac{2}{n-2}}(\eta|u|^\beta)^2$.
Here $C$ depends on the bounded lower-order geometric coefficients.
The $L^2$ trace inequality gives
\[
 \int_\Gamma(\eta u^\beta)^2
 \le\varepsilon\|\nabla(\eta u^\beta)\|_2^2
       +C_\varepsilon\|\eta u^\beta\|_2^2.
\]
This absorbs the bounded face term.
For one fixed $1<\beta\le n/(n-2)$, H\"older gives
\[
 \int_\Gamma u^{\frac{2}{n-2}}(\eta u^\beta)^2
 \le\|u\|_{L^{\frac{2(n-1)}{n-2}}(\supp\eta\cap\Gamma)}^{\frac{2}{n-2}}
                         \|\eta u^\beta\|_{L^{\frac{2(n-1)}{n-2}}(\Gamma)}^2.
\]
Choose the half-ball so that
\[
 C\beta^2\|u_j\|_{L^{2(n-1)/(n-2)}(\supp\eta\cap\Gamma)}^{2/(n-2)}<\frac12.
\]
The trace theorem transfers strong $H^1$ convergence to the
critical trace space. Since the limit is smooth, shrinking the
patch makes the displayed norm small uniformly for all late $j$.
Absorption in \eqref{eq:critical-robin-test} gives
\[
 u_j\in L^{\frac{2(n-1)\beta}{n-2}}(\Gamma'),
 \qquad u_j^{\frac2{n-2}}\in L^{(n-1)\beta}(\Gamma').
\]
For a coefficient $c\in L^p(\Gamma')$, $p>n-1$, interpolation gives
\[
 \int_{\Gamma'}|c|w^2
 \le\varepsilon\|w\|_{L^{\frac{2(n-1)}{n-2}}(\Gamma')}^2
                     +C_{\varepsilon,p,\|c\|_p}\|w\|_{L^2(\Gamma')}^2.
\]
Apply this inequality in the subsequent power tests on nested
half-balls. The ordinary Moser iteration now gives a uniform
local supremum bound. The preceding H\"older estimate gives uniform
convergence. The positive limit then permits iteration of
\eqref{eq:buffered-robin-estimate}, giving smooth convergence
on smaller patches. The
critical coefficient is absorbed at the first, fixed power.
Subsequent steps use the resulting subcritical coefficient.

On the Euclidean cylindrical end, let $v_c+w$ be an exact
solution with $w\in X_\delta$. On translates of each fixed strip,
$w\to0$ in $H^1$. The absorption argument following
\eqref{eq:critical-robin-test} and H\"older compactness give $w\to0$
uniformly on smaller strips. Thus $v_c+w\ge\frac12\min v_c$ for large $T$.
The secant coefficient is then smooth with uniform bounds, and
\eqref{eq:buffered-robin-estimate} gives
\begin{equation}\label{eq:weighted-strip-regularity}
 \|w\|_{C^k([T,T+1]\times\overline S)}
 \le C_k\|w\|_{H^1([T-1,T+2]\times S)}
 \le C_ke^{-\delta T}\|w\|_{X_\delta}.
\end{equation}
The same estimate applies to differences of the exact solutions.

\section{Integration and projection formulas}\label{app:integration-formulas}

Write $\Gamma(z)$ and $\mathrm B(z,w)$ for Euler's gamma and beta functions.
Euler's beta integral \cite[Section 1.1]{AAR} evaluates the
direct metric terms. Gegenbauer projections below evaluate the
Robin correction and boundary shape terms.
For multi-indices $\alpha$ and $p>(n-1)/2+|\alpha|$, the tangential
moments used throughout the paper are
\begin{equation}\label{eq:beta-rule}
 \int_{\R^{n-1}}{x'}^{2\alpha}(|{x'}|^2+(1+{x_n})^2)^{-p}\dd {x'}
 =\frac{\prod_i\Gamma(\alpha_i+\tfrac12)
          \Gamma(p-\tfrac{n-1}{2}-|\alpha|)}{\Gamma(p)}
    (1+{x_n})^{n-1+2|\alpha|-2p}.
\end{equation}
The remaining normal integration uses
\[
 \int_0^\infty {x_n}^r(1+{x_n})^{-p}\dd {x_n}
 =\frac{\Gamma(r+1)\Gamma(p-r-1)}{\Gamma(p)}\quad\text{for }r>-1,\ p>r+1.
\]
Odd tangential moments vanish by reflection.

We use the Gegenbauer normalization and orthogonality of
\cite[Section 4.7]{Szego}. Write the factorials
\[
 (z)_k=\prod_{j=0}^{k-1}(z+j),\qquad
 k^{\underline r}=\prod_{j=0}^{r-1}(k-j),
\]
for $k,r\in\mathbb N_0$, with empty products equal to one.
For tangential harmonic degree $\ell$, put
$\lambda=\ell+(n-2)/2$.
For $\lambda+\tfrac12>\nu$, the projection and norm formulas are
\begin{equation}\label{eq:gegenbauer-beta}
 \begin{aligned}
 I_\nu(k)&=\int_{-1}^1(1-t^2)^{\lambda-1/2}
                 (1+t)^{-\nu}C_k^\lambda(t)\dd t\\
 &=(-1)^k2^{2\lambda-\nu}
       \mathrm B(\lambda+\tfrac12-\nu,\lambda+\tfrac12)
       \frac{(2\lambda)_k(\nu)_k}{k!(2\lambda-\nu+1)_k},\\
 h_k&=\mathrm B(1/2,\lambda+\tfrac12)
       \frac{(2\lambda)_k}{k!}\frac\lambda{k+\lambda}.
 \end{aligned}
\end{equation}
Rodrigues' formula \cite[equation (4.7.12)]{Szego}, followed by $k$
integrations by parts, reduces $I_\nu(k)$ to Euler's beta integral.
The endpoint terms vanish since $\lambda+\tfrac12>\nu$. The expression for $h_k$
is the norm in \cite[equation (4.7.15)]{Szego}.
For $r\in\mathbb N_0$, $b>0$ and $c>b+r+1$, Gauss summation
and its differentiated form \cite[Sections 2.1--2.2]{AAR}
give
\begin{equation}\label{eq:factorial-moment}
 \sum_{k=0}^\infty k^{\underline r}\frac{(b)_k}{(c)_k}
 =\frac{(c-1)r!(b)_r}
        {\prod_{j=0}^r(c-b-1-r+j)}.
\end{equation}

Taking $b=1$, $c=v+1$ and using
$\mathrm B(k+1,v)=v^{-1}(1)_k/(v+1)_k$, together with
$k^{\underline r}(k+\alpha)=k^{\underline{r+1}}+(r+\alpha)k^{\underline r}$,
also gives, for $v>r+2$,
\begin{equation*}
 \sum_{k=0}^\infty k^{\underline r}(k+\alpha)\mathrm B(k+1,v)
 =r!\{(r+\alpha)\mathrm B(r+1,v-r-1)
                    +(r+1)\mathrm B(r+2,v-r-2)\}.
\end{equation*}
These formulas are used for the metric Robin pairings and the
Euclidean shape energy. After cancellation of the common positive
factor, the entries below are rational functions of the dimension.
We use the recurrence $\mathrm B(x+1,y)=x\mathrm B(x,y)/(x+y)$
together with the displayed integral formulas.

\section{Polynomial families for the counterexamples}

We give the Robin pairings and sign calculations used in
Section~\ref{subsec:successive-examples}. The profile-selection
argument is Lemma~\ref{lem:profile-selection}.

\subsection{Metric polynomials}\label{app:metric-wells}

We evaluate the correction from Lemma~\ref{lem:jet-inverse} in the
inverse-metric convention \eqref{eq:inverse-profile-source}, retaining
its constant negative mode.
By \eqref{eq:corrector-covariance}, it suffices to compute at the
unit bubble after transforming the polynomial.
In this subsection write $D=|y'|^2+(1+y_n)^2$.
For a tangential homogeneous harmonic polynomial $P_\ell$ of
degree $\ell$, put
$\phi=P_\ell R(z,y_n)D^{-n/2}$, where $z=|y'|^2$. Differentiation gives
\[
 -\Delta\phi=P_\ell D^{-(n+2)/2}\mathcal T_\ell R,
\]
where
\[
 \begin{split}
 \mathcal T_\ell R={}&-D\{4zR_{zz}+[2(n-1)+4\ell]R_z+R_{y_ny_n}\}\\
 &+4nzR_z+2n(y_n+1)R_{y_n}+2n(\ell-1)R.
 \end{split}
\]
Seek the particular solution in
\[
 R\in\Span\{z^iy_n^j:i,j\ge0,\ 2i+j\le d_*+2-\ell\},
\]
where $d_*$ is the largest metric degree of the polynomial.

For a trace-free tangential polynomial $H$, the source polynomial is
\begin{equation}\label{eq:profile-source-polynomial}
 D^{(n+2)/2}F_{-H}
 =n(n-2)H_{ij}y_iy_j
       -(n-2)D(\partial_iH_{ij})y_j
       +{\frac{n-2}{4(n-1)}}D^2\partial_i\partial_jH_{ij}.
\end{equation}
For $H=f(z,y_n)A$, put $P_2=(y')^tAy'$.
The source polynomial then reduces to
\[
 P_2\{n(n-2)f-2(n-2)Df_z+4{\frac{n-2}{4(n-1)}}D^2f_{zz}\}.
\]
Angular components can be extracted without choosing a numerical
basis. For a homogeneous polynomial $P$ of degree $k$ in $n-1$
variables, its top harmonic component is
\[
 \sum_{j=0}^{\lfloor k/2\rfloor}
 \frac{(-1)^j\Gamma(k+(n-1)/2-j-1)}
      {4^jj!\Gamma(k+(n-1)/2-1)}|y'|^{2j}\Delta_{y'}^jP.
\]
Applying $\Delta_{y'}$ verifies harmonicity. Subtraction and
iteration give the complete Fischer decomposition. The
matrix entries of $\mathcal T_\ell$ are specified by
\begin{equation*}\begin{split}
 \mathcal T_\ell(z^iy_n^j)
 &={}-[4i(i-1)+(2(n-1)+4\ell)i]z^{i-1}y_n^jD\\
 &-j(j-1)z^iy_n^{j-2}D+4ni z^iy_n^j
       +2njz^i(y_n+1)y_n^{j-1}+2n(\ell-1)z^iy_n^j.
\end{split}\end{equation*}
Terms with a negative exponent and zero prefactor are omitted.
Equating monomial coefficients determines a particular
solution, with any homogeneous harmonic ambiguity set to zero.
A different such choice changes its Robin completion by the
opposite harmonic function and leaves the completed
pairing unchanged. For $d_*\le4$ and angular degree
$\ell=2$, this space has at most nine monomials.
For $\phi_i=P_\ell R_iD^{-n/2}$, the boundary operator
$\partial_\nu-nU^{2/(n-2)}$ cancels the derivative of $D^{-n/2}$.
Its residual is therefore
\[
 r_i^{\partial}
 =(\partial_\nu-nU^{2/(n-2)})\phi_i
 =P_\ell b_i(z)(1+z)^{-n/2},
 \qquad b_i(z)=-\partial_{y_n}R_i(z,0).
\]
For $H=y_nA$ and $H=y_n zA$ the respective particular solutions are
\[
 \begin{gathered}
 R_0=\frac{n-2}{4}(y_n-1),\\
 R_1=\frac{(n-5)y_n^3+(n-9)y_n^2+(4n-11)y_n z
                 -(n+3)y_n-(2n-3)z-n+1}{24}.
 \end{gathered}
\]
Their residuals have coefficients
$b_0=-(n-2)/4$ and $b_1=((n+3)-(4n-11)z)/24$.
All angular harmonic components below use this monomial rule. Their residual
polynomials $b_i$ have degree at most one.

Complete each particular solution by the harmonic function
with Robin datum $-r_i^\partial$.
Compute first for a fixed angular harmonic $P_\ell$, and sum
the orthogonal harmonic components. The angular norm is retained
explicitly below.
For a profile basis $(H_i)$, write $F_i=F_{-H_i}$ and
$\lambda=\ell+(n-2)/2$. If $t$ is the last coordinate on the
compactified sphere, define the residual projection coefficients by
\[
 a_{ik}=\int_{-1}^1(1-t^2)^{\lambda-1/2}(1+t)^{-\ell}
 b_i\!\left(\frac{1-t}{1+t}\right)C_k^\lambda(t)\dd t.
\]
The norms $h_k$ are those in \eqref{eq:gegenbauer-beta}.
Green's identity gives the completed correction pairing
\begin{equation}\label{eq:complete-robin-pairing}
 R_{ij}=\int F_i\phi_j-\int_\Gamma r_j^{\partial}\phi_i
    +2^{-n}\left(\int_{\Sph^{n-2}}P_\ell^2\right)
      \sum_{\substack{k\ge0\\\ell+k\ne1}}
                       \frac{a_{ik}a_{jk}}{h_k(\ell+k-1)}.
\end{equation}
For the harmonic correction, the Cayley change of variables preceding
\eqref{eq:ball-robin-form} gives ball Robin datum
$-D^{n/2}r_i^{\partial}/2$. Its contribution to the quadratic form has factor
\[
 2^{2-n}\left(\frac12\right)^2=2^{-n}.
\]
Green symmetry gives the first two terms in
\eqref{eq:complete-robin-pairing}.
Jacobi compatibility removes every term with $\ell+k=1$.
The constant term $\ell=k=0$ has $(\ell+k-1)^{-1}=-1$ and is retained.
Formula~\eqref{eq:factorial-moment} evaluates the remaining series,
and \eqref{eq:beta-rule} gives the direct entries. Thus
\[
 \mathcal C(H_i,H_j)=J_2(H_i,H_j)-\frac12R_{ij}.
\]

In the rest of this subsection, divide $\mathcal C$ and its derivatives
by the common positive factor
$\mathrm B((n-1)/2,(n-1)/2)\int_{\Sph^{n-2}}((y')^t{A}y')^2$.
We suppress this factor in those formulas. Multiplication restores
the main-text normalization and preserves every sign and critical point.

For the eighteenth-dimensional profile, the energy has scale
degree at most eight. Its scale second derivative and center
Hessian entries have degree at most six.
Use their Bernstein expansions on $I=[15/16,1]$.
For a polynomial $G$ of degree at most six,
its Bernstein coefficients are
\[
 \beta_j^I(G)=\sum_{r=0}^j
 \frac{\binom jr}{\binom6r}\frac{G^{(r)}(15/16)}{16^rr!}.
\]
They give the expansion
\[
 G(15/16+t/16)=\sum_{j=0}^6\beta_j^I(G)\binom6j t^j(1-t)^{6-j}.
\]

The basis functions are nonnegative and sum to one, so
$\min_j\beta_j^I(G)\le G\le\max_j\beta_j^I(G)$ on $I$.

For \eqref{eq:profile18}, reflection symmetry diagonalizes the center
Hessian. The two $A$ directions, the remaining $Y$ directions,
and the $Z$ direction correspond to $\xi_1,\xi_3,\xi_{17}$.
The only additional angular decomposition needed in
\eqref{eq:complete-robin-pairing} is
\[
 |Y|^2((y')^t{A}y')
       =\frac{20}{21}|y'|^2((y')^t{A}y')+P_4,
 \quad \Delta P_4=0,\quad
 \frac{\int P_4^2}{\int((y')^t{A}y')^2}=\frac{40}{10143}.
\]
Scale derivatives use angular degrees two and four. First center
variations also use degrees one and three, with their kernel
projection removed. Substituting the profile in
\eqref{eq:profile-source-polynomial} and then
\eqref{eq:complete-robin-pairing} gives polynomials of degree at most
six for all four Hessian entries. The moment formulas give the negative value
\begin{equation}\label{eq:18-profile-bounds}
 \mathcal C_{18}(1,0)<-\frac45.
\end{equation}
At the interval endpoints the scale derivatives satisfy
\[
 \partial_s\mathcal C_{18}(15/16,0)<-\frac1{25},\qquad
 \partial_s\mathcal C_{18}(1,0)>\frac1{16}.
\]
For $0\le j\le6$, the Bernstein coefficients satisfy
\begin{gather*}
 \beta_j^I(\partial_{ss}\mathcal C_{18})>\frac9{10},\\
 \beta_j^I(\partial_{\xi_i\xi_i}\mathcal C_{18})>\frac75
                 \quad\text{for }i=1,3,\\
 \beta_j^I(\partial_{\xi_{17}\xi_{17}}\mathcal C_{18})>\frac16.
\end{gather*}
For these bounds, insert the polynomial particular solution into
the displayed formula for $\beta_j^I$ and use
\eqref{eq:factorial-moment} for its harmonic completion.

For $n\ge19$, the positive Gram matrices are identified in the proof
of Lemma~\ref{lem:metric-wells}. We evaluate their dilation and
center-Hessian entries in the common positive normalization above.

At $n=19$ test the coefficient vector of the factorized profile
\[
 H=y_n(1-y_n)(3-|y'|^2)A,
\]
with the second copy zero. The complete pairing
\eqref{eq:complete-robin-pairing} gives the strict bound
\begin{equation}\label{eq:nineteen-negative-dilation}
 q(c)<-\frac16.
\end{equation}

For $n\ge20$, take the basis $(y_nA,y_n|y'|^2A)$ and its positive
Gram matrix $G$. Take $c_1=(3,-1)$ and $c_2=0$, corresponding to
$y_nA(3-|y'|^2)$.
Substitution of $3R_0-R_1$ in \eqref{eq:complete-robin-pairing} gives
\[
 \frac{q(c)}{G_{00}}
 =-\frac{8p(n)}{(n-8)(n-7)(n-6)(n-5)(n-3)(n-2)^2(n^2-2n-2)}.
\]
Writing $y=n-20\ge0$, we obtain the positive decomposition
\[
 \begin{aligned}
 p(n)={}&n^4y(6y^3+237y^2+3320y+15805)\\
 &+n^3(16404y+148177)+138541n^2\\
 &+50229y+904896>0.
 \end{aligned}
\]
Thus $q(c)<0$. Completion of the square gives
\begin{equation}\label{eq:cubic-metric-dilation}
 4G_{01}^2-3G_{00}G_{11}
 =\frac{G_{00}}2\left\{
 2G_{00}\left(3-\frac{2G_{01}}{G_{00}}\right)^2-q(c)\right\}>0.
\end{equation}

The mixed Gram entry used there is
\begin{equation}\label{eq:metric-positive-coupling}
 -\mathcal C(y_nA,y_n|y'|^2A)
 =\frac{(n-2)(n^2-2n-4)(n^3-7n^2+10n+2)}
               {64(n-6)(n-5)(n-4)(n-3)(n-1)}>0.
\end{equation}
Every factor is positive for $n\ge19$.

\subsection{Degree-six boundary polynomials}\label{app:euclidean-profile}

For $n\ge23$, use the Gram matrix from the proof of
Lemma~\ref{lem:euclidean-well}. We evaluate its entries using the
projection formulas in Appendix~\ref{app:integration-formulas}.
For tangential harmonic degree two, take $\lambda=(n+2)/2$ in
\eqref{eq:gegenbauer-beta}. The radial factor of
$\widehat{Q_1|z|^{2j}}$ is
\[
 (1-t)^j(1+t)^{-j-1}
 =\sum_{b=0}^j(-1)^b\binom jb2^{j-b}(1+t)^{-(j+1-b)}.
\]
Its Gegenbauer projection, divided by $I_1(k)$, is the polynomial
$P_j(k)$. Formula~\eqref{eq:gegenbauer-beta} gives
\[
 \begin{aligned}
 P_0(k)&=1,\\
 P_1(k)&=\frac{2(k+1)(k+n+1)}{n-1}-1,\\
 P_2(k)&=\frac{2(k+1)(k+2)(k+n)(k+n+1)}{(n-3)(n-1)}
       -\frac{4(k+1)(k+n+1)}{n-1}+1.
 \end{aligned}
\]
The full spherical degree of the $k$th mode is $k+2$.
Multiplying the two projections by its shape eigenvalue and
dividing by its Gegenbauer norm yields, for $0\le i,j\le2$,
\begin{equation}\label{eq:euclidean-gram-evaluator}
 \begin{split}
 \frac{G_{ij}}{\int_{\Sph^{n-2}}Q_1^2\,\dd\omega}
 ={}&\frac{2(n+2)\mathrm B(1/2,(n+3)/2)}{(n+1)^2}\\
 &\times\sum_{k=0}^\infty\frac{k!}{(n+2)_k}
 (k+\tfrac{n+2}{2})(k+1)(k+2)(k+n+1)P_i(k)P_j(k).
 \end{split}
\end{equation}
The polynomial in the summand has degree at most twelve.
Expanding it in falling factorials and applying
\eqref{eq:factorial-moment} with $b=1$ and $c=n+2$ evaluates
every entry as a rational function of $n$, up to the common
positive factor. The convergence condition is $n>12$.
The common factor cancels in the following ratios.

For every $n\ge23$, take
\[
 f(z)=Q_1(z)\left(1-\frac{|z|^2}{4}\right)^2,
 \qquad c_1=(1,-1/2,1/16),\quad c_2=0.
\]
Its dilation derivative is
\[
 q(c)=2G_{00}-4G_{01}+\tfrac34G_{02}+\tfrac32G_{11}
                    -\tfrac12G_{12}+\tfrac5{128}G_{22}.
\]
To verify its sign uniformly in $n$, put
\[
 D_n=(n-1)(n^2-2n-2)\prod_{j=5}^{12}(n-j)>0.
\]
Applying \eqref{eq:factorial-moment} to
\eqref{eq:euclidean-gram-evaluator} and collecting powers of $n-23$
gives the coefficient comparison
\[
 -D_n\left(1+128\frac{q(c)}{G_{00}}\right)
 =\sum_{j=0}^{11}a_j(n-23)^j,
 \qquad \min_{0\le j\le11}a_j=26.
\]
The coefficients $a_j$ are integers. Since $n-23\ge0$, the
right-hand side is strictly positive, and therefore
\begin{equation}\label{eq:euclidean-negative-dilation}
 \frac{q(c)}{G_{00}}<-\frac1{128}
 \qquad\text{for }n\ge23.
\end{equation}

The mixed Gram entry is
\begin{equation}\label{eq:euclidean-positive-coupling}
 \frac{G_{01}}{G_{00}}
 =\frac{(n-1)^4+(n-1)^3-19(n-1)^2-7(n-1)+48}{(n-5)(n-6)((n-1)^2-3)}>0.
\end{equation}
Its numerator is $(n-1)^2((n-1)^2-19)+(n-1)((n-1)^2-7)+48>0$.

\subsection{Polynomials with totally geodesic boundary}\label{app:geodesic-profile}

Use the tensors $W^\alpha$, $K_\alpha$, the contraction matrices
$\Omega_\alpha$ and the coefficient polynomials $P_{c_\alpha}$
from Section~\ref{subsec:successive-examples}.

The centered direct energy reduces to a radial integral.
After multiplying the two coefficient polynomials, define,
for $0\le k\le4$,
\[
 a_k(c_\alpha)=\sum_{i+j=k}(n+1+2(i+j)+2ij)c_{\alpha i}c_{\alpha j}.
\]
Here $0\le i,j\le2$.
The centered radial polynomial is
\[
 I_{c_\alpha}(x)=\sum_{k=0}^4\frac{a_k(c_\alpha)x^{k+2}}{n-6-2k}
                      \prod_{j=0}^k\frac{n-1+2j}{n-5-2j}.
\]

The source cancellation and orthogonality were proved in
Lemma~\ref{prop:geodesic-profile}.
The auxiliary functional in \cite[Section 4]{AlmarazBlowup} equals
$2\mathcal C$ with the normalization \eqref{eq:action}.
Hence \cite[Proposition 4.4, equations (4.4)--(4.5)]{AlmarazBlowup} gives
\[
 \mathcal C(H;s,0)=-k_n\sum_{\alpha=1}^2 I_{c_\alpha}(s^2),
 \qquad k_n>0.
\]

For every $n\ge35$, take
\[
 P_{c_1}(v)=\left(1-\frac v3\right)^2,\qquad P_{c_2}=0.
\]
Substituting $c_1=(1,-2/3,1/9)$ in the displayed formula for
$I_{c_1}$ and collecting powers of $n-35$ yields
\[
 I_{c_1}'(1)
 =-\frac{32(n-1)(n+1)}{27\prod_{j=5}^{14}(n-j)}
       \sum_{j=0}^7 b_j(n-35)^j,
 \qquad \min_{0\le j\le7}b_j=2.
\]
All coefficients $b_j$ are integers, and the denominator is
positive for $n\ge35$. Since $c_2=0$, this proves
\begin{equation}\label{eq:geodesic-negative-dilation}
 q(c)=2k_nI_{c_1}'(1)<0
 \qquad\text{for }n\ge35.
\end{equation}

In the following center-Hessian integrals, the coefficient polynomials
and their derivatives are evaluated at $r^2$. Define
\begin{align*}
 B_\alpha&=\int_0^\infty\!\!\int_0^\infty
 \frac{r^{n+4}(2P_{c_\alpha}P_{c_\alpha}'
                       +r^2(P_{c_\alpha}')^2)}
 {(r^2+(1+x_n)^2)^n}\,dr\,dx_n,\\
 C_\alpha&=\int_0^\infty\!\!\int_0^\infty
 \frac{r^{n+4}(P_{c_\alpha}')^2}
 {(r^2+(1+x_n)^2)^{n-1}}\,dr\,dx_n.
\end{align*}
By
\cite[Proposition 4.6, equation (4.6)]{AlmarazBlowup},
\begin{equation}\label{eq:geodesic-center-hessian}
 D^2_{\xi\xi}\mathcal C(H;1,0)
 =\frac{(n-2)^2|\Sph^{n-2}|}{(n-1)(n+1)(n+3)}
 \sum_{\alpha=1}^2
 \left\{-B_\alpha\Omega_\alpha+
 \left(-\frac{B_\alpha}4+
             \frac{n+3}{8(n-1)}C_\alpha\right)I\right\}.
\end{equation}
The matrices $\Omega_1,\Omega_2$ are diagonal in the fixed
basis chosen before Lemma~\ref{prop:geodesic-profile}, so
\eqref{eq:geodesic-center-hessian} is diagonal in that basis.
For the test
$P_{c_1}(v)=v$, $P_{c_2}=0$, radial integration by parts gives
\[
 B_1=\frac{3(n+5)}{2(n-1)}C_1>0,\qquad
 -\frac14+\frac{n+3}{8(n-1)}\frac{C_1}{B_1}
             =-\frac{n+6}{6(n+5)}<0.
\]
These identities show that the center Hessian is negative definite
for this test profile, as used in Lemma~\ref{prop:geodesic-profile}.

\section{Metric variations}\label{app:metric-variations}

Let $h$ be a smooth compactly supported Fermi-preserving tensor and
take $g_\eta=\exp(\eta h)$. Since $h_{in}=0$, the normal remains
$-e_n$, and the metric densities satisfy
\[
 \begin{aligned}
 g_\eta^{-1}&=I-\eta h+\tfrac12\eta^2h^2+O(\eta^3),\\
 \dd V_{g_\eta}&=e^{\eta\tr h/2}\,\dd x,\qquad
 \dd\sigma_{g_\eta}=e^{\eta\tr h/2}\,\dd x',\\
 h_{g_\eta}&=-\frac{\eta}{2(n-1)}\partial_{x_n}\tr h
                         \quad\text{on }\Gamma.
 \end{aligned}
\]
The linear scalar-curvature term is
\[
 R_{g_\eta}
 =\eta\bigl(\partial_i\partial_jh_{ij}-\Delta\tr h\bigr)+O(\eta^2).
\]
For a positive function $v$, substitution into \eqref{eq:action}
therefore gives
\begin{equation}\label{eq:full-first-action-variation}
 \begin{aligned}
 J_1(h,v)={}&\frac12\int_{\R^n_+}
       (\tfrac12\tr h\,\delta^{ij}-h_{ij})\partial_i v\partial_j v
 +\frac{n-2}{8(n-1)}\int_{\R^n_+}
       (\partial_i\partial_jh_{ij}-\Delta\tr h)v^2\\
 &-\frac{n-2}{8(n-1)}\int_\Gamma
       (\partial_{x_n}\tr h)v^2
 -\frac{(n-2)^2}{4(n-1)}\int_\Gamma
       (\tr h)v^{\frac{2(n-1)}{n-2}}.
 \end{aligned}
\end{equation}
Differentiation in $v$ at $U$ produces the sources in
\eqref{eq:metric-correction-source}. In the integration by parts,
the conormal and boundary-area terms cancel because
\[
 \begin{aligned}
 &\bigl((\tfrac12\tr h\,I-h)\nabla U\bigr)\cdot\nu
       -\frac{n-2}{2}(\tr h)U^{\frac n{n-2}}\\
 &\qquad=\frac12(\tr h)
       \bigl(\partial_\nu U-(n-2)U^{\frac n{n-2}}\bigr)=0.
 \end{aligned}
\]
Consequently,
\[
 D_fJ_1(h,U)[\zeta]
 =-\int_{\R^n_+}F_h\zeta-\int_\Gamma F_h^\partial\zeta.
\]

For the direct quadratic coefficient, suppose $\tr h=0$.
Then $\dd V_{g_\eta}=\dd x$, $\dd\sigma_{g_\eta}=\dd x'$ and
$h_{g_\eta}=0$. Expanding the Christoffel formula to second order gives
\[
 R_{g_\eta}
 =\eta\partial_i\partial_jh_{ij}+\eta^2R_2(h)+O(\eta^3),
\]
where
\[
 R_2(h)=-\tfrac14|\partial h|^2+\tfrac12|\diver h|^2
                    -\partial_i(h_{ik}\partial_jh_{jk}).
\]
The inverse-metric expansion and \eqref{eq:action} now imply
\[
 J_2(h)=\frac14\int_{\R^n_+}|h\nabla U|^2
       +\frac{n-2}{8(n-1)}\int_{\R^n_+}R_2(h)U^2.
\]
The flux of the last term in $R_2(h)$ vanishes on $\Gamma$ because
$h_{nk}=0$. Compact support therefore permits integration by parts,
which yields \eqref{eq:direct-form}.
For polynomial tensors of maximal degree $d$, the spherical flux
and cutoff errors are $O(R^{2d+2-n})$. Thus the same formula extends
to $d<(n-2)/2$ by first cutting off and then letting $R\to\infty$.

\section{Auxiliary estimates for the eighteenth-dimensional construction}
\label{app:many-bubble-estimates}

All quantities are those defined in
Section~\ref{subsec:eighteen-examples}, with $g=g_N$.
The constants below are independent of $N$ and
$\theta\in\Theta^N$.

\subsection{Green potential estimates}

If $p_k$ is a closest center to $x$, then for $i\ne k$,
$|x-p_i|\asymp |x-p_k|+d_{ik}$ with $d_{ik}\asymp 2^{-\min(i,k)}$.
Consequently, consecutive distances are comparable whenever both
indices differ from $k$, and
\[
 \frac{U_{i+1}(x)}{U_i(x)}
 \le C\Bigl(\frac{\lambda_{i+1}}{\lambda_i}\Bigr)^8
 \quad\text{for }i,i+1\ne k.
\]
Summing separately over $i<k$ and $i>k$, and retaining $U_k$, gives
$W_N\le C\max_{1\le j\le N}U_j$ on $\overline{\R^{18}_+}$.

For positive families $\widetilde U_i=\widetilde U_i(x;\theta_i)$ satisfying
$\widetilde U_i\asymp U_i$ and $|D_{\theta_i}\widetilde U_i|\le CU_i$, set
\[
 \widetilde W_N=\sum_{i=1}^N\widetilde U_i,\qquad
 E_N^\partial=-16\Bigl(\widetilde W_N^{9/8}
                              -\sum_{i=1}^N\widetilde U_i^{9/8}\Bigr)\le0.
\]

\begin{lemma}\label{lem:bridge-potentials}
The following estimates hold uniformly in $N\ge1$.
\begin{subequations}\label{eq:bridge-estimates}
\begin{align}
 G_{g,\Gamma}(U_j^{9/8})&\asymp U_j,\\
 G_{g,\Gamma}\Bigl[W_N^{9/8}
 \mathbf1_{\Gamma\setminus\bigcup_jB_{R\lambda_j}(p_j)}\Bigr]
       &\le CR^{-1}W_N,\\
 G_{g,\Gamma}\bigl[(U_i+U_j)^{9/8}-U_i^{9/8}-U_j^{9/8}\bigr]
 &\le C\frac{\sqrt{\lambda_i\lambda_j}}{d_{ij}}(U_i+U_j)\quad\text{for }i\ne j,\\
 \int_\Gamma U_i^{9/8}U_j\,\dd x'&\le Ce_{ij}\quad\text{for }i\ne j.
\end{align}
\end{subequations}
Moreover,
\begin{equation}\label{eq:bridge-moments}
\begin{gathered}
 \sup_{\R^{18}_+}\frac{G_{g,\Gamma}[-E_N^\partial]}{W_N}
       +\sup_\Gamma\frac{d_N(-E_N^\partial)}{W_N}\le C\delta_M,\\
 \int_\Gamma U_j(-E_N^\partial)\,\dd x'+
       \int_\Gamma W_N|D_{\theta_j}E_N^\partial|\,\dd x'
                         \le Ce_j,\\
 \int_\Gamma
   \bigl(\widetilde W_N^{1/8}-\widetilde U_j^{1/8}\bigr)U_jW_N\,\dd x'
                         \le Ce_j.
\end{gathered}
\end{equation}
\end{lemma}

\begin{proof}
Since $g_{in}=0$ ($i<18$) and $g_{nn}=1$, the even reflection
$\widetilde g_{ij}(x',x_n)=g_{ij}(x',|x_n|)$ converts the weak Neumann problem
on $\R^{18}_+$ into $-\diver(\widetilde g^{-1}\nabla u)+\frac{4}{17}\widetilde R_{g_N}u=2f^\partial\delta_{\{x_n=0\}}$ on $\R^{18}$.
Differentiating \eqref{eq:metric18-blocks} twice gives
\[
 |D^2f_0(y)|\le Ca_0,
 \qquad
 |D^2f_j(y)|\le Ca_0\lambda_j(\lambda_j+|y-p_j|)
 \quad\text{for }y\in\supp f_j,
\]
and hence
\[
 \sup_{x\in\R^{18}}\int_{\R^{18}}|x-y|^{-16}\Bigl|\frac{4}{17}R_{g_N}(y)\Bigr|\dd y
 \le Ca_0\left(1+\sum_{j=1}^N
              \lambda_j(\lambda_j+\rho_j)\rho_j^2\right)
 \le Ca_0\left(1+\sum_{j=1}^\infty\rho_j^4\right)\le Ca_0.
\]
For small $a_0$, the resolvent series and the two-sided Green bound of \cite{LSW,GW}
imply $G_g(x,y)\asymp|x-y|^{-16}$, which yields the first line of \eqref{eq:bridge-estimates}.

For a single bubble, rescaling by $x=p_j+\lambda_jX$ reduces the second bound in
\eqref{eq:bridge-estimates} to
\[
 \int_{|Y|>R}|X-Y|^{-16}
 (1+|Y|^2)^{-9}\,\dd Y
 \le CR^{-1}(1+|X|)^{-16}\qquad\text{with }Y\in\R^{17},
\]
which follows by splitting into $|Y|<|X|/2$, $|Y-X|<|X|/2$ and their complement.
Since $W_N\le C\max_jU_j$ implies $W_N^{9/8}\le C\sum_jU_j^{9/8}$, summing over $j$
proves the second bound in \eqref{eq:bridge-estimates}.

For the pair bounds with $i\ne j$, uniform equivalence on $\Theta$ allows us to set
$(s_i,\xi_i)=(s_j,\xi_j)=(1,0)$.
A Kelvin inversion centered at a boundary endpoint of the semicircle orthogonal to
$\Gamma$ through $p_i+\lambda_ie_n$ and $p_j+\lambda_je_n$, followed by a translation
and dilation, maps $U_i,U_j$ to concentric bubbles $U_{\delta,0},U_{1,0}$, where
\[
 \delta+\delta^{-1}
 =\frac{d_{ij}^2+\lambda_i^2+\lambda_j^2}{\lambda_i\lambda_j},
 \qquad 0<\delta\le1,\qquad
 \delta\asymp\frac{\lambda_i\lambda_j}{d_{ij}^2}.
\]
Conformal covariance of the flat Riesz potential and of the boundary integral
$\int_{\R^{17}}U_i^{9/8}U_j\,\dd x'$ (since $9+8=17$) reduces the third and fourth
bounds in \eqref{eq:bridge-estimates} to the concentric pair $(U_{\delta,0},U_{1,0})$.
On $\Gamma$, with $r=|x'|$,
\[
 U_{\delta,0}^{1/8}=\frac{\delta}{r^2+\delta^2},
 \qquad
 U_{1,0}^{1/8}=\frac{1}{1+r^2},
\]
and
\[
 (U_{\delta,0}+U_{1,0})^{9/8}
 -U_{\delta,0}^{9/8}-U_{1,0}^{9/8}
 \le C\min\bigl(U_{\delta,0}^{1/8}U_{1,0},
 U_{\delta,0}U_{1,0}^{1/8}\bigr)\le M_\delta(r),
\]
where
\[
 M_\delta(r)=C\begin{cases}
 \delta^{-1},&0<r<\delta,\\
 \delta r^{-2},&\delta\le r<\sqrt\delta,\\
 \delta^8r^{-16},&\sqrt\delta\le r<1,\\
 \delta^8r^{-18},&r\ge1.
 \end{cases}
\]
At $x'=Re_1\in\Gamma$, splitting $\R^{17}$ into $\{|y|<R/2\}$, $\{|y-Re_1|<R/2\}$ and
their complement gives
\[
 \begin{aligned}
 \int_{\R^{17}}|Re_1-y|^{-16}M_\delta(|y|)\,\dd y
 &\le C\left(R^{-16}\int_0^{R/2}M_\delta(r)r^{16}\,\dd r
       +R M_\delta(R)+\int_{R/2}^\infty M_\delta(r)\,\dd r\right)\\
 &\le C\sqrt\delta\,\bigl(U_{\delta,0}(Re_1,0)+U_{1,0}(Re_1,0)\bigr).
 \end{aligned}
\]
Harmonic extension extends this bound to $\R^{18}_+$, and direct radial integration gives
\[
 \int_{\R^{17}}M_\delta(|y|)\bigl(U_{\delta,0}(y,0)+U_{1,0}(y,0)\bigr)\dd y
 +\int_{\R^{17}}U_{\delta,0}^{9/8}U_{1,0}\,\dd y\le C\delta^8,
\]
proving the third and fourth lines of \eqref{eq:bridge-estimates}.

For the pointwise bound in \eqref{eq:bridge-moments}, since
$d_N\le\min(\lambda_i+|x-p_i|,\lambda_j+|x-p_j|)$ and
$\max(\lambda_i+|x-p_i|,\lambda_j+|x-p_j|)\ge d_{ij}/2$,
\[
 d_N\bigl[(U_i+U_j)^{9/8}-U_i^{9/8}-U_j^{9/8}\bigr]
 \le Cd_N(U_iU_j)^{1/16}(U_i+U_j)
 \le C\frac{\sqrt{\lambda_i\lambda_j}}{d_{ij}}(U_i+U_j).
\]
Subadditivity of $s\mapsto s^{1/8}$ gives
\[
 -E_N^\partial\le C\sum_{i<j}
 \bigl[(\widetilde U_i+\widetilde U_j)^{9/8}
                  -\widetilde U_i^{9/8}-\widetilde U_j^{9/8}\bigr],
\]
so summing over $i\ne j$ proves the first line of \eqref{eq:bridge-moments}.
Finally, $|D_{\theta_j}E_N^\partial|\le CU_j(\widetilde W_N^{1/8}-\widetilde U_j^{1/8})$.
Partition $\Gamma=\bigcup_{k=1}^N\Omega_k$ by
\[
 \Omega_k=\bigl\{x\in\Gamma:\widetilde U_k(x)=\max_{1\le i\le N}\widetilde U_i(x)\bigr\}.
\]
On $\Omega_j$, $\widetilde W_N\asymp U_j$ and concavity of $s^{1/8}$ give
\[
 \widetilde W_N^{1/8}-\widetilde U_j^{1/8}
        \le CU_j^{-7/8}\sum_{i\ne j}U_i,\qquad
 -E_N^\partial\le CU_j^{1/8}\sum_{i\ne j}U_i,
\]
so on $\Omega_j$,
\[
 U_j(-E_N^\partial)+W_N|D_{\theta_j}E_N^\partial|
 +\bigl(\widetilde W_N^{1/8}-\widetilde U_j^{1/8}\bigr)U_jW_N
 \le C\sum_{i\ne j}U_j^{9/8}U_i.
\]
On $\Omega_k$ with $k\ne j$, $\widetilde W_N\asymp U_k$, so the same three quantities
are bounded by $CU_jU_k^{9/8}$. Integrating over $\Gamma$ and applying the fourth line
of \eqref{eq:bridge-estimates} yields the remaining two bounds of \eqref{eq:bridge-moments}.
\end{proof}

\subsection{The single-bubble local error}

Let $r_j$ be the energy remainder in
\eqref{eq:single-bubble-energy}, and $R_j=\rho_j/\lambda_j$.

\begin{lemma}\label{lem:single-bubble-localization}
Uniformly for $\theta_j\in\Theta$ and $j\ge1$,
\begin{equation}\label{eq:18-local-error}
 \norm{r_j}_{C^1(\Theta)}
 \le C\bigl(R_j^{-4}+\mu_j^{1/8}+\mu_j(1+|\log\lambda_j|)^{5/9}\bigr),
 \qquad |\hat c_j|\le C\mu_j^2.
\end{equation}
In particular, $\lim_{M\to\infty}\sup_{j\ge1}\norm{r_j}_{C^1(\Theta)}=0$.
\end{lemma}
\begin{proof}
Fix $j$, write $\theta_j=(s,\xi)$, and set $U=U_{s,\xi}$,
$Q=Q_{s,\xi}$ and $\psi=\psi_{-H_{18};s,\xi}^{\perp}$.
Rescale by $x=p_j+\lambda_jy$ ($r=|y|$) and define
\[
 \widehat\psi_j(y)=\lambda_j^8\psi_j(p_j+\lambda_jy),\qquad
 \widehat g(y)=g^{(j)}(p_j+\lambda_jy),\qquad
 \widehat g^{-1}=\exp(fA)=\diag(e^f,e^{-f},1,\ldots,1).
\]
Then $f=\mu_jh_{18}$ on $\{r<R_j\}$, $\|f\|_\infty\le Ca_0$,
$\supp f\subset\mathcal B_{C/\lambda_j}$, and
\[
 |D^kf(y)|\le C_k\mu_j(1+|y|)^{4-k}.
\]
The common fourth-order metric term is included in $\widehat g$
throughout this rescaling.
For $\widehat g^{-1}=\diag(e^f,e^{-f},1,\ldots,1)$, the nonzero Christoffel symbols
involving $\beta\in\{3,\ldots,18\}$ are
\[
 \Gamma^1_{1\beta}=-\tfrac12f_\beta,\quad
 \Gamma^2_{2\beta}=\tfrac12f_\beta,\quad
 \Gamma^\beta_{11}=\tfrac12e^{-f}f_\beta,\quad
 \Gamma^\beta_{22}=-\tfrac12e^ff_\beta,
\]
and the scalar curvature is
\begin{equation}\label{eq:exact-diagonal-curvature}
 R_{\widehat g}=-\partial_{11}e^f-\partial_{22}e^{-f}
                  -\frac12\sum_{\beta=3}^{18}f_\beta^2.
\end{equation}
The determinant-one block in the first two coordinates contributes
the first two terms. For each $\beta\ge3$, the derivative terms
cancel in the scalar contraction, leaving $-f_\beta^2/2$.
Since $(\widehat g)_T|_\Gamma=I$ and $h_{\widehat g}=0$, integrating the
$\partial_{11}e^f$ and $\partial_{22}e^{-f}$ terms by parts twice and using,
for $D=|y'-\xi|^2+(y_n+s)^2$ and $\beta=1,2$,
\[
 \tfrac12(\partial_\beta U)^2-\frac2{17}
                     \partial_{\beta\beta}(U^2)
             =\frac{64}{17}s^{16}D^{-17},
\]
we obtain the exact identity
\begin{equation}\label{eq:exact-diagonal-action}
 J_{\widehat g}(U)-J_\delta(U)
 =\frac{128}{17}\int_{\R^{18}_+}(\cosh f-1)s^{16}D^{-17}\,\dd y
       -\frac1{17}\int_{\R^{18}_+}|\nabla_{3,\ldots,18}f|^2U^2\,\dd y.
\end{equation}

Since $B_{\widehat g}U=16U^{9/8}$ on $\Gamma$, the rescaled residual at $U$ is
$(-L_{\widehat g}U,0)$. On $\{r<R_j\}$, the linear part of $-L_{\widehat g}U$ is
$-\mu_jF_{-H_{18};s,\xi}$. Since $|F_{-H_{18};s,\xi}|=O((1+r)^{-14})$ and
the quadratic part of $L_{\widehat g}U$ is $O(\mu_j^2(1+r)^{-10}\mathbf1_{\{r<C/\lambda_j\}})$,
their $L^{9/5}(\R^{18}_+)$ norms satisfy
\[
 \|F_{-H_{18};s,\xi}\|_{L^{9/5}(r>R_j)}
 \le C\left(\int_{R_j}^\infty r^{17-14(9/5)}\,\dd r\right)^{5/9}
                            =CR_j^{-4},
\]
and
\[
 \norm{-\mu_j^{-1}L_{\widehat g}U+F_{-H_{18};s,\xi}}
              _{C^1_\theta(\mathcal D^{1,2})^*}
              \le C\bigl(R_j^{-4}+\mu_j(1+|\log\lambda_j|)^{5/9}\bigr).
\]
The exponent in the quadratic-source integral is
$17-10(9/5)=-1$, which accounts for the logarithmic factor.
By \eqref{eq:correction-decay}, $|D_y^kD_\theta^\beta\psi|\le C(1+r)^{-12-k}$
for $k\le2$ and $|\beta|\le1$, so
\[
 |(\Delta-L_{\widehat g})D_\theta^\beta\psi|
             \le C\mu_j(1+r)^{-10}\mathbf1_{\{r<C/\lambda_j\}}.
\]
Taking the energy-dual norm gives
\[
 \|(\Delta-L_{\widehat g})\psi\|_{C^1_\theta(\mathcal D^{1,2})^*}
             \le C\mu_j(1+|\log\lambda_j|)^{5/9}.
\]
For the nonlinear boundary term, $|D_\theta U|\le CU$ and
$\|\widehat\psi_j\|_{\mathcal D^{1,2}}+\|D_\theta\widehat\psi_j\|_{\mathcal D^{1,2}}\le C\mu_j$ give
\[
 |D_\theta N_U(\widehat\psi_j)|
 \le C\bigl(|\widehat\psi_j|^{1/8}|D_\theta\widehat\psi_j|
                               +|\widehat\psi_j|^{9/8}\bigr).
\]
The critical trace inequality therefore implies
\[
 \mu_j^{-1}\|N_U(\widehat\psi_j)\|_{C^1_\theta(\mathcal D^{1,2})^*}
 \le C\mu_j^{1/8}.
\]
Before applying the inverse, project the rescaling of $\psi$ onto
$\mathcal X^{(j)}$. Its flat Dirichlet orthogonality, the decay of
$\psi$ and the Jacobi fields, and the uniformly invertible Gram
matrix bound the change in the normalized correction by
$O_{C^1_\theta\mathcal D^{1,2}}(\mu_j+R_j^{-8})$.
Dividing \eqref{eq:single-projected-equation} by $\mu_j$ and subtracting
the projected flat correction equation now gives an equation on the
same orthogonal complement. The uniformly bounded single-bubble
energy inverse and the preceding residual estimates imply
\begin{equation}\label{eq:endpoint-response-convergence}
 \norm{\widehat\psi_j/\mu_j-\psi}_{C^1_\theta\mathcal D^{1,2}}
 \le C\bigl(R_j^{-4}+\mu_j^{1/8}+\mu_j(1+|\log\lambda_j|)^{5/9}\bigr).
\end{equation}

In \eqref{eq:exact-diagonal-action}, the quadratic contribution from
$\{r>R_j\}$ is $O(\mu_j^2R_j^{-8})$
and the quartic remainder $\cosh f-1-\frac12f^2=O(f^4)$ contributes at most
\[
 C\mu_j^4\int_1^{C/\lambda_j}r^{-1}\,\dd r
                              \le C\mu_j^4(1+|\log\lambda_j|).
\]
Taylor expansion of $J_{\widehat g}(U+\widehat\psi_j)$ gives
\begin{equation*}
 \begin{split}
 \mu_j^{-2}\bigl(J_{\widehat g}(U+\widehat\psi_j)-J_\delta(U)\bigr)
 &=\mu_j^{-2}\bigl(J_{\widehat g}(U)-J_\delta(U)\bigr)
       +\mu_j^{-2}\langle-L_{\widehat g}U,\widehat\psi_j\rangle\\
 &\quad+\frac1{2\mu_j^2}Q(\widehat\psi_j)
       +\frac1{2\mu_j^2}(Q_{\widehat g}-Q_\delta)
                              (\widehat\psi_j,\widehat\psi_j)+O_{C^1}(\mu_j^{1/8}).
 \end{split}
\end{equation*}
On the fixed correction $\psi$, the metric form difference satisfies
\[
 \|(Q_{\widehat g}-Q_\delta)(\psi,\psi)\|_{C^1(\Theta)}
 \le C\bigl(\mu_j+\mu_j^2(1+|\log\lambda_j|)\bigr).
\]
The linear metric contribution has radial integral
$\mu_j\int_1^\infty r^{-5}\,\dd r$, while the quadratic contribution
has radial integral $\mu_j^2\int_1^{C/\lambda_j}r^{-1}\,\dd r$.
The form difference is uniformly bounded on the energy space.
Together with \eqref{eq:endpoint-response-convergence}, this bounds
the $C^1$ error from replacing $\widehat\psi_j/\mu_j$ by $\psi$
in the quadratic terms. The residual estimate replaces the linear
source by $-F_{-H_{18};s,\xi}$. Since
$Q(\psi)=\int F_{-H_{18};s,\xi}\psi$, the linear and flat quadratic
terms sum to $-\tfrac12Q(\psi)$. Thus \eqref{eq:reduced-form} yields
\[
 J_2(-H_{18};U)-\tfrac12Q(\psi)=\mathcal C_{18}(s,\xi)
\]
up to the error in \eqref{eq:18-local-error}. Under the prescribed
scales, the normalized direct errors
$R_j^{-8}+\mu_j^2(1+|\log\lambda_j|)$ are absorbed by the right-hand
side of that estimate.

Finally, testing \eqref{eq:single-projected-equation} with $D_{\theta_j}\widetilde U_j$ gives
\[
 \sum_\alpha \hat c_{j\alpha}
       Q_{g^{(j)}}(Z_{j\alpha},D_{\theta_j}\widetilde U_j)
 =D_{\theta_j}J_{g^{(j)}}(\widetilde U_j)
 =\mu_j^2\bigl(D\mathcal C_{18}(\theta_j)+Dr_j(\theta_j)\bigr)
 =O(\mu_j^2).
\]
The correction derivative is $O(\mu_j)$ in energy, and the cutoff
errors vanish as $R_j\to\infty$. Consequently, the matrix
$Q_{g^{(j)}}(Z_{j\alpha},D_{\theta_j}\widetilde U_j)$ converges
uniformly on $\Theta$ to the positive flat Jacobi Gram matrix.
Its uniformly bounded inverse therefore gives
$|\hat c_j|\le C\mu_j^2$.
\end{proof}

\subsection{Interaction moment estimates}

For $i\ne k$, \eqref{eq:patch-operator-envelope} and $|\widetilde U_i|\le CU_i$ give
\begin{equation}\label{eq:cross-patch-pointwise}
 |\mathcal D_k\widetilde U_i|+
       |\mathcal D_kD_{\theta_i}\widetilde U_i|
 \le C(a_0\lambda_k)d_k
             \lambda_i^{8}d_{ik}^{-16}\mathbf1_{\mathcal B_{2\rho_k}(p_k)}.
\end{equation}

Integrating $d_k$ against $|x-y|^{-16}$ on $\mathcal B_{2\rho_k}(p_k)$ and using
$\rho_k\le cd_{ik}$ gives $G_g(|\mathcal D_k\widetilde U_i|)\le C(a_0\lambda_k)\rho_k^3U_i$.

\begin{lemma}\label{lem:interaction-moments}
For the projected correction $w_N$ and multipliers $\eta_j$ in
\eqref{eq:many-projected-equation}, the following estimates hold
uniformly for $1\le j\le N$.
\begin{align}
 \int_{\R^{18}_+}U_j|E_N|\,\dd x
 +\int_{\R^{18}_+}W_N|D_{\theta_j}E_N|\,\dd x
 &\le Ce_j,\label{eq:cross-metric-moments}\\
 \int_\Gamma U_j^{9/8}|w_N|\,\dd x'+|\eta_j|
 &\le Ce_j.\label{eq:local-trace-moment}
\end{align}
\end{lemma}

\begin{proof}
\emph{Step 1. Mixed metric integrals.}
In $E_N=\sum_{k\ne i}\mathcal D_k\widetilde U_i$, fix $j$ and estimate
$\int U_j|\mathcal D_k\widetilde U_i|$ in the three cases $k=j$, $i=j\ne k$,
and $i,j,k$ distinct.
For $k=j\ne i$, \eqref{eq:cross-patch-pointwise} and
\[
 \int_{\mathcal B_{2\rho_j}^+(p_j)}U_jd_j\,\dd x
 \le C\lambda_j^8\int_0^{2\rho_j}\frac{r^{17}}{(\lambda_j+r)^{15}}\,\dd r
 \le C\lambda_j^8\rho_j^3
\]
give $\int U_j|\mathcal D_j\widetilde U_i|\le C(a_0\lambda_j)\rho_j^3e_{ij}$.
For $i=j\ne k$, \eqref{eq:cross-patch-pointwise} gives
\begin{equation*}
 \sum_{k\ne j}\int U_j|\mathcal D_k\widetilde U_j|\,\dd x
 \le C\lambda_j^{16}\sum_{k\ne j}
       (a_0\lambda_k)\rho_k^{19}d_{jk}^{-32}
                  \le C\lambda_j^{16}\rho_j^{-16}\le Ce_j,
\end{equation*}
since $d_{jk}\asymp 2^{-k}$ for $k<j$ and $d_{jk}\asymp 2^{-j}$ for $k>j$ imply
\[
\begin{gathered}
     \sum_{k<j}(a_0\lambda_k)\rho_k^{19}d_{jk}^{-32}
    \le C\sum_{k<j}\lambda_k2^{13k}\le C\le C\rho_j^{-16},\\
    \sum_{k>j}(a_0\lambda_k)\rho_k^{19}d_{jk}^{-32}
    \le C2^{32j}\sum_{k>j}2^{-19k}
    \le C2^{13j}\le C\rho_j^{-16}.
\end{gathered}
\]
For $i,j,k$ distinct, $\rho_k\le c\min(d_{ik},d_{jk})$ implies $\rho_kd_{ij}\le Cd_{ik}d_{jk}$, so
\[
 \int U_j|\mathcal D_k\widetilde U_i|\,\dd x
 \le C(a_0\lambda_k)\rho_k^{19}\lambda_i^8\lambda_j^8
          d_{ik}^{-16}d_{jk}^{-16}
       \le C(a_0\lambda_k)\rho_k^3e_{ij}.
\]
Summing over $i\ne k$ and repeating the bound for
$D_{\theta_j}\widetilde U_j$ proves \eqref{eq:cross-metric-moments}.
Similarly, \eqref{eq:patch-operator-envelope} and \eqref{eq:many-spatial-derivatives}
give $|\mathcal D_kw_N|\le C\eps_M(a_0\lambda_k)d_kW_N$, and hence
\[
 \int_{\R^{18}_+} U_j\left|\sum_{k\ne j}\mathcal D_kw_N\right|\dd x\le C\eps_Me_j.
\]
For $i\ne j$, the radial cutoff defining $Z_{i\alpha}$ has zero
normal derivative on $\Gamma$. The Jacobi boundary condition therefore implies
$|B_{g_N}Z_{i\alpha}|\le CU_i^{9/8}$, and
\[
 \int_{\R^{18}_+} U_j|L_{g_N}Z_{i\alpha}|\,\dd x
 +\int_\Gamma U_j|B_{g_N}Z_{i\alpha}|\,\dd x'\le Ce_{ij}\quad\text{for }i\ne j,
\]
along with the same bound for $D_{\theta_i}Z_{i\alpha}$.

\emph{Step 2. A single-bubble adjoint estimate.}
Let $\mathcal L_j=\mathcal L_{g^{(j)},\widetilde U_j}$. We claim that if
$w\in\mathcal X^{(j)}$ solves
\[
 \mathcal L_jw=(f,f^{\partial})+
      \sum_\alpha\eta_\alpha
                     Q_{g^{(j)}}(Z_{j\alpha},\cdot),
\]
then
\[
 \int_\Gamma U_j^{9/8}|w|\,\dd x'+|\eta|
       \le C\left(\int_{\R^{18}_+} U_j|f|\,\dd x+\int_\Gamma U_j|f^{\partial}|\,\dd x'\right).
\]
Indeed, for smooth $\sigma_\eps\to\operatorname{sign}w$ on $\Gamma$ with $|\sigma_\eps|\le1$,
Lemma~\ref{lem:relative-many-inverse} with $N=1$ provides $z_\eps\in\mathcal X^{(j)}$ solving
\[
 \mathcal L_jz_\eps=(0,U_j^{9/8}\sigma_\eps)+\sum_\alpha\zeta_{\eps,\alpha}Q_{g^{(j)}}(Z_{j\alpha},\cdot),
 \qquad |z_\eps|\le CU_j,\qquad \|z_\eps\|_{\mathcal D^{1,2}}\le C.
\]
Passing to a weak limit $z\in\mathcal X^{(j)}$ with $|z|\le CU_j$ and pairing with $w$ gives
\[
 \int_\Gamma U_j^{9/8}|w|\,\dd x'
          =Q_{g^{(j)},\widetilde U_j}(w,z)
          =\langle(f,f^{\partial}),z\rangle
          \le C\left(\int_{\R^{18}_+} U_j|f|\,\dd x+\int_\Gamma U_j|f^{\partial}|\,\dd x'\right),
\]
and testing the equation for $w$ with $Z_{j\alpha}$ bounds $|\eta|$.

Since $g_N=g^{(j)}$ on $\supp Z_{j\alpha}$, we have $w_N\in\mathcal X^{(j)}$ and
\[
 \mathcal L_jw_N=(f_j^*,f_j^{*,\partial})+\sum_\alpha\eta_{j\alpha}Q_{g^{(j)}}(Z_{j\alpha},\cdot),
\]
where
\begin{align*}
 f_j^*&=-E_N+(L_{g_N}-L_{g^{(j)}})w_N
                  -\sum_{i\ne j,\alpha}\eta_{i\alpha}L_{g_N}Z_{i\alpha},\\
 f_j^{*,\partial}&=-E_N^\partial+N_{\widetilde W_N}(w_N)
       +18\bigl(\widetilde W_N^{1/8}-\widetilde U_j^{1/8}\bigr)w_N
                  +\sum_{i\ne j,\alpha}\eta_{i\alpha}B_{g_N}Z_{i\alpha}.
\end{align*}
Set $M_j=\int_\Gamma U_j^{9/8}|w_N|\,\dd x'$. Using
\[
 \int_\Gamma U_j|N_{\widetilde W_N}(w_N)|\,\dd x'
 \le C\eps_M\int_\Gamma U_j\widetilde W_N^{1/8}|w_N|\,\dd x'
 \le C\eps_MM_j+C\eps_M^2e_j,
\]
together with \eqref{eq:bridge-moments}, \eqref{eq:cross-metric-moments} and the
single-bubble adjoint estimate gives
\[
 M_j+|\eta_j|
 \le Ce_j+C\eps_MM_j+C\sum_{i\ne j}e_{ij}|\eta_i|.
\]
Since $\sum_{i\ne j}e_{ij}e_i\le(\sup_ie_i)e_j$, setting
$T_N=\max_{1\le j\le N}(M_j+|\eta_j|)/e_j$ yields
\[
 T_N\le C+C\bigl(\eps_M+\sup_ie_i\bigr)T_N.
\]
Choosing $M$ so that $C(\eps_M+\sup_ie_i)<1/2$ bounds $T_N$
uniformly and proves \eqref{eq:local-trace-moment}.
\end{proof}

\par\medskip
\noindent\textbf{Acknowledgements.}
The author thanks Professor Jingang Xiong for his continued advice
and support in choosing the research topic and writing this paper,
and for his guidance in the author's studies more generally.

\par\medskip
\noindent\textbf{AI usage statement.}
The author selected the overall approach and used GPT-6.1 Sol to carry
out extensive calculations and work out the detailed counterexample
constructions. The model also assisted with language editing.
The author reviewed, verified, and substantially revised all adopted
material and takes full responsibility for the correctness and final
content of this paper.


\begin{thebibliography}{99}
\small
\setlength{\itemsep}{2pt}
\setlength{\parsep}{0pt}

\bibitem{ADN}
S. Agmon, A. Douglis, and L. Nirenberg,
Estimates near the boundary for solutions of elliptic partial
differential equations satisfying general boundary conditions. II,
\emph{Comm. Pure Appl. Math.} \textbf{17} (1964), 35--92.

\bibitem{Almaraz}
S. de Moura Almaraz,
An existence theorem of conformal scalar-flat metrics on manifolds
with boundary,
\emph{Pacific J. Math.} \textbf{248} (2010), 1--22.

\bibitem{AlmarazBlowup}
S. de Moura Almaraz,
Blow-up phenomena for scalar-flat metrics on manifolds with boundary,
\emph{J. Differential Equations} \textbf{251} (2011), 1813--1840.

\bibitem{AAR}
G. E. Andrews, R. Askey, and R. Roy,
\emph{Special Functions},
Cambridge University Press, Cambridge, 1999.

\bibitem{Brendle}
S. Brendle,
Blow-up phenomena for the Yamabe equation,
\emph{J. Amer. Math. Soc.} \textbf{21} (2008), 951--979.

\bibitem{BrendleMarques}
S. Brendle and F. C. Marques,
Blow-up phenomena for the Yamabe equation II,
\emph{J. Differential Geom.} \textbf{81} (2009), 225--250.

\bibitem{CJSX}
L. Caffarelli, T. Jin, Y. Sire, and J. Xiong,
Local analysis of solutions of fractional semi-linear elliptic
equations with isolated singularities,
\emph{Arch. Ration. Mech. Anal.} \textbf{213} (2014), 245--268.

\bibitem{Escobar}
J. F. Escobar,
Conformal deformation of a Riemannian metric to a scalar flat
metric with constant mean curvature on the boundary,
\emph{Ann. of Math.} (2) \textbf{136} (1992), 1--50.

\bibitem{GW}
M. Gr\"uter and K.-O. Widman,
The Green function for uniformly elliptic equations,
\emph{Manuscripta Math.} \textbf{37} (1982), 303--342.

\bibitem{HJWXZ1}
Z.-C. Han, Q. Jiang, H.-Y. Wang, J. Xiong, and L. Zhang,
Isolated singularities of solutions to the Yamabe equation in
higher dimensions I,
preprint (2026), \href{https://arxiv.org/abs/2609.28323}{arXiv:2609.28323}.

\bibitem{HJWXZ2}
Z.-C. Han, Q. Jiang, H.-Y. Wang, J. Xiong, and L. Zhang,
Isolated singularities of solutions to the Yamabe equation in
higher dimensions II,
preprint (2026), \href{https://arxiv.org/abs/2610.24703}{arXiv:2610.24703}.

\bibitem{JLXBoundary}
T. Jin, Y. Liao, and J. Xiong,
Classification and asymptotics for the boundary Yamabe equation
with an isolated singularity,
preprint (2026), \href{https://arxiv.org/abs/2611.07636}{arXiv:2611.07636}.

\bibitem{CompanionI}
Y. Liao,
Isolated singularities of solutions to the Yamabe equation with boundary
in higher dimensions I,
preprint (2026), \href{https://arxiv.org/abs/2611.20647}{arXiv:2611.20647}.

\bibitem{LM}
Y. Liao and Y. Ma,
Isolated singularities of solutions to the Yamabe equation with
boundary in dimension 3 and 4,
preprint (2026), \href{https://arxiv.org/abs/2607.16650}{arXiv:2607.16650}.

\bibitem{Lieberman}
G. M. Lieberman,
Higher regularity for nonlinear oblique derivative problems in Lipschitz domains,
\emph{Ann. Scuola Norm. Sup. Pisa Cl. Sci.} (5)
\textbf{1} (2002), 111--151.

\bibitem{LSW}
W. Littman, G. Stampacchia, and H. F. Weinberger,
Regular points for elliptic equations with discontinuous coefficients,
\emph{Ann. Scuola Norm. Sup. Pisa Cl. Sci.} (3)
\textbf{17} (1963), 43--77.

\bibitem{Szego}
G. Szeg\H{o},
\emph{Orthogonal Polynomials}, 4th ed.,
American Mathematical Society, Providence, RI, 1975.

\end{thebibliography}
\end{document}